\documentclass[a4paper,11pt]{article}
\small 
\normalsize

\usepackage[T1]{fontenc}
\usepackage[english]{babel}
\usepackage{times}
\usepackage{ulem}

\usepackage[all]{xy}

\usepackage{stmaryrd}
\usepackage{graphicx}
\usepackage{paralist}
\usepackage{amsfonts}
\usepackage[leqno]{amsmath}

\usepackage{amssymb,mathrsfs}
\usepackage{color}
\usepackage{enumerate,vmargin}
\usepackage{tocvsec2}
\usepackage{verbatim}
\usepackage[labelfont=bf,font=it]{caption}
\usepackage[numbers,sort&compress]{natbib}

\setmarginsrb{2.9cm}{2.6cm}{2.9cm}{1.7cm}{0cm}{0mm}{0cm}{9mm}

\usepackage{hyperref}
\hypersetup{
    unicode      = false,     
    pdftoolbar   = true,      
    pdfmenubar   = true,      
    pdffitwindow = true,      
    pdfnewwindow = true,      
    colorlinks   = true,      
    linkcolor    = blue,      
    citecolor    = red,      
    filecolor    = blue,      
    urlcolor     = green       
}

\definecolor{vio}{rgb}{.5,.1,.5}
\definecolor{grey}{rgb}{.5,.5,.5}
\definecolor{zhanc}{rgb}{.9,0.5,0}
\definecolor{bernc}{rgb}{.6,0,.8}
\definecolor{thomc}{rgb}{.6,0,.8}
\definecolor{darkred}{RGB}{165,42,42}

\def \tt#1{{\color{darkred}#1\color{black}}}

\def\thesection{\arabic{section}}
\renewcommand{\theequation}{\thesection.\arabic{equation}}

\newtheorem{theorem}{Theorem}[section]
\newtheorem{proposition}[theorem]{Proposition}
\newtheorem{conjecture}[theorem]{Conjecture}
\newtheorem{definition}[theorem]{Definition}
\newtheorem{example}[theorem]{Example}
\newtheorem{remark}[theorem]{Remark}

\newtheorem{lemma}[theorem]{Lemma}

\newcommand{\noi}{\noindent}

\newcommand{\un}{{\bf 1}}

\newcommand{\p}[1]{{\mathbb P}\left(#1\right)}

\newcommand{\e}{{\mathbb E}}

\newcommand{\bbN}{\mathbb N}

\newcommand{\bbR}{\mathbb R}

\def\cq{$\hfill \square$}
\def\cqfd{$\hfill \blacksquare$}
\def\ino{ \! \in \! }

\def\r{{\mathbb R}}
\def\e{{\mathbb E}}
\def\p{{\mathbb P}}

\def\P{{\bf P}}

\def\ee{\mathrm{e}}
\def\d{\, \mathrm{d}}

\def\leqo{\! \leq \! }
\def\leko{\! < \! }
\def\geqo{\! \geq \! }
\def\geko{\! > \! }
\def\eqo{\! = \! }

\def\epp{\varepsilon}

\def\ff{{\mathbf{F}}}

\def\fg{{\mathrm{F}}}
\def\sgn{\varepsilon}
\def\inv{\mathrm{inv}}
\def\psib{\overline{\psi}}
\def\phib{\overline{\phi}}

\newcommand{\eqnsection}{
\renewcommand{\theequation}{\arabic{section}.\arabic{equation}}
    \makeatletter
    \csname  @addtoreset\endcsname{equation}{section}
    \makeatother}
\eqnsection

\usepackage[refpage,noprefix]{nomencl}                          
\makeglossary
\makenomenclature                                           

\graphicspath{{.}{figures/}}

\title{ \textbf{Hipster random walks, random series-parallel graphs and random homogeneous systems}}

\author{Xinxing \textsc{Chen} 
\thanks{School of Mathematical Sciences, Shanghai Jiaotong University, 200240 Shanghai, China. Email: chenxinx@sjtu.edu.cn Partially supported by NSFC grant No.~12271351.}
\and Thomas \textsc{Duquesne}
\thanks{LPSM, Sorbonne Universit\'e, Campus Pierre et Marie Curie, 4 place Jussieu, F-75252 Paris Cedex 05, France.
Email: thomas.duquesne@sorbonne-universite.fr}
\and Zhan \textsc{Shi}
\thanks{State Key Laboratory of Mathematical Sciences, AMSS, Chinese Academy of Sciences, 
100190 Beijing, China.
Email: shizhan@amss.ac.cn}
}

\begin{document}

\maketitle

\begin{abstract} 
We study a class of random processes that we call \emph{random homogeneous systems}. They extend   
distance and resistance of \emph{critical random series-parallel graph} introduced by Hambly \& Jordan~\cite{hambly-jordan} 
and Pemantle \cite{pemantle}, and \emph{Hipster random walks} introduced by Addario-Berry, Cairns, Devroye, Kerriou \& Mitchell~\cite{addario-berry-cairns-devroye-kerriou-mitchell}.  
Our main result asserts that under suitable general assumptions, these systems converge weakly, upon an appropriate normalization, to the probability distribution with density $\frac34 \, (1-x^2) \, {\bf 1}_{\{ x\in (-1, \, 1)\} }$. 
It gives an affirmative answer to a conjecture of Hambly and Jordan~\cite{hambly-jordan} and further conjectures of Addario-Berry et al.~\cite{addario-berry-cairns-devroye-kerriou-mitchell} and Derrida~\cite{bernard}, whereas for the hipster random walk, we recover a previous result of Addario-Berry, Cairns, Devroye, Kerriou and Mitchell \cite{addario-berry-cairns-devroye-kerriou-mitchell}.

\medskip


\medskip


\end{abstract}

\section{Introduction}
Below we  
use the notations $\r_+^* \! := \! (0, \, \infty)$, $\r_+ \! :=\!  [0, \infty)$ and $\r_-\! :=\!  (-\infty, 0]$; $\bbN$ stands for the set of nonnegative integers and $\bbN^*\! := \! \bbN\backslash \{ 0\}$ for the set of positive integers. 
Let us also mention that all the random variables (r.v.s for short) that are mentionned in this article are defined on the same probability space $(\Omega, \mathscr F , \p)$.

Let $(X_n)_{  n\in \bbN }$ be positive r.v.s that are defined recursively as follows. Let $X_0$ be an arbitrary 
$\r_+^*$-valued r.v. For all $n\ino \bbN$, let $X'_n$ denote an independent copy of $X_n$ and let $\ff \!  :\!   (\r_+^*)^2 \! \to \! \r_+^*$ be a random function which is independent of $(X_n, X'_n)$ and whose law is fixed. Then 
\begin{equation}
    X_{n+1} \; {\buildrel \mathrm{law} \over =} \; \ff (X_n, \, X'_n) .
    \label{syshom}
\end{equation}
\noindent 
Here $\overset{\textrm{law}}{=}$ means that the laws of the r.v.s are equal. 
The letter $\ff$ is kept in boldface to remind its randomness.
The law of $X_n$ is thus characterized by $n$, by the law of $X_0$ and by the law of $\ff$. 
The goal of the paper is to study the asymptotic behaviour of the law of $X_n$  as $n$ goes to $\infty$ under specific  assumptions on $\ff$ that are discussed below. 
Before getting into more details let us provide a brief overview of the content of this introduction.

\smallskip

\noindent 
$-$ First, we specify the type of functions $\ff $ that we are considering. Namely $1$-homogeneous 
continuous functions that are nondecreasing in each of their coordinates and either greater than or less than the $\max$ or $\min$ functions (see Definition \ref{classH} below).

\smallskip

\noindent 
$-$ This general framework unifies two previously studied models: 
\begin{compactenum}
\item[$(i)$] {\it Distances} and {\it resistances} on {\it series-parallel graphs}, as introduced by Hambly and Jordan~\cite{hambly-jordan} (Auffinger and Cable \cite{auffinger-cable} 
solve the distances in the critical case and C., Derrida, D.~and S.~\cite{serpar_distance} study the distances in the slightly super-critical case);  
\item[$(ii)$] {\it Hipster random walks}, as introduced and studied by Addario-Berry, Cairns, Devroye, Kerriou 
and Mitchell \cite{addario-berry-cairns-devroye-kerriou-mitchell}. 
\end{compactenum}
We review these models and related results. 

\smallskip

\noindent 
$-$ Next, we discuss Derrida's conjectures on the possible behaviours of $X_n$ within our framework. Roughly speaking  these conjectures state that $\log X_n$ is either bounded, of order $n$, $n^{1/2}$ or $n^{1/3}$.

\smallskip

\noindent
$-$ Finally, we present our main result, Theorem \ref{t:main}, which proves Derrida's conjecture when $\log X_n$ is of order $n^{1/3}$. 
Theorem \ref{t:main} entails 
the conjecture on resistances of series-parallel graphs (due to Addario-Berry, Cairns, Devroye, Kerriou 
and Mitchell \cite{addario-berry-cairns-devroye-kerriou-mitchell}) and extends the result of the same authors on symmetric hipster random walks.

\smallskip

\noindent 
$-$ It is important to mention that an article by Morfe \cite{morfe}, pre-published on Arxiv a few days before ours, addresses the same model and solves Derrida's conjectures using PDE methods that are completely different from ours. 
See the comments at the end of the introduction.

\medskip

\noindent
{\large \textbf{The class of random functions $\mathbf{F}$ that are considered}}. $\; $
Equation \eqref{syshom} defines a random iterative system, and can also 
be viewed as a random system defined on a hierarchical lattice. Random iterative 
systems were well studied in probability theory, see for example 
Shneiberg \cite{shneiberg}, Li and Rogers \cite{li-rogers}, Wehr and Woo \cite{wehr-woo}, Jordan \cite{jordan}, Goldstein \cite{goldstein}. 
On the other hand, hierarchical lattices have an important literature in theoretical physics 
(Berker and Ostlund \cite{berker-ostlund}, Clark \cite{clark}, Griffiths and Kaufman \cite{griffiths-kaufman}, Derrida, De Seze and Itzykson \cite{derrida-deseze-itzykson}, Derrida and Griffiths \cite{derrida-griffiths} Tremblay and Southern \cite{tremblay-southern}), 
and the hierarchical structure leads to random recursions in the form of \eqref{syshom}. 
The recursions are often simple to formulate, but hard to analyze, leaving many open problems, 
such as the innocent-looking Derrida--Retaux model \cite{derrida-retaux}. The seminal paper Aldous and Bandyopadhyay \cite{aldous-bandyopadhyay} 
gives a general account of fixed points of recursive distributional equations and contains several exactly solvable families of systems with a max-type random function $\ff$.

As already mentioned our goal is to investigate the asymptotic behavior of $X_n$ as in (\ref{syshom}) when $\mathbf F$ belongs to the specific class of functions $\mathscr H$ that is defined as follows. 
\begin{definition}
\label{classH}
{\rm 
 Let $\mathscr H$ be the space of continuous functions 
$F \!: \!  (\r_+^*)^2 \!  \to \! \r_+^*$ satisfying the following properties. 
\begin{compactenum}

\smallskip

\item[$(a)$] $F$ is nondecreasing in each coordinate and $1$-homogeneous. Namely for all $x,y,z\ino \r^*_+$, 
$\min \big( F(x+z,y), F(x,y+z) \big)\geqo F(x,y) $ and $F(zx,zy)\eqo zF(x,y)$.

\smallskip

\noindent
\item[$(b)$] Either $F(x,y) \geqo 
\max (x,y) $ for all $x,y\ino \r^*_+$ and we set $\varepsilon (F)\eqo 1$, 
or $F(x,y) \leqo \min (x,y) $ for all $x,y\ino \r^*_+$ and we set $\varepsilon (F)\eqo  -1$.

\smallskip

\noindent
\item[$(c)$] For all $x\ino \r_+^*$, if $\varepsilon (F)\eqo 1$, $\lim_{y\to 0} F(x,y)\eqo \lim_{y\to 0} F(y,x)\eqo x$ 
and if $\varepsilon (F)\eqo -1$, $\lim_{y\to \infty} F(x,y)\eqo  \lim_{y\to \infty} F(y,x)\eqo x$. 

\smallskip

\end{compactenum}
Note that $F$ can be extended continuously on $\r_+^2$ by setting $F(0,0)\eqo 0$ and for all $x\ino \r^*_+$,  
\begin{equation}
\label{contiexthomog}
 F(x,0)\eqo F(0,x) \eqo  \lim_{y\to 0} F(x,y)\eqo \lim_{y\to 0} F(y,x)\eqo x\un_{\{ \epp (F)= 1\}}\; .
\end{equation}
\noindent
We also set  $\mathscr H_+ \! =\! \big\{ F \! \in\!  \mathscr H \! : \varepsilon (F)\! =\!  1\big\}$ and 
$\mathscr H_- \! = \! \big\{ F \! \in \! \mathscr H : \varepsilon (F) \! =\!  -1\big\}$. \cq 
 } 
\end{definition}

Note that $\mathscr H$ is a Borel subset of $C^0(\bbR_+^2,  \bbR)$ equipped with the Polish topology of uniform convergence on every compact subset of $\bbR_+^2$. By an $\mathscr H$-valued random function $\ff$, we mean an $\mathscr F$-Borel measurable r.v.~$\ff\! : \Omega \to C^0(\bbR_+^2,  \bbR)$ such that $\p (\ff \ino \mathscr H)\eqo 1$.

\begin{example}
\label{1homogex} 
{\rm 
$(a)$ It is convenient to set 
$\varphi_+(x,y)\! =\!  \max (x,y)$  and $ \varphi_-(x,y) \! =\!  \min (x,y)$ for all $x,y\ino \r_+$. Note that $\varphi_+ \ino \mathscr H_+$ and that $\varphi_- \ino \mathscr H_-$. 

\smallskip

\noindent
$(b)$ Let $\alpha \ino \r\backslash \{ 0\}$. For all $x,y\ino \r_+^*$ we set $F_{\! \alpha} (x,y)\eqo (x^{1/\alpha}+y^{1/\alpha})^{\alpha}$. Then $ F_{\! \alpha} \ino \mathscr H_+$ if $\alpha \geko 0$ and $ F_{\! \alpha} \ino \mathscr H_-$ if $\alpha \leko 0$. \cq 
 }
\end{example}
In order to connect our framework with hipster random walks (see below) it is convenient to provide a representation in the logarithmic scale of functions of $\mathscr H$ thanks to the following space. 
\begin{definition}
\label{classG} 
{\rm 
 Let $\mathscr G$ be the space of continuous functions $G \!: \! \r \mapsto \r_+$ satisfying the following properties.
\begin{compactenum}

\smallskip

\item[$(a)$] On $\r_+$, $G$ is $1$-Lipschitz nonincreasing and on $\r_-$, $G$ is $1$-Lipschitz nondecreasing.  

\smallskip

\item[$(b)$]  $\lim_{z\to -\infty} G(z)\eqo \lim_{z\to \infty} G(z)\eqo 0$. \cq 
\end{compactenum}
} 
\end{definition}
Then there is a one-to-one correspondence between $F\ino \mathscr H$ and 
$(\varepsilon, G)\ino \{ -1, 1\} \times \mathscr G$ given by 
\begin{equation}
\label{correspond}
\varepsilon (F)= \varepsilon \quad \textrm{and} \quad \log F(\mathrm{e}^x, \mathrm{e}^y) = \varphi_\varepsilon (x,y)+ \varepsilon G(\varepsilon (x-y)) , \quad x,y\in \r, 
\end{equation}
where we recall from Example \ref{1homogex} $(a)$ that $\varphi_+ (x,y)\eqo \max (x,y)$ and $\varphi_- (x,y)\eqo \min (x,y)$, and where $\varphi_\epp$ means $\varphi_+$ if $\epp \eqo 1$ and $\varphi_-$ if $\epp=-1$.  
See Proposition \ref{corrprop} for a more detailed statement. 
\begin{remark}
\label{GLipex} 
 {\rm 
$(a)$ We observe that $F\ino \{ \varphi_+, \varphi_-\}$ if and only if $G$ as in (\ref{correspond}) 
is identically null. 

\smallskip
 \noindent
$(b)$ Let $\alpha \ino \r\backslash \{ 0\}$, let $F_{\! \alpha}$ be as in Example \ref{1homogex} and let $G_{\! \alpha} \ino \mathscr G$ be associated with $F_{\! \alpha}$ as in (\ref{correspond}). Then $G_{\! \alpha} (z) \eqo |\alpha| \log (\mathrm e^{ z/|\alpha|} +1) \! -\! (z)_+$, where $(z)_+\eqo \max (z, 0)$ is the positive part of $z$. \cq 
 } 
\end{remark}

\noi
\textbf{Notation.} Unless otherwise specified, in this article $\mathbf F$ is an $\mathscr H$-valued random function whose law is fixed  
and in what follows, we denote by $\mathbf G$ its $\mathscr G$-valued corresponding function as defined in (\ref{correspond}). 

\bigskip

\noindent
{\large \textbf{Previous models.}} One of the motivations for introducing the general framework of random $\mathscr H$-valued functions is to unify several models that have been previously investigated. Namely the series-parallel graphs and the hipster random walks that we discuss now. 
Understanding the first model was, in fact, the original motivation of this work.

\begin{example}
\label{ex:graph}
{\rm 
{\it(The random series-parallel model).} We fix a parameter $p\in [0, \, 1]$ and we consider a sequence $(\mathtt{Graph}(n))_{n\in \bbN }$ of random graphs defined recursively as follows: $\mathtt{Graph}(0)$ is simply composed of two vertices 
${\mathtt a}$ and ${\mathtt z}$ connected by a (non-oriented) edge; for each $n\ino \bbN$, $\mathtt{Graph}(n+1)$ is obtained from $\mathtt{Graph}(n)$ by replacing each edge of $\mathtt{Graph}(n)$ either with probability $p$ by two edges in series, or with probability $(1\! -\! p)$ by two parallel edges (all the replacements being independent, and independent of everything in $\mathtt{Graph}(n)$). The resulting sequence of random graphs is called the series-parallel random graphs with parameter $p$.
We assign each edge of $\mathtt{Graph}(n)$ with a unit length  and a unit resistance and we denote by $\Delta_n$ and $R_n$ respectively the distance and the effective resistance between vertices ${\mathtt a}$ and ${\mathtt z}$ on $\mathtt{Graph}(n)$.
Here $\Delta_0\eqo R_0\eqo 1$. 

By considering the possible two situations for $\mathtt{Graph}(1)$, i.e.~two edges in series or two edges in parallel, we check 
that for all $n\ino \bbN$,
\begin{equation}
   \Delta_{n+1}{\buildrel \hbox{\scriptsize law} \over =}  (\Delta_n\! + \Delta'_n)  \mathscr{E} + \min (\Delta_n ,\Delta'_n) (1\! -\! \mathscr{E}) \; \,  \textrm{and} \; \,    R_{n+1}
     {\buildrel \hbox{\scriptsize law} \over =} 
    (R_n\!  + R'_n) \mathscr{E} + \frac{1\! -\! \mathscr E}{\frac{1}{R_n}  + \frac{1}{R'_n}} ,
    \label{def_recurrence_eff_conduct}
\end{equation}
\noindent where $\Delta_n, \Delta'_n, \mathscr E$ (resp.~$R_n, R'_n, \mathscr E$) are independent,  
$\Delta_n$ and $\Delta'_n$  (resp.~$R_n$ and $R'_n$) have the same law and $\mathscr{E}$ is a Bernoulli$(p)$ r.v.: $\p(\mathscr{E}\eqo1) \eqo p \eqo 1\! -\!  \p(\mathscr{E} \eqo 0)$. 
Using the notations of Example \ref{1homogex} we easily check the following. 
\begin{compactenum}

\smallskip

\item[$-$] The \emph{$p$-distance model} corresponds to a random $\mathscr H$-valued function 
$\mathbf F$ whose law is given by $ \p (\mathbf F\! =\!  F_1) \! =\!  p$ and $\p (\mathbf F \! =\!  \varphi_-) \! =\!  1\! -\! p$.

\smallskip

\item[$-$] The \emph{$p$-resistance model} corresponds to a random $\mathscr H$-valued function $\mathbf F$ whose law is given by $ \p (\mathbf F \! =\!  F_1) \! =\!  p$ and $\p (\mathbf F \! =\!  F_{-1}) \! =\!  1\! -\! p$. 
\end{compactenum}

\smallskip

It has been proved by Hambly and Jordan~\cite{hambly-jordan}, among many other results, that $\Delta_n$ and 
$R_n$ experience a phase transition at $p\eqo p_c \! :=\!  \frac12$. More precisely, as $n\! \to \! \infty$, if $p\geqo p_c$, then 
$\Delta_n \! \to \! \infty$ and  $R_n\!  \to \!  \infty$ weakly, and if $p\leko  p_c$, then $\Delta_n $ tends to a finite r.v.~$\Delta_\infty$  and $R_n \! \to \! 0$ weakly. 
For $p\eqo p_c $ and for the distance model Hambly and Jordan~\cite{hambly-jordan} have conjectured a subexponential growth of $\Delta_n$ and Auffinger and Cable \cite{auffinger-cable} have proved the following. 
\begin{theorem}
\label{t:distance}
{\bf (Auffinger and Cable.~\cite{auffinger-cable})} Let $(\Delta_n)_{n\in \bbN }$ be as in \eqref{def_recurrence_eff_conduct} with $p\eqo p_c \eqo \frac12$. Then for all $y\in \r$,
\begin{equation}
    \lim_{n\to \infty} \p \Big( \frac{\Delta_n}{(c_\mathtt{dis}\, n)^{1/2}} \le y \Big)
    =
    \int_{-\infty}^y 2x\, {\bf 1}_{\{ 0< x <1 \} } \, \d x \, ,
    \label{wcv_distancer}
\end{equation}
\noindent where $c_\mathtt{dis} :=\pi^2/3$.
\end{theorem}

 C., Derrida, D.~and S.~in \cite{serpar_distance} study distances in the slightly super-critical cases. Namely, for all $p\geko p_c \eqo \frac12$ a simple subadditive argument shows that $\frac{1}{n}\log  \e [\Delta_n]\!  \to \! \alpha (p) \ino (0, \infty)$. Then it is proved in \cite{serpar_distance} that 
$$ \lim_{p\to p_c} \frac{\alpha (p)}{\sqrt{p-p_c}}= \frac{\pi}{\sqrt{6}} \; . $$ 
  Concerning the resistance model Hambly and Jordan made the following conjecture.
\begin{conjecture}
\label{conj:hambly-jordan}
{\bf (Hambly and Jordan~\cite{hambly-jordan}).} If $p\eqo p_c\eqo \frac12$, then as $n\! \to \! \infty$,
$\lim_{n\to \infty}|\log R_n| \eqo \infty$ weakly. 
\end{conjecture}

Let us assume $p\eqo p_c\eqo \frac12$. 
Then $R_n$ and $1/R_n$ have the same distribution and Conjecture \ref{conj:hambly-jordan} 
is thus equivalent to saying that $\lim_{n\to \infty} \p(R_n \ino [x, \, y]) \eqo 0$ for all $0\leko x\leqo y\leko \infty$. 
Furthermore, it is briefly mentioned in \cite{hambly-jordan} that $R_n$ would be asymptotically exponentially 
large or exponentially small in this case, but later on, Addario-Berry et 
al.~\cite{addario-berry-cairns-devroye-kerriou-mitchell} argued that they would rather expect $\frac{\log R_n}{n^{1/3}}$ to converge weakly. More precisely, here is their prediction. 
\begin{conjecture}
\label{conj:addario-berry}
{\bf (Addario-Berry et al.~\cite{addario-berry-cairns-devroye-kerriou-mitchell}).} If $p\eqo p_c\eqo \frac12$, there is $c_\mathtt{eff}\ino (0, \, \infty)$ such that for all $y\ino \r$,
\begin{equation}
    \lim_{n\to \infty} \p \Big( \frac{\log R_n}{(c_\mathtt{eff}\, n)^{1/3}} \le y \Big)
    =
    \int_{-\infty}^y \!\! \tfrac{3}{4}  \big(1\! -\! x^2 \big) \, {\bf 1}_{\{ |x| <1 \} } \, \d x \, .
    \label{wcv_resistance}
\end{equation}
\end{conjecture}

\noindent
It is clear that Conjecture \ref{conj:addario-berry} implies Conjecture \ref{conj:hambly-jordan}.
\begin{conjecture}
\label{conj:zeta(3)}
{\bf (Derrida~\cite{bernard}).}
In Conjecture \ref{conj:addario-berry}, $c_\mathtt{eff}\eqo  9 \zeta(3)$ where $\zeta(3)\eqo \sum_{n=1}^\infty n^{-3}$.
\end{conjecture}

\noindent
Our main result yields Conjectures \ref{conj:addario-berry} and \ref{conj:zeta(3)} as a special case. 
Moreover we show that the probability measure $\frac34 \, (1-x^2) \, {\bf 1}_{\{ |x| <1 \} } \, \d x$ 
is actually the limiting distribution of a whole family of random homogeneous systems.\hfill$\Box$
} 

\end{example}
\begin{example}
\label{ex:hipster}
{\rm 
{\it (Hipster random walks).} Let $V_0$ and $U_0$ be integer-valued r.v.s. For $n\ino \bbN$ we set 
 \begin{equation}
      U_{n+1}
     \! = \! 
     \xi_n U_n+(1\! -\! \xi_n)U'_n + D_n \, {\bf 1}_{\{ U_n=U'_n\}} \; \textrm{and} \; V_{n+1}
     =
     \xi_n V_n+(1\! -\! \xi_n)V'_n + C_n \, {\bf 1}_{\{ V_n=V'_n\}} \, , \label{hipster}
 \end{equation}
where $\xi_n$, $U_n$, $U'_n$ and $D_n$ (resp.~$\xi_n$, $V_n$, $V'_n$ and $C_n$)
are independent r.v.s with $\P(\xi_n\eqo 0) \eqo \P(\xi_n\eqo 1) \eqo \frac12$, $\P(D_n\eqo -\! 1) \eqo  \P(D_n\eqo 1) \eqo \frac12$,  $\P(C_n\eqo 1) \eqo  \P(C_n\eqo 0) \eqo  \frac12$ and $U'_n$ and $U_n$ (resp.~$V'_n$ and $V_n$) have the same law.  The process $(U_n)_{n\in \bbN}$ (resp.~$(V_n)_{n\in \bbN}$) was introduced and called the {\it symmetric hipster random walk} (resp.~the {\it lazy hispter random walk}) 
 by Addario-Berry, Cairns, Devroye, Kerriou 
and Mitchell \cite{addario-berry-cairns-devroye-kerriou-mitchell} who proved the following. 
\begin{theorem}
\label{t:hipster}
{\bf (Addario-Berry et al.~\cite{addario-berry-cairns-devroye-kerriou-mitchell}).} Let $V_n $ and $U_n$, $n\ino \bbN$, be as in \eqref{hipster}. Then for all $y\ino \r$,
\begin{equation}
    \lim_{n\to \infty} \p \Big( \frac{V_n}{(c_\mathtt{laz}\, n)^{1/2}} \le y \Big)
    =
    \int_{-\infty}^y \!\! 2x {\bf 1}_{\{ 0< x <1 \} } \, \d x \, \quad \textrm{and}    \label{wcv_lazyhipster}
\end{equation}
\begin{equation}
    \lim_{n\to \infty} \p \Big( \frac{U_n}{(c_\mathtt{hip}\, n)^{1/3}} \le y \Big)
    =
    \int_{-\infty}^y \!\! \tfrac{3}{4}  \big( 1\! -\! x^2 \big)  {\bf 1}_{\{ |x| <1 \} } \, \d x \, ,
    \label{wcv_hipster}
\end{equation}
\noindent where $c_\mathtt{laz} \eqo 2$ and  $c_\mathtt{hip}\eqo \frac92$.
\end{theorem}

  Similarities between \eqref{wcv_resistance} and \eqref{wcv_hipster} are immediately observed. In fact, Addario-Berry, Cairns, Devroye, Kerriou 
and Mitchell \cite{addario-berry-cairns-devroye-kerriou-mitchell} argued, in a heuristic manner, that the effective resistance of the critical random series-parallel graph and the symmetric hipster random walk should have similar asymptotics; this had led from Theorem \ref{t:hipster} to Conjecture \ref{conj:addario-berry} in \cite{addario-berry-cairns-devroye-kerriou-mitchell}, though the value of the constant $c_\mathtt{eff}$ in \eqref{wcv_resistance} was not clear. See also Addario-Berry, Beckman and Lin  \cite{addario-berry-beckman-lin} for related models.

Our main result also entails Theorem \ref{t:hipster}. We have not been able to justify the heuristics in \cite{addario-berry-cairns-devroye-kerriou-mitchell} arguing that the effective resistance of the critical random series-parallel graph in Example \ref{ex:graph} and the hipster random walk in Example \ref{ex:hipster} should have similar asymptotics; rather, we investigate a family of random homogeneous systems that includes both $(U_n)_{ n\in \bbN}$ in Example \ref{t:hipster} and $(R_n)_{n\in \bbN}$ in Example \ref{ex:graph}, and we prove that these systems converge weakly, upon a suitable normalization, to the same distribution $\frac34 \, (1-x^2) \, {\bf 1}_{\{ |x| <1 \} } \, \d x$. As such, we have not established the possible relation between $(U_n)_{ n\in \bbN}$  and $(R_n)_{ n\in \bbN}$ as suggested in \cite{addario-berry-cairns-devroye-kerriou-mitchell}, but we have proved that they belong to a general family of random systems to which our main result applies.

{\it Indeed}, let us briefy explain how these models fit into our framework. We denote by $G_{\mathtt{hip}} (z) = \max (0, 1\! -\! |z|) $, $z\in \r$, which belongs to $\mathscr G$. Here $X_n$ is as in (\ref{syshom}) and  $\mathbf F$ and 
$(\varepsilon (\mathbf F), \mathbf G)$ are as in (\ref{correspond}). It is not difficult to check the following. 

\begin{compactenum}

\smallskip

\item[$(i)$] $\log X_n\overset{\textrm{(law)}}{=} V_n$ with 
$\p (\mathbf G \eqo  G_{\mathtt{hip}}  \, ; \,   \varepsilon(\mathbf F) \eqo 1)\eqo  
\p (\mathbf G\!  \equiv\!  0 \,  ; \,  \varepsilon(\mathbf F) \! =\!  -\! 1)\eqo \frac12$.   

\smallskip

\item[$(ii)$] $\log X_n\overset{\textrm{(law)}}{=} U_n$ with 
$\p (\mathbf G \eqo  G_{\mathtt{hip}}  \, ; \,   \varepsilon(\mathbf F) \eqo 1)\eqo  \p (\mathbf G \eqo  G_{\mathtt{hip}}  \, ; \,   \varepsilon(\mathbf F) \eqo -\! 1)\eqo \frac12$. \cq 
\end{compactenum}
 } 
\end{example}

\medskip

\noindent
{\large \textbf{Possible orders of magnitude of $X_n$}}. $\;$ Let us first mention that most of the statements presented below are conjectural. 
Since we take the infimum or the supremum of a very large number of random variables to obtain $X_n$, the tails of $X_0$ at $0^+$ and $\infty$ and the tail of $\mathbf G(0)$ at $\infty$ actually play a quite complicated role, which 
may blur the neat features of the model that we want to present here. 
To avoid these unnecessary complications, we make the following bounded support assumptions. 
\begin{equation}
\label{boundedsupport}
 \exists c\in (1, \infty) : \quad \p \big( \mathbf G (0)  \! \leq \! c) = \p \big( |\log X_0| \! \leq \! c \big)=1 \; .  
\end{equation}
This immediately implies that 
\begin{equation}
 \forall n\in \mathbb N, \qquad  \p \big( |\log X_n|  \! \leq \! c(n+1) \big) = 1 \; . 
\end{equation}
Let us set 
\begin{equation} 
\label{defiparamp}
p= \p (\varepsilon (\mathbf F)  \! =\!  1) \; .
\end{equation}
To avoid trivial cases, we also assume that 
\begin{equation} 
\label{nontrivhyp}
 \p \big( \mathbf F \notin \{ \varphi_+, \varphi_-\} \big)= \p \big( \textrm{$\mathbf G$ is not null}\big) >0\; .
 \end{equation} 
Let us first consider the cases where $\mathbf G$ is trivial on one side. A simple 
percolation argument on the complete binary tree implies that there is a $\r_+$-valued r.v.~$X_\infty$ such that 
$$ X_n \overset{\textrm{(law)}}{\underset{n\to \infty}{-\!\!\!-\!\!\! \longrightarrow}} X_\infty \quad \textrm{if $p \! >\! \tfrac12$ and $\p (\mathbf F\! =\!  \varphi_+ |\varepsilon (\mathbf F)= 1)= 1$ or if 
$p \! <\! \tfrac12 $ and $\p (\mathbf F \! = \!  \varphi_- | \varepsilon (\mathbf F) = -\! 1)= 1$. }$$

We next justify why we require $\mathbf F$ to satisfy Definition \ref{classH} $(c)$. More precisely, we denote by $\mathscr H^o$ the class of functions satisfying Definition \ref{classH} $(a)$ and $(b)$ (only) and we denote by $\mathscr G^o$ the class of functions that satisfy Definition \ref{classG} $(a)$ (only). Then (\ref{correspond}) establishes a one-to-one correspondence between $F\in \mathscr H^o$ and $(\varepsilon, G) \in \{ -1, 1\} \times \mathscr G^o$ and $F$ satisfies Definition \ref{classH} $(c)$ if and only if $G$ satisfies Definition \ref{classG} $(b)$. Let $\mathbf F$ be a random $\mathscr H^o$-valued function and let 
$\mathbf G \in  \mathscr G^o$ be as in (\ref{correspond}). We assume that there is a constant $\eta >0$ such that a.s.
$$  \mathbf G (-\infty)= \lim_{z\to -\infty} \mathbf G (z) >\eta \quad \textrm{and} \quad
 \mathbf G (\infty)= \lim_{z\to \infty} \mathbf G (z) >\eta .$$
Then heuristic arguments and numerical simulations indicate that $\log X_n$ is genuinely of order $n$ for all $p\in [0,1]$ and we observe no transition in the order of magnitude. 

More precisely, we conjecture that if $p\! >\!  \frac12$, there is a deterministic positive speed $v_+$ such that $\frac1n \log X_n \to v_+$, if $p\! <\!  \frac12$,  there is a deterministic negative speed $v_-$ such that $\frac1n \log X_n \to v_-$ and if $p\! = \! \frac12$, there is a random variable $U$ that is uniformly distributed 
on an open interval containing $0$ such that $\frac1n \log X_n \to U$. So there is a transition at $p\! = \! \frac12$ but it does not affect the order of magnitude of $\log X_n$ which is $n$.

From now on we assume that $\mathbf F$ a.s.~belongs to $ \mathscr H$ and we discuss conjectures raised by Derrida and 
based on a discrete non-linear PDE derived from (\ref{syshom}). To simplify our statements we assume that $\mathbf F$ is symmetric, i.e.~that a.s.~$\mathbf F (x,y)\! =\!  \mathbf F(y,x)$ for all $x,y\in \r^*_+$.   
In addition to $p$ as defined in (\ref{defiparamp}) we introduce the following two parameters: 
\begin{equation}
\label{alphapm}
\alpha_+ \eqo   \int_0^\infty \! \!\! \e \big[ \mathbf G(z)  \, \big| \, \varepsilon (\mathbf F)\! = \! 1\big] \mathrm dz \quad \textrm{and} \quad \alpha_- \eqo   \int_0^\infty \! \!\!  \e \big[ \mathbf G(z)  \, \big| \, \varepsilon (\mathbf F)\! = \! -1\big] \mathrm dz \, 
\end{equation}
which are supposed to be finite. To properly prove the following conjectures additional moments are needed but our purpose here only consists in providing a broad guideline for possible convergences. 
\begin{conjecture}
\label{conj4orders} {\bf (Derrida~\cite{bernard})}

\begin{compactenum}

\smallskip

\item[$(a)$] If $p\! > \! \frac12$ and if $\, \p (\mathbf F\! =\!  \varphi_+ \, |\, \varepsilon (\mathbf F) \! = \! 1)\! < \!  1$, 
then there is $c_+\in \r_+^*$ such that $\e [ \log X_n] \! \geq \!  c_+ n $.  
If $p\! < \! \frac12$ and if $\p (\mathbf F\! =\!  \varphi_- \, |\, \varepsilon (\mathbf F) \! = \! -1)\! < \!  1$, 
then there is $c_-\in \r_+^*$ such that $-\e [\log X_n] \! \geq \!  c_- n $. 

\smallskip

\item[$(b)$] If $p \! =\!  \frac12$ and $\alpha_+  \! \neq \!  \alpha_-$, then $ n^{-1/2}\log X_n$ converges in law to $\mathtt{Cst} \sqrt{U}$, where $U$ is uniformly distributed on $[0, 1]$ and where $\mathtt{Cst}$ is a constant which depends on the law of $\mathbf F$. 

\smallskip

\item[$(c)$] If $p \! =\!  \frac12$ and $\alpha_+  \! = \!  \alpha_-$, then $ n^{-1/3}\log X_n$ converges in law to $\mathtt{Cst} V$, where the density of $V$ is $\frac{3}{4} (1\! -\! x^2)\mathbf 1_{[-1, 1]} (x)$ and where $\mathtt{Cst}$ is a constant which depends on the law of $\mathbf F$. 
\end{compactenum}
\end{conjecture}
The parameters $(\alpha_+, \alpha_-, p)$ range in $D\! :=\! \r_+^2 \! \times \! [0, 1] $ when $\mathbf F$ is distributed according to all the possible laws on $\mathscr H$. Then the previous conjectures assert the following

\begin{compactenum}

\smallskip

\item[$-$] The $n^{1/3}$ cases 
correspond to the half-line 
$\mathtt L\! =\!  \big\{ ( t,t, \frac{_{_1}}{^{^2}} ) ; t \! \in\!  \r_+^* \big\} $. 

\smallskip

\item[$-$] The $n^{1/2}$ cases correspond to the slit quadrant 
$ \mathtt P\! =\!  \big( \r_+^2  \! \times \! \{ \frac{_{_1}}{^{^2}} \} \big)\backslash (\mathtt L \cup \{ ( 0,0, \frac{_{_1}}{^{^2}}) \} ) $.

\smallskip

\item[$-$] The $n$ cases correspond to $\mathtt D^o\cup \mathtt B_+\cup \mathtt B_-$ where   
$ \mathtt D^o \! = \! \big( (\r^*_+)^2 \! \times \! [0, 1]\big)   \backslash (  (\r^*_+)^2 \! \times  \!  \{ \frac{_{_1}}{^{^2}} \} )$,  

\noindent  
$\mathtt B_+ \eqo  \big\{ (t,0, p); p \! \in\!  (\frac{_{_1}}{^{^2}} , 1], t \! \in\!  \r_+^* \big\}$ and 
$\mathtt B_- \eqo  \big\{ (0,t, p); p \! \in \!  [0,\frac{_{_1}}{^{^2}} ), t \! \in \! \r^*_+ \big\} \; .$

\smallskip

\item[$-$] The bounded cases correspond to $\mathtt B'_+\cup \mathtt B'_-$ where $\mathtt B'_+ \eqo \big\{ (t,0, p); p \! \in \!  [0, \frac{_{_1}}{^{^2}} ), t\! \in \! \r_+^* \big\}$ and $\mathtt B'_- \eqo  \big\{ (0,t, p); p\! \in \!   (\frac{_{_1}}{^{^2}} , 1],t\! \in \! \r_+^* \big\}$. 

\smallskip

\end{compactenum}
Let us mention that when $p$ varies the distance and the resistance models do not allow generic transitions because there is only one degree of freedom in the parameters of these models.  
{\it Indeed} in the distance model, we see that $\alpha_+ \geko 0$ does not depend on $p$ and that $\alpha_-\! = \! 0$. So $p \!\in \!  [0, \frac12)$  corresponds to bounded cases, $p\! = \! \frac12$  corresponds to an $n^{1/2}$ case and $p\ino (\frac12, 1]$ corresponds to $n$ cases: there is no possible transition to any $n^{1/3}$ case. 
In the resistance cases we see that $\alpha_+\! = \! \alpha_- \! >\! 0$ do not depend on $p$. So $p \! \in \! (0,1) \backslash \{ \frac12\}$ corresponds to  $n$ cases and $p\! = \! \frac12$  corresponds to an $n^{1/3}$ case: there is no possible transition to any $n^{1/2}$ case.

\bigskip 

\noindent 
{\large \textbf{Statement of the main result}}. $\; $ We assume that $X_n$ is as in (\ref{syshom}) with $\mathbf F \ino \mathscr H$ a.s. Under 
moment conditions we prove Derrida's conjecture $(c)$, namely the $n^{1/3} $ case. 
Let us mention that \emph{we do not assume that $\mathbf F$ is symmetric}. To express our assumptions we first need to introduce the following notations. Let $F \! \in \! \mathscr H$. For all $x,y\! \in \! \r_+^*$, we set 
\begin{equation} 
 \label{F_inv}
  F^\inv (x,  y)\eqo 
    \frac{1}{F(\frac1x \, ,  \frac1y)},
  \;     F^*(x,y)
 \eqo  \begin{cases}
        F (x,y) & \!\!\!\! \!\!\!\!\!\!  \hbox{\rm  if $F\ino \mathscr{H}_+$,}
        \\
        F^\inv(x,y) &\!\!\!\!   \hbox{\rm if $F\ino \mathscr{H}_-$,}
      \end{cases}   \; \textrm{and} \;   F^\# (x, \, y) \eqo F(y, \, x).
\end{equation}
We note that $F^\#, F^\inv\ino \mathscr H$, that $\sgn(F^\#)\eqo \sgn(F)\eqo - \sgn(F^\inv)$, 
that $F^* \in \mathscr{H}_+$ and that $F \mapsto (F^*, \, \sgn (F))$ defines a bijection between $\mathscr H$ and $\mathscr{H}_+\!  \times\!  \{ 1,  -\! 1\}$ (see Proposition \ref{corrprop} for more details). Then for any $F \! \in \!  \mathscr{H}$, 
we define for all $x\in \r_+^*$,
\begin{equation}
    T_F(x) =
    \inf \{ y\ino \r: \, F^* (\ee^{-x}, \, \ee^{-y}) \leko  1 \} \in [0, \, \infty) \, .
    \label{TF}
\end{equation}
We see that $T_F \eqo T_{F^*}$. The fact that $T_F(x) \in [0, \, \infty)$  comes from 
$\lim_{y\to \infty} F^* (\ee^{-x}, \, \ee^{-y}) = \ee^{-x} \leko 1$ and from $F^* (\ee^{-x}, \, 1) \geqo \varphi_{+}(\ee^{-x}, \, 1) \eqo 1$. By definition, for all $x, y\ino \r_+^*$ 
\begin{equation}
T_F(x) \leko  y
\,  \Longleftrightarrow  \, F^*(\ee^{-x}, \, \ee^{-y}) \leko  1
\, \Longleftrightarrow \, 
    \begin{cases}
        F(\ee^{-x}, \, \ee^{-y}) \leko 1 & \!\!\!\! \hbox{ \rm if $\sgn(F)\eqo 1$}, \\
        F(\ee^x, \, \ee^y) \geko  1 & \!\!\!\!  \hbox{ \rm if $\sgn(F)\eqo -1$}\, .
    \end{cases}
    \label{T<y}
\end{equation}

 Next for all $a \! \in \! \r_+$, $b \! \in \!  \r_+^*$ and $\eta \ino (0 ,1)$, we introduce the following quantities which are well-defined but possibly infinite.  
\begin{equation}
    \Gamma_{{\! F}}^{{a,b}}
    \! :=\! 
    \int_0^\infty \!\!\! \! x^a \, (T_F (x))^b \d x, \quad \mathtt{c} (F)     \! :=\! 
   \Gamma_{\! F}^{1,1}\! + \tfrac{1}{2}\Gamma_{\! F}^{0,2}  \quad \textrm{and} \quad M^{_{(\eta)}}_{\! F}\! :=\! 
  \Gamma_{\! F}^{1+\eta,1}\! + \Gamma_{\! F}^{0,2+\eta} .
     \label{Gamma}
\end{equation}
In Lemma \ref{Gammabasicprop} $(viii)$ we prove that if $ M^{_{(\eta)}}_{\! F}\leko \infty$, then for all $(a,b)\ino \r_+ \! \times \! [1, \infty)$ such that 
$a+b\leko 2+ \eta$, there is a constant $c_{a,b,\eta} \ino (0, \infty)$, which only depends on $a,b$ and $\eta$, such that 
\begin{equation}
\label{momentmono} 
\Gamma_{{\! F}}^{{a,b}} \leq c_{a,b,\eta} \max \big(1, M^{_{(\eta)}}_{\! F} \big) \; .
\end{equation}
In particular, it implies that $\mathtt c (F) \leqo (c_{1,1,\eta}+ c_{0, 2, \eta})    \max \big(1, M^{(\eta)}_{\! F} \big)$. 
It is now possible to state the main result of the paper.
\begin{theorem}
\label{t:main}
Let $X_0$ be an $\r_+^*$-valued r.v.~and let $(X_n)_{n\in \bbN}$ be as in (\ref{syshom}) with $\p (\mathbf F \! \in\!  \mathscr H) \! = \! 1$.  Let $\eta\ino (0, 1)$. 
We assume that 
\begin{equation}
\label{maintheohyp}
\p \big(\mathbf F \! \notin \! \{ \varphi_+, \varphi_-\} \big)\geko 0, \quad \e \big[ M^{_{(\eta)}}_{\! \mathbf F}\big]\leko \infty \quad  \textrm{and} \quad \e[\sgn (\ff) ] \eqo \e \big[ \Gamma_{\! \ff}^{0,1} \sgn (\ff) \big]\eqo 0 .
\end{equation}
We set 
\begin{equation}
     c_* := \tfrac92 \e \big[ \mathtt c(\ff) \big] \; .
    \label{cste}
\end{equation}
Then $ c_* \ino \r_+^*$ and for all $y\in \r$,
 \begin{equation}
     \lim_{n\to \infty} \p \Big( \frac{\log X_n}{(c_* n)^{1/3}} \le y \Big)
     =
     \int_{-\infty}^y \!\!\! \tfrac{3}{4}  \big( 1\! -\! x^2 \big)  {\bf 1}_{\{ |x| <1 \} } \, \d x \, .
     \label{wcv}
 \end{equation}
\end{theorem}
\noi
\textbf{Proof.} See Remarks \ref{cstarfiniteandpos} for a proof of $c_*\ino \r_+^*$ and see Section \ref{proofmainthsec} for the proof of (\ref{wcv}). \cqfd 
\begin{remark}
\label{r:main} 
 {\rm 
 $(i)$ If $\ff$ is symmetric, i.e.~if $\ff\eqo \ff^\#$, then 
$\e[\sgn (\ff) ] \eqo \e \big[ \Gamma_{\! \ff}^{0,1} \sgn (\ff) \big]\eqo 0$ means that $p$ as in (\ref{defiparamp}) is equal to $\frac12$ and that 
$\alpha_+$ and $\alpha_-$ as in (\ref{alphapm}) are equal.  

\smallskip

\noindent
$(ii)$ Note from (\ref{T<y}) that $T_\ff \eqo T_{\ff^*}\eqo  T_{\ff^\inv}$. Therefore $\ff$ satisfies (\ref{maintheohyp}) if and only if $\ff^\inv$ also satisfies (\ref{maintheohyp}). Moreover, note that $1/X_n$ is derived from $1/X_0$ as in (\ref{syshom}) with  
$\mathbf F$ is replaced by $\ff^\inv$. This shows that in order to prove Theorem \ref{t:main} we only need to prove that (\ref{maintheohyp}) implies 
 \begin{equation}
   \forall y\ino \r , \quad   \lim_{n\to \infty} \p \Big( \frac{\log X_n}{(c_* n)^{1/3}} \geko  y \Big)
     \geqo 
     \int^{\infty}_y \!\!\! \tfrac{3}{4}  \big( 1\! -\! x^2 \big)  {\bf 1}_{\{ |x| <1 \} } \, \d x \, .
     \label{halfwcv}
 \end{equation}

\noi
$(iii)$ Lemma \ref{Gammabasicprop} $(iii)$ asserts for all $a, b\ino [1, \infty)$ that 
\begin{equation}
 \label{IPP}
a\Gamma_{\! \ff^{\#}}^{a-1,b} \eqo b \Gamma_{\! \ff}^{b-1,a}\, .
\end{equation}
Since $\ff^\# (X_n, X'_n)$ has the same distribution as $\ff (X_n, X'_n)$, Theorem \ref{t:main} applies also to $\ff^\#$ in place of $\ff$ and leads to the same conclusion \eqref{wcv}. \cq 
 } 
\end{remark}

\begin{example}
\label{ex:symmetric}
 {\rm 
{\it (Symmetric cases and the resistance cases)} 
In the setting of Theorem \ref{t:main}, if for all $x, y\ino \r_+^*$, $\ff(x, \, y) \eqo \ff (y, \, x)$ a.s.~then 
$\e \big[ \Gamma_{\! \ff}^{0,2} \big] \! = \! 2\e \big[ \Gamma_{\! \ff}^{1,1} \big]$ by \eqref{IPP}. Therefore 
$c_* \! = \! 9 \e \big[ \Gamma_{\! \ff}^{1,1} \big]$ in this case.

For instance let $\boldsymbol{\alpha}$ be an $(\r \backslash \{ 0\})$-valued r.v.~and let us consider 
$F_{\! \boldsymbol{\alpha}} (x,y)\eqo  (x^{1/\boldsymbol{\alpha}} \! + y^{1/\boldsymbol{\alpha}} )^{\boldsymbol{\alpha}}$, $x,y\ino \r_+^*$, as introduced in Example \ref{1homogex} $(b)$. 
For all $x\in \bbR$, we set $\mathtt{s} (x)\! :=\!  {\bf 1}_{\{ x\geq 0\} } \! -\!  {\bf 1}_{\{ x<0\} }$, which is taken as the \emph{sign of $x$}. 
Then (\ref{maintheohyp}) is equivalent to assume that 
$$ \e\big[ |\boldsymbol{\alpha}|^{3+\eta}\big] <\infty \quad \textrm{and} \quad \e[\mathtt{s}  (\boldsymbol{\alpha})]\eqo   \e[\boldsymbol{\alpha}^2  \mathtt{s} (\boldsymbol{\alpha})] \eqo 0, $$
Under this assumption (\ref{wcv}) holds true with $c_* \eqo 9 \zeta(3) \e \big[ |\boldsymbol{\alpha}|^3 \big]$. In particular 
the effective resistance of the critical random series-parallel graph corresponds to 
$\p(\boldsymbol{\alpha} \eqo1) \eqo \p( \boldsymbol{\alpha} \eqo -1)\eqo \frac12$ 
and this proves Conjectures \ref{conj:addario-berry} and \ref{conj:zeta(3)} with $c_{\mathtt{eff}} \eqo 9 \zeta(3)$. \cq 
 } 
\end{example}
\begin{example}
\label{Hipstersuite}
{\rm 
 (\emph{Symmetric hipster RW}) As already mentioned the symmetric hipster random walk 
corresponds to $\mathbf F$ such that $\p (\mathbf G \eqo G_{\mathtt{hip}}\,  ; \,  \varepsilon(\mathbf F) \eqo 1)\eqo  
\p (\mathbf G \! \equiv\!  0 \,  ; \,  \varepsilon(\mathbf F) \eqo -1)\eqo\frac12$ where $G_{\mathtt{hip}} (z) \eqo  \max (0, 1\! -\! |z|) $, $z \! \in \! \r$. 
Namely $\p (\mathbf F\! = \! \fg^{\mathtt{hip}}_+)\! = \! \p (\mathbf F\! = \! \fg^{\mathtt{hip}}_-)\! = \! \frac12$, where for all $x,y \! \in \! \r$, 
$$\log \fg^{\mathtt{hip}}_{\pm} (\mathrm e^x, \mathrm e^y )\! = \! \varphi_{\pm} (x,y) \pm G_{\mathtt{hip}} (\pm |x\! -\! y|) \; .$$ 
We observe that  $\fg^{\mathtt{hip}}_- \eqo (\fg^{\mathtt{hip}}_+)^\inv$ and that 
$T_{\fg^{\mathtt{hip}}_+} (t)\eqo  T_{\fg^{\mathtt{hip}}_-} (t) \eqo  {\bf 1}_{(0, \, 1]}(t)$ for $t\in \r_+^*$. 
Hence 
$$\Gamma_{\! \fg^{\mathtt{hip}}_\pm}^{0,2} \eqo   2\Gamma_{\! \fg^{\mathtt{hip}}_\pm}^{1,1} \eqo 1 \quad 
\textrm{and} \quad c_*\eqo  \tfrac{9}{2} \, \e \big[ \mathtt c \big( \ff^{\mathtt{hip}}\big) \big]\eqo   \tfrac{9}{2} .$$
Theorem \ref{t:main} gives in this case the same conclusion as Theorem \ref{t:hipster}, with $
c_{\mathtt{hip}}\! =  \!  \frac92 $, as expected. \cq 
 } 
\end{example}

\medskip

\noindent
{\large \textbf{A new paper by P.~Morfe}}. $\; $ Theorem \ref{t:main} is (almost) identical to what has been independently established by Morfe~\cite{morfe}, whose manuscript was posted on arXiv a few days before ours. Studying the recursive distributional equation satisfied by the cumulative distribution function of the logarithm of homogeneous systems, Morfe~\cite{morfe} proved that upon a suitable normalization, the cumulative distribution function converges, in an appropriate function space, to a limiting function satisfying a certain partial differential equation (PDE). Since it was known (\cite{addario-berry-cairns-devroye-kerriou-mitchell}) that this PDE characterizes the cumulative distribution function of the law with density $\frac34  (1\! -\! x^2)  {\bf 1}_{\{ |x| <1 \} }$, this yields the desired conclusion. The method in Morfe~\cite{morfe} relies on an adroit analysis of the recursive distributional equation by means of PDE techniques (extending the approach of Barles and Souganidis~\cite{barles-souganidis}). Our method is more bare-handed. The basic idea, going back at least to \cite{auffinger-cable} in the study of the distance in the critical random series-parallel graph, consists in establishing that the cumulative distribution function of homogeneous systems is bounded from above by an explicit (cumulative distribution) function and also from below by another explicit (cumulative distribution) function. This is done by an induction argument on $n$. While we are guided by the statement of Conjecture \ref{conj:zeta(3)} to find appropriate candidates for the explicit functions,
it is technically involving to prove that the comparison is preserved under the operation (in law) $X_n \mapsto \ff (X_n, \, X'_n)$, due to the fact that the random function $\ff$ takes values in both $\mathscr{H}_+$ and $\mathscr{H}_-$ with positive probability. This is the main difficulty in our proof. 

Let us also mention a few other differences between the work of Morfe~\cite{morfe} and ours. The integrability assumption in \cite{morfe} is slightly weaker; it requires only the finiteness of $\e \big[ \Gamma_{\! \ff}^{1,1} \big]$ and $\e \big[ \Gamma_{\! \ff}^{0,2} \big]$ in Theorem \ref{t:main} (in our notation). On the other hand, Morfe~\cite{morfe} restricted attention to symmetric functions in $\mathscr{H}_+$ and $\mathscr{H}_-$, whereas in our paper, the symmetry of functions plays no particular role, leading only to simplification of some arguments and some constants (for the latter, see for example the beginning of the discussions in Example \ref{ex:symmetric}). Also, Morfe~\cite{morfe} was able to treat the two cases $\e[ \Gamma_{\! \ff}^{0,1} \sgn (\ff) ]\eqo 0$ and $\e[ \Gamma_{\! \ff}^{0,1} \sgn (\ff) ] \! \neq \!  0$ (assuming $\e[ \sgn (\ff) ]\eqo 0$ as in our paper) in a relatively unified setting, whereas our method, even though it also applies to both cases, requires different calculations, and for this reason we plan to write the proof for the case $\e[ \Gamma_{\! \ff}^{0,1} \sgn (\ff) ] \! \neq \!  0$ in a separate paper where the calculations are comfortably shorter than those in the present paper.

\medskip

\noindent
{\large \textbf{Acknowledgement}}. $\;$  We warmly thank B.~Derrida for the many fruitful discussions that led to this work. In particular he devised 
the general framework of homogeneous functions adopted in this article in order to unify previous models and 
he is the source of Conjecture \ref{conj4orders}.

\medskip

\noindent
{\large \textbf{Organization of the paper}}. $\;$  The rest of the paper is as follows. 
In Section \ref{prelimsec} we state and prove basic estimates that are needed in the proof of Theorem \ref{t:main}. The proof of this theorem is done at Section \ref{proofmainsec} and it is a consequence of a general result: Theorem \ref{p:Lambda_cond=>lb}, which is stated and proved in Section \ref{s:proof1}. Its proof is not technical. Theorem \ref{p:Lambda_cond=>lb} holds under a specific condition (\ref{Lambdahyp}) which is proved in Proposition \ref{p:main}. The proof of Proposition \ref{p:main} is done in Section \ref{Lambdachecksec}: it is the technical part of the proof of Theorem \ref{t:main} where all the estimates from Section \ref{prelimsec} are used.

\section{Preliminary estimates}
\label{prelimsec}
\subsection{Basic properties of functions of $\mathscr H$}
\label{basichomogsec}

We recall the space $\mathscr H$ from Definition \ref{classH} and the space $\mathscr G$ from Definition \ref{classG}. 
We also recall from (\ref{F_inv}) the definition of $F^\inv$, $F^*$ and $F^\#$ and from Example \ref{1homogex} $(a)$ that $\phi_+ (x,y)\eqo \max (x,y)$ and $\varphi_- (x,y)\eqo \min (x,y)$.  
We first prove the following proposition which states the one-to-one correspondence given by (\ref{correspond}). 
\begin{proposition}
\label{corrprop} The following holds true. 
\begin{compactenum}

\smallskip

\item[$(i)$] There is a one-to-one correspondence between $F\ino \mathscr H$ and 
$(\varepsilon, G)\ino \{ -1, 1\} \! \times \! \mathscr G$ given by (\ref{correspond}): namely, $\varepsilon (F)\eqo  \varepsilon$
and $\log F(\mathrm{e}^x, \mathrm{e}^y) \eqo \varphi_\varepsilon (x,y)+ \varepsilon G \big(\varepsilon (x\! -\! y) \big)$ for all  $x,y\ino \r$. In this case $\epp \eqo \epp (F)$ and it is convenient to use the notation $G_F$ for $G$.

\smallskip

\item[$(ii)$] For all $F\ino \mathscr H$, $\epp (F^*)\eqo 1$, $\epp (F)\eqo \epp (F^\#)\eqo -\epp (F^\inv)$ and for 
$z\ino \bbR$, $G_F (z) \eqo G_{F^\inv} (z)\eqo G_{F^*} (z) \eqo G_{F^\# } (-z)$. 
\end{compactenum}
\end{proposition}
\noi
\textbf{Proof.} Let $F\ino \mathscr H$. For all $z\ino \bbR$, we set $G_F (z)\eqo \epp (F) \log F\big( \mathrm e^{ \epp (F) z}, 1 \big)\! -\! (z)_+$, where as before we denote the positive part of $z$ by $(z)_+\! :=\! \max (0, z)$. Since $F$ is $1$-homogeneous, for all  $x,y\ino \bbR$, we get $\log F(\mathrm e^x, \mathrm e^y)\eqo \log F(\mathrm e^{x-y}, 1)+y $, which easily implies (\ref{correspond}). 
We also see that $G_F$ is continuous 
and we easily check that $G_{F^\inv}\eqo G_F$ and thus $G_{F^*}\eqo G_F$. Since $F$ is $1$-homogeneous and since $\epp (F^\#)\eqo \epp (F)$, we also 
get $G_{F^\#} (-z)\eqo \epp (F)\log F\big( 1, \mathrm e^{-\epp (F)z }\big) \! -\! (-z)_+\eqo -z \! -\! (-z)_+ +  \epp (F) \log F \big( \mathrm e^{\epp (F) z}, 1\big)\eqo G_F(z)$. 

We next prove that $G_F\ino \mathscr G$. We first consider the case where $\epp (F)\eqo 1$. Let $z,h\ino \bbR_+$. 
Since $F$ is nondecreasing in its first coordinate, we first get 
$ G_F(z)\!-\! G_F(z+h) \eqo  \log F( \mathrm e^{z}, 1) \! -\!   \log F( \mathrm e^{z+h}, 1) + h  \leqo h$ and since $F$ is $1$-homogeneous we next see that $ G_F(z)\!-\! G_F(z+h) \eqo  \log F( \mathrm e^{z}, 1) \! -\!   \log F( \mathrm e^{z}, \ee^{-h}) $, which is a nonnegative quantity  since $F$ is nondecreasing in its second coordinate. Namely $0\leqo G_F (z)\! -\! G_{F} (z+h) \leqo h$. Similarly, we see that $ G_F(-z)\!-\! G_F(-z-h) \eqo  \log F( \mathrm e^{-z}, 1) \! -\!   \log F( \mathrm e^{-z-h}, 1) \geqo 0$, 
since $F$ is nondecreasing in its first coordinate. Since $F$ is $1$-homogeneous, we also get 
$ G_F(-z)\!-\! G_F(-z-h) \eqo  \log F( \mathrm e^{-z}, 1) \! -\!   \log F( \mathrm e^{-z}, \ee^h) +h $, which is $\leqo h$ since  $F$ is nondecreasing in its second coordinate. 
Thus $G_F$ satisfies Definition \ref{classG} $(a)$.

We next observe that Definition \ref{classH} $(c)$ implies $\lim_{z\to -\infty} G_F(z)\eqo 0$ if $\epp (F)\eqo 1$. 
We can also apply this result to $F^\#$ since $\epp (F^\#)\eqo \epp (F)\eqo 1$ and we get $\lim_{z\to \infty} G_F(z)\eqo \lim_{z\to -\infty} G_{F^\#} (z)\eqo 0$. This proves that $G_F\ino \mathscr G$ if $\epp (F)\eqo 1$. 

If $\epp (F)\eqo -1$, then $\epp (F^\inv) \eqo 1$ and the previous result implies that $G_{F^\inv} \ino \mathscr G$. Since $G_F\eqo G_{F^\inv}$, this shows that $G_F\ino \mathscr G$. 

Conversely, let $(\varepsilon, G)\ino \{ -1, 1\} \! \times \! \mathscr G$ and let $F$ be defined by (\ref{correspond}). We immediately see that $F$ is a continuous $1$-homogeneous function which satisfies Definition \ref{classH} $(b)$ and $(c)$. Since $G$ is Lipschitz on $\bbR$, Lebesgue's differentiation theorem asserts that $G$ is differentiable Lebesgue almost everywhere on $\bbR$ and $G(z)\! = \! G(0)+ \int_0^z G'(x) \mathrm d x$, where $G'$ stands for the derivative of $G$ and where we use the usual convention that $ \int_0^z G'(x) \mathrm d x \eqo -\int_{z}^{0} G'(x) \mathrm d x$ if $z\leko 0$. By Definition \ref{classG} $(a)$ we get $\un_{\bbR_+} (z)+ G'(z) \geqo 0$, for almost all $z\ino \bbR$. 
Let $y\ino \bbR$. We note that $\frac{\mathrm d}{\mathrm dx} \varphi_\epp (x,y)\eqo \un_{\bbR_+} \big( \epp (x\! -\! y) \big)$ for all $x\ino \bbR\backslash \{ y\}$ . Therefore $\frac{\mathrm d}{\mathrm dx} \log F( \mathrm e^x, \mathrm e^{y}) \eqo  \un_{\bbR_+} \big( \epp (x\! -\! y) \big) + G' \big( \epp (x\! -\! y) \big) \geqo 0$ for almost all $x\ino \bbR$. Similarly 
$\frac{\mathrm d}{\mathrm dy} \log F( \mathrm e^x, \mathrm e^{y}) \geqo  0$, for almost all $y\ino \bbR$. Thus $F$ is nondecreasing with respect to each coordinates. This proves that $F\ino \mathscr H$ and it completes the proof of the proposition. \cqfd 

\medskip

Next we fix $F\ino \mathscr H$ and we recall from (\ref{TF}) the definition of $T_F$ and from (\ref{Gamma}) the definition of $\Gamma^{a,b}_{\! F}$. The following lemma provides basic properties of $T_F$ and
$\Gamma^{a,b}_{\! F}$, which are used throughout the article. 
\begin{lemma}
\label{Gammabasicprop} Let $F\ino \mathscr H$. We set $G\eqo G_F\ino \mathscr G$ as in (\ref{correspond}). We fix $(a,b)\ino \bbR_+ \! \times \! (0, \infty)$ and $\eta \ino (0, \infty)$. We define $A \eqo \big\{ (x,y)\ino \bbR_+^2: G(y\! -\! x) \geq \min (x,y) \big\} $ which is not symmetric in general and for all measurable functions $Q\! :\! \bbR_+^2 \! \to \! \bbR_+$, 
we write  
$\iint_A Q(x,y) \, \mathrm dx \mathrm dy $ to mean $\iint_{\bbR^2_+} Q(x,y) \un_A (x,y) \,  \mathrm dx \mathrm dy$. 
We also recall that $M^{_{(\eta)}}_{\! F} \eqo \Gamma^{1+\eta, 1}_{\! F}\! + \Gamma^{0, 2+\eta}_{\! F}$. 
Then the following holds true. 
\begin{compactenum}

\smallskip

\item[$(i)$] The function $x\ino \r^*_+ \mapsto T_F (x) \ino [0, \infty)$ is nonincreasing and its right-limit at $x$ is given by 
$T_F (x+)\eqo \inf \{ y\ino \bbR_+ : F^*(\ee^{-x}, \ee^{-y}) \leqo 1 \}$.  
Moreover for all $(x, y)\ino \bbR_+\! \times \! \bbR$, 
$T_F(x) \geqo y$ if and only if $G(y\! -\! x) \geqo \min (x,y)$.

\smallskip

\item[$(ii)$] $\Gamma^{a,b}_{\! F}\eqo b\iint_{A} x^ay^{b-1} \mathrm dx \mathrm d y $. Moreover 
$\Gamma^{a,b}_{\! F}\eqo 0$ if and only if $F\! \ino \! \{ \varphi_+, \varphi_-\}$, which is equivalent to $G \equiv 0$. 

\smallskip

\item[$(iii)$] $(a+1) \Gamma^{a,b+1}_{\! F^\#}\!  \eqo (b+1) \Gamma^{b, a+1}_{\! F}$ .

\smallskip

\item[$(iv)$] $\iint_A \max (x,y)^{1+\eta}  \mathrm dx \mathrm d y\leqo  M^{_{(\eta)}}_{\! F}\; $ and  
$\; \iint_A \max (x,y)^{\eta}   (x+y)\mathrm dx \mathrm d y\leqo  2M^{_{(\eta)}}_{\! F}$.

\smallskip

\item[$(v)$] For all $\gamma \ino (0, \infty)$, $\iint_{A\cap [\gamma, \infty)^2}  \mathrm dx \mathrm d y\eqo \int_\bbR \big( G(z) \! -\! \gamma)_+ \mathrm dz  \leqo   \int_\bbR G(z)  \mathrm dz \eqo \Gamma^{0,1}_{\! F}$.

\smallskip

\item[$(vi)$] $ T_F (G(0)) \geqo G(0) \geqo T_F(G(0)+)$.

\smallskip

\item[$(vii)$] $\sup_{x\in \bbR^*_+} x^{a+1} T_F(x)^b \leqo (a+1) \Gamma^{a,b}_{\! F}$. In particular $G(0)^{a+b+1}\leqo (a+1) \Gamma^{a,b}_{\! F}$.

\smallskip

\item[$(viii)$] Let $q\ino (0, \infty)$. We assume that $\Gamma^{q, 1}_{\! F}\! + \Gamma^{0, 1+q}_{\! F}\leko \infty$. 
Then for all $(a,b)\ino \bbR_+ \! \times \! [1, \infty)$ such that $a+b\leko 1+ q$, there exists a constant $k \ino (0, \infty)$, such that $\Gamma^{a,b}_{\! F} \leq k\max \big( 1, \Gamma^{q, 1}_{\! F}\! + \Gamma^{0, 1+q}_{\! F}\big)$.  
In particular, if $M^{_{(\eta)}}_{\! F}\leko \infty$, then for all $(a,b)\ino \bbR_+ \! \times \! [1, \infty)$ such that $a+b\leko 2+ \eta$, there exists $c_{a,b,\eta} \ino (0, \infty)$, which only depends on $a,b$ and $\eta$, such that 
$\Gamma^{a,b}_{\! F} \leq c_{a,b,\eta} \max \big( 1,M^{_{(\eta)}}_{\! F}\big)$.

\smallskip

\item[$(ix)$]  We assume $M^{_{(\eta)}}_{\! F}\leko \infty$ and $G(0)\! \neq \! 0$, i.e.~$F\! \notin \! \{ \varphi_-, \varphi_+\}$ (by $(ii)$ of the present lemma). Then $\Gamma^{a,b}_{\! F} \leqo M^{_{(\eta)}}_{\! F} G(0)^{a+b-2-\eta}$ for all $(a,b)\ino [0, 1+\eta] \! \times \! [1, 2+ \eta]$.  
\end{compactenum}
\end{lemma}
\noi
\textbf{Proof.} To prove $(i)$, 
we first observe for all $x,y,h\ino (0, \infty)$ that 
$F^*(\ee^{-x}, \ee^{-y})\geqo  F^*(\ee^{-x-h}, \ee^{-y})$, since $F^*$ is nondecreasing in its coordinates. This entails that $x\ino (0, \infty) \mapsto T_F (x) \ino [0, \infty)$ is nonincreasing. 

We next deduce from (\ref{T<y}) that $T_F(x) \geqo y$ is equivalent to $F^*(\ee^{-x}, \ee^{-y})\geqo 1$ which actually means 
$G(y\! -\! x) \geqo \min (x,y)$. 

Let us compute the right-limit of $T_F$ at $x$. 
We set $\mathcal T_{\! x}\eqo \inf \{ y\ino \bbR : F^*(\ee^{-x}, \ee^{-y}) \leqo 1 \}$. 
Since $F^*$ is continuous, for all $z\ino \bbR^*_+$, $F^*( \ee^{-z}, \ee^{-T_F(z)}) \leqo 1$ (this inequality is actually an equality).
If $z \! \downarrow \! \downarrow \! x$, it entails $F^*( \ee^{-x}, \ee^{-T_F(x+)}) \leqo 1$ and thus $T_F(x+)\geqo \mathcal T_{\! x}$. The continuity of $F$ also implies,  
$F^*( \ee^{-x}, \ee^{-\mathcal T_x}) \leqo 1$. Since $F$ is $1$-homogeneous, for all $h\ino \bbR_+^*$ we get 
$F^* (\ee^{-x-h}, \ee^{-(\mathcal T_x +h) })\eqo \ee^{-h} F^*( \ee^{-x}, \ee^{-\mathcal T_x}) \leko 1$, which implies 
$T_F (x+h) \leqo \mathcal T_{\! x} + h$ and thus $T_F (x+) \leqo  \mathcal T_{\! x}$ as $h\! \downarrow \! 0$. This entails $T_F (x+) \eqo \mathcal T_{\! x}$, which completes the proof of $(i)$. 

\smallskip

Let us prove $(ii)$. By $(i)$, for all $x,y\ino \bbR^*_+$, $T_F(x) \geqo y$ if and only if $(x,y)\ino A$. By Fubini we thus get 
$\Gamma^{a,b}_{\! F}\eqo b\iint_{\bbR_+^2}  x^ay^{b-1}\un_{\{ T_F(x) \geq y \}} \mathrm dx \mathrm dy \eqo  b\iint_{A}  x^ay^{b-1} \mathrm dx \mathrm dy$. 

Let us prove the second statement of $(ii)$. We first deduce from (\ref{correspond}) that $F\! \ino \! \{ \varphi_+, \varphi_-\}$ is equivalent to $G \equiv 0$. Let us suppose that $G\equiv 0$. Then $A \! \subset \! (\bbR_+\! \times \! \{ 0\})\cup (\{0\} \! \times \!  \bbR_+)$. Therefore $A$ is Lebesgue negligible and $\Gamma^{a,b}_{\! F}\eqo 0$. Conversely, let us assume that $G$ is not identically null. Therefore $G(0)\geko 0$ and by continuity there is $\delta \ino \r_+^*$ such that $[0, \delta)^2 \!\subset \! \{ (x,y)\ino \bbR^2_+ : 
G(y\! -\! x) \geko \min (x,y) \} \! \subset \! A$. Then the first statement of $(ii)$ implies that $0\leko \tfrac{1}{a+1} \delta^{a+b+1} \leqo \Gamma^{a,b}_{\! F}$, which completes the proof of $(ii)$. 

\smallskip

 To prove $(iii)$ we introduce the function $\mathtt{s}\! :\! (x,y) \ino  \bbR^2_+ \!  \mapsto \! (y,x)$. By Proposition \ref{corrprop} $(ii)$, $G_{F^\#}\eqo G (-\, \cdot)$. Therefore $\mathtt s(A)\eqo \big\{ (x,y)\ino \bbR^2_+: G_{F^\#} (y\! -\! x) \geqo \min (x,y) \big\}$ and $(iii)$ is a consequence of $(ii)$ and the change of variable $(x',y')\eqo \mathtt s(x,y)$.

\smallskip

To prove $(iv)$, we note that $\max (x,y)^{1+\eta} \leqo x^{1+\eta}\! + y^{1+\eta}$. By $(ii)$, 
$\iint_A \max (x,y)^{1+\eta}  \mathrm dx \mathrm d y \leqo \Gamma^{1+\eta, 1}_{\! F} \! + \frac{1}{2+ \eta}  
\Gamma^{0,2+\eta}_{\! F} \leqo M^{_{(\eta)}}_{\! F}$. 
We next observe that $\max (x,y)^{\eta} (x+y)\leqo 2 \max (x,y)^{1+\eta} $ and we use the previous upper bound to complete the proof of $(iv)$.  

\smallskip

  We next prove $(v)$. To simplify we set $I(\gamma)\eqo \iint_{A\cap [\gamma, \infty)^2}  \mathrm dx \mathrm d y$ 
and we define 
$$A_\gamma\eqo \big\{ (x',y')\ino (\bbR^*_+)^2: G (y'\! -\! x')\! -\! \gamma \geqo \min (x',y') \big\}\; .$$ We observe that $(x,y)\ino A\cap (\gamma, \infty)^2$ if and only if $(x\! -\! \gamma, y\! -\! \gamma) \ino A_\gamma$. Therefore 
$I(\gamma)\eqo  \iint_{A_\gamma}  \mathrm dx \mathrm d y$ because $I(\gamma)\eqo \iint_{A\cap [\gamma, \infty)^2}  \mathrm dx \mathrm d y\eqo  \iint_{A\cap (\gamma, \infty)^2}  \mathrm dx \mathrm d y$. 
We next set 
$B_+ \eqo \big\{ (x,y)\ino \bbR^2_+: x\geko y \big\}$ and $B_-\eqo \big\{ (x,y)\ino \bbR^2_+: y\geqo x \big\}$. Taking $z\eqo y\! -\! x$, we first get $\iint_{A_\gamma \cap B_+} \mathrm dx \mathrm d y \eqo \int_{\bbR_-} \mathrm dz\int_{\bbR+} 
\mathrm dy \un_{\{ (G(z)-\gamma)_+ \geq y\}}\eqo  \int_{\bbR_-} (G(z)\! -\! \gamma)_+ \mathrm dz$. Similarly, we also get 
$\iint_{A_\gamma \cap B_-} \mathrm dx \mathrm d y \eqo \int_{\bbR_+} \mathrm dz\int_{\bbR+} 
\mathrm dx \un_{\{ (G(z)-\gamma)_+ \geq x\}}\eqo  \int_{\bbR_+} (G(z)\! -\! \gamma)_+ \mathrm dz$. 
Thus $I(\gamma) \eqo  \int_\bbR \big( G(z) \! -\! \gamma)_+ \mathrm dz  $. Taking $\gamma\eqo 0$, we see that $A\eqo A_0$ and by $(ii)$ we get $\Gamma^{0,1}_{\! F}\eqo \iint_{A}  \mathrm dx \mathrm d y\eqo \int_{\bbR} G(z) \mathrm dz$.

\smallskip

Let us prove $(vi)$. Since 
$G\big( G(0)\! -\! G(0)\big)\eqo G(0) \eqo \min (G(0), G(0))$, by $(i)$ we first 
get $T_F(G(0))\geqo G(0)$. By $(i)$ again, we also see that $T_F(x+)\eqo \inf\{ y\ino \bbR_+: G(y\! -\! x) \leqo \min (x,y)\} $. Therefore $T_F (G(0)+)\leqo G(0)$.

\smallskip

To prove $(vii)$ we observe that $y^a T_F(y)^b \geqo y^a T_F(x)^b$ 
for all positive real numbers $x\geqo y$, since $T_F$ is nonincreasing by $(i)$. 
Thus $\Gamma^{a,b}_{\! F} \geqo \int_0^x  y^a T_F(y)^b \mathrm dy \geqo 
\frac{1}{a+1} x^{a+1}T_F(x)^b$, which implies the first inequality in $(vii)$. 
The second one is a consequence of $(vi)$ which implies 
that $G(0)^{a+b+1} \leq G(0)^{a+1} T_F(G(0))^b\leqo \sup_{x\in \bbR^*_+} x^{a+1} T_F(x)^b$. 

\smallskip

To prove $(viii)$ we assume that $M\! :=\! \Gamma^{q, 1}_{\! F}\! + \Gamma^{0, 1+q}_{\! F}\leko \infty$, where $q$ is some 
positive real number. Let 
$(a,b)\ino \bbR_+ \! \times \! [1, \infty)$ such that $a+b\leko 1+q$. 
We first suppose $x\ino [0, 1]$. By $(vii)$, $xT_F(x)^{1+q} \leqo \Gamma^{0, 1+q}_{\! F}$. Thus $x^aT_F(x)^b\leqo \big( \Gamma^{0, 1+q}_{\! F}\big)^{b/(1+q)} x^{a- \frac{b}{1+q}}$. Observe that $k_{a,b,q} \! :=\! a-\frac{b}{1+q} \geko a\! -\! 1\geqo -1$ since $b\leko 1+q$ and $a\geqo 0$.  Therefore we get 
$$ \int_0^1 \!\!\! x^aT_F(x)^b \, \mathrm dx \leq \tfrac{1}{k_{a,b,q} +1} \big( \Gamma^{0, 1+q}_{\! F}\big)^{b/(1+q)}  \leq  \, \tfrac{1}{k_{a,b,q} +1} \max \big( 1,  M\big)$$
because $b/(1+q) \leqo 1$.

We next fix $x\ino [1, \infty)$ and we set $\kappa \eqo 2+q \! -\! a\! -\! b \geko 1$ and $\theta \eqo (1+q \! -\! b)/q$. Note that 
$0\leqo a/q\leqo \theta \leqo 1 $. We then observe the following 
\begin{eqnarray*}
x^aT_F(x)^b &\eqo & x^{-\kappa} \big( x^{1+q} T_F(x)\big)^\theta \big( x T_F(x)^{1+q}\big)^{1-\theta}  \\
&\leq & x^{-\kappa} (1+q)^\theta  \big( \Gamma^{q, 1}_{\! F} \big)^\theta \big(\Gamma^{0, 1+q}_{\! F} \big)^{1-\theta} \leqo  x^{-\kappa} (1+q)^\theta  \max \big( 1, M\big) .
\end{eqnarray*}
We then get $\int_1^\infty x^aT_F(x)^b \, \mathrm dx \leq \frac{(1+q)^\theta}{\kappa -1} \max \big( 1, M\big) $, which completes the proof of $(viii)$.

\smallskip

We finally prove $(ix)$. We first suppose that 
$x\geko  G(0)$. We then see that 
$$ 
\frac{x^a T_F (x)^b }{G(0)^{a+b-2-\eta} \, x^{1+\eta} T_F(x)}\eqo \Big(\frac{G(0)}{x}\Big)^{\! 1+\eta-a} \Big(\frac{T_F(x)}{G(0)} \Big)^{\! b-1} \leqo 1$$ 
because $T_F (x) \leqo T_F (G(0)+) \leqo G(0)$ (by $(i)$ and $(vi)$), because $b\geqo 1$ and because $a\leqo 1+ \eta$.  

If $x \leko G(0)$, 
$$ \frac{x^a T_F (x)^b }{G(0)^{a+b-2-\eta} \,  T_F(x)^{2+\eta}}\eqo \Big( \frac{x}{G(0) }\Big)^{a} \Big( \frac{G(0)}{T_F(x)} \Big)^{2+\eta -b} \leqo 1 $$ 
because $T_F (x) \geqo T_F(G(0))\geqo G(0)$ (by $(i)$ and $(vi)$) and because $b\leqo 2+ \eta$. Therefore 
$$ \Gamma^{a,b}_{\! F} \leq  G(0)^{a+b-2-\eta} \int_{G(0)}^\infty \!\!\!\! x^{1+\eta} T_F(x) \, \mathrm dx + 
G(0)^{a+b-2-\eta} \int_0^{G(0)} \!\!\!\! T_F(x)^{2+\eta} \, \mathrm dx ,$$
which entails the desired inequality.  \cqfd 

\begin{remark}
\label{cstarfiniteandpos} 
 {\rm 
Let $\ff$ be an $\mathscr H$-valued random function satisfying the first 
two conditions of (\ref{maintheohyp}), namely, $\p \big(\mathbf F \! \notin \! \{ \varphi_+, \varphi_-\} \big)\geko 0$ and 
$\e \big[ M^{_{(\eta)}}_{\! \mathbf F}\big]\leko \infty$. We explain here why it implies that $c_* \ino \r^*_+$ as stated in Theorem \ref{t:main}. First, by Lemma \ref{Gammabasicprop} 
$(ii)$, a.s.~$\Gamma^{{a,b}}_{{\! \ff}} \geko 0$ 
for all $(a,b)\ino \bbR_+ \! \times \! (0, \infty)$, which in particular implies that $\e [\Gamma^{0,1}_{\! \ff} ]$, $\e [\Gamma^{a,b}_{\! \ff} ]$ and $\e [\Gamma^{a,b}_{\! \ff} ]$ are positive. By Lemma \ref{Gammabasicprop} $(viii)$, we also get 
$ \e [\Gamma^{0,1}_{\! \ff} ]\leqo c_{0,1,\eta } (1+ \e \big[ M_{\ff}^{_{(\eta)}}])$. Similarly, 
$ \e [\Gamma^{1,1}_{\! \ff} ]\leqo c_{1,1,\eta } (1+ \e \big[ M_{\ff}^{_{(\eta)}}])$ and $ \e [\Gamma^{0,2}_{\! \ff} ]\leqo c_{0,2,\eta } (1+ \e \big[ M_{\ff}^{_{(\eta)}}])$. Therefore $c_* \ino \r^*_+$. \cq
 } 
\end{remark}

\subsection{Basic properties of $\Lambda_{\psi, F}$}
\label{basicLambdasec}
Here we discuss basic properties of the distribution function of $\log F(\ee^{Y} ,\ee^{Y'})$ 
where $F\ino \mathscr H$ and $Y,Y'$ are two independent r.v.s~with the same probability density $\psi$. 
This is done through a function $\Lambda_{\psi,F}$ which is defined as follows and whose probabilistic interpretation is discussed in Lemma \ref{elempropLambda} below. 
\begin{definition}
\label{Lambdadef} 
{\rm 
 $(a)$ We denote by $\mathscr D$ the space of functions $\psi\! : \! \bbR \! \to \!\bbR_+$ which are probability densities, i.e.~Borel measurable functions such that $\int_{\bbR}\!  \psi (v)\,  \mathrm dv\eqo 1$.  We denote the corresponding distribution function by $\psib$. Namely, $\psib(v)\eqo \int_{-\infty}^v \! \psi (w) \, \mathrm dw$, for all 
$v\ino \bbR$. 

\smallskip

\noi
$(b)$ Let $\psi \ino \mathscr D$ and $F\ino \mathscr H$. We set $G\eqo G_F\eqo G_{F^*}$ and 
$A\eqo \{ (x,y)\ino \bbR_+^2\! : \! G(y\! -\! x) \geqo \min (x,y)\}$ as in (\ref{correspond}) and Proposition \ref{corrprop}, respectively, and 
\begin{equation}
\label{Lambdadefeff}
 \Lambda_{\psi, F} (v)= \epp(F)\iint_A \!\! \psi \big( v\! -\! \epp (F)x\big)\psi \big( v\! -\! \epp (F)y\big) \, \mathrm dx \mathrm dy \; ,
\end{equation}
for all $v\ino \bbR$. \cq 
 } 
\end{definition}
\begin{lemma}
\label{elempropLambda} Let $F\ino \mathscr H$ and let 
$Y,Y'$ be two independent r.v.s~with the same probability density $\psi \ino \mathscr D$. Then the following holds true. 
\begin{compactenum}

\smallskip

\item[$(i)$] For all $v\ino \bbR$, $\; \p \big(\! \log F\big(\ee^Y, \ee^{Y'}\big) \leqo v \big)= 2\psib(v)\un_{\{ \epp(F)= -1\}} + \epp(F) \psib(v)^2-\Lambda_{\psi,F} (v)$.

\smallskip

\item[$(ii)$] Let 
$\sigma', \sigma'' \ino \bbR$ be such that $\sigma'\leko \sigma''$. 
We suppose that $\psi$ is supported in $[\sigma', \sigma'']$. 
Let us furthermore assume that $\epp (F)\eqo 1$. Then 
$v\ino \bbR\! \mapsto \psib(v)^2 \! -\! \Lambda_{\psi, F} (v)$ is nondecreasing, $\Lambda_{\psi, F}$ is null on $(-\infty, \sigma']$, $\Lambda_{\psi, F}$ is nonincreasing on $[ \sigma'' , \infty)$ and $0\leqo \Lambda_{\psi, F} (v)\leqo \psib(v)^2$ for all $v\ino \bbR$.

\smallskip

\item[$(iii)$] Let $q\ino [0, 1)$. Let $\ff$ be a random $\mathscr H$-valued function such that $\e [\epp (\ff)]\eqo 0$. 
Let $Z,Z'$ be two independent 
$[-\infty, \infty)$-valued r.v.s whose law is given by $\p(Z\leqo v)\eqo \p(Z'\leqo v)\eqo q+ (1\! -\! q) \psib(v)$, $v\ino \bbR$, which implies $\p(Z\eqo -\infty)\eqo \p(Z'\eqo -\infty)\eqo q$. We assume that $\ff$, $Z$ and $Z'$ are independent. Then   
\begin{equation}
\label{centredLambda}
\forall v\ino \bbR, \qquad  \p \big(\! \log \ff \big(\ee^Z, \ee^{Z'}\big) \leqo v \big)= q + (1\! -\! q) \psib(v) -(1\! -\! q)^2 \e \big[ \Lambda_{\psi, \ff} (v) \big] 
\end{equation}
where $\ff$ is extended continuously on $\bbR_+^2$ as specified in (\ref{contiexthomog}) and where we adopt the conventions $\ee^{-\infty}\eqo 0$ and $\log (0)\eqo -\infty$. In particular, if $q\eqo 0$, the right-hand side of (\ref{centredLambda}) is simply $\psib (v) \! -\!  \e \big[ \Lambda_{\psi, \ff} (v) \big]$. 
 \end{compactenum}
\end{lemma}
\noi
\textbf{Proof.} We first prove $(i)$. To simplify notation we set $\epp \eqo \epp(F)$ and we set 
$G\eqo G_F\eqo G_{F^*}$ as in Proposition \ref{corrprop} $(ii)$. Note that  $\log F\big(\ee^Y, \ee^{Y'}\big) \leqo v $ if and only if 
$v\! -\! \varphi_\epp (Y,Y') \! -\! \epp G( \epp (Y\! -\! Y') ) \geqo 0$. We then 
set $X\eqo \epp v\! -\! \epp Y$ and $X'\eqo \epp v\! -\! \epp Y'$ and we observe that $G( \epp (Y\! -\! Y') )
\eqo G(X'\! -\! X)$ and that $v\! -\! \varphi_\epp (Y,Y') \eqo \epp \min (X,X')$. 
Therefore 
\begin{eqnarray}
\p \big(\! \log F\big(\ee^Y, \ee^{Y'}\big) \leqo v \big) \!\!\!\! & =& \!\!\!\! \p \big(\! 
 -\! \epp \big( G(X'\! -\! X) \! -\! \min (X,X')\big) \geqo 0\, \big) \nonumber \\ 
 = \!  \un_{\{ \epp = -1\}} \!\!\!\! \!\!\!\!  & & \!\!\!\! \!\!\!\!  \p \big( G(X'\!\!  -\! X) \geqo \min (X,X') \big) +
 \un_{\{ \epp = 1\}}  \p \big( G(X'\!\!  -\! X) \leqo \min (X,X') \big) \label{Lamstep}
\end{eqnarray}
We note $X$ and $X'$ admit the probability density $z\ino \bbR \! \mapsto \! \psi (v\! -\! \epp z)$. 
By a change of variable we see that $\p \big( G(X'\!\!  -\! X) \! = \!  \min (X,X') \big) \eqo 0$. Indeed, 
\begin{eqnarray*}
\p \big( G(X'\!\!  -\! X) \!\!\! \!\!  & =& \!\!\! \! \!  \min (X,X') \big) = \p \big( G(X'\!\!  -\! X) \eqo X ; X'\geqo X \big)+ 
\p \big( G(X'\!\!  -\! X) \eqo X' ; X\geko X' \big) \\
 \!\!\! \!\!  & =& \!\!\! \!   \int_{\bbR^*_+} \!\!\!\!\! \mathrm d y\!  \int_{\bbR}  \!\! \mathrm d z \,  \psi (v\! -\! \epp z)
 \psi (v\! -\! \epp z \! -\! \epp y) \un_{\{ z= G(y) \}} \\
 & & \quad +  \int_{\bbR_-}\!\!\!\!\!   \mathrm d y \int_{\bbR}\!\!  \mathrm d z' \, \psi (v\! -\! \epp z') 
\psi (v\! -\! \epp z' + \epp y) \un_{\{ z'= G(y) \}}  \; =\, 0.
\end{eqnarray*}
By (\ref{Lamstep}), we thus get 
\begin{eqnarray}
\p \big(\! \log F\big(\ee^Y, \ee^{Y'}\big) \leqo v \big) \!\!\!\! & =& \un_{\{ \epp = 1\}} \! -\! \epp \p \big( G(X'\!\!  -\! X) \geqo \min (X,X') \big) \nonumber \\
\!\!\!\! & =& \!\!\!\! \un_{\{ \epp = 1\}} \! -\! \epp \p \big( G(X'\!\!  -\! X) \geqo \min (X,X')\geqo 0 \big) -\epp \big( 1\! -\! \p(\min (X,X') \geqo 0)\big)  \nonumber\\
\!\!\!\! & =& \!\!\!\!  \un_{\{ \epp = 1\}} \! -\! \epp \p \big( (X, X')\ino A\big) -\epp + \epp \p(X\geq0)^2, \label{Lambdaprop1}
\end{eqnarray}
where we recall that $A\eqo \{ (x,y)\ino \bbR_+^2\! : \! G(y\! -\! x) \geqo \min (x,y) \}$. As already mentioned, $X$ and $X'$ admit the probability density $z\ino \bbR \! \mapsto \! \psi (v\! -\! \epp z)$, so we get $\epp \p \big( (X, X')\ino A\big) \eqo \Lambda_{\psi, F} (v)$. Since $\p (X\geqo 0)\eqo \un_{\{ \epp =-1\}} + \epp \psib(v)$, this implies the desired result thanks to (\ref{Lambdaprop1}).

Let us prove $(ii)$. We assume that $\epp \eqo \epp (F)\eqo 1$ and that $\psi$ has support in 
$[\sigma', \sigma'' ]$.  By (\ref{Lambdadefeff}) the function  $\Lambda_{\psi, F}$ only takes nonnegative values.  If $v \leko \sigma'$, then $\psi(v\! -\! z)\eqo 0$ for all $z\ino \bbR_+$ and (\ref{Lambdadefeff})
implies that $\Lambda_{\psi, F} (v)\eqo 0$. Moreover by $(i)$, $ \psib(v)^2 \! -\! \Lambda_{\psi, F} (v)\eqo \p \big(\! \log F\big(\ee^Y, \ee^{Y'}\big) \leqo v \big)$ is nondecreasing in $v$. Next, if $v\geqo \sigma''$, we get 
$\Lambda_{\psi, F} (v)\eqo 1\! -\! \p \big(\! \log F\big(\ee^Y, \ee^{Y'}\big) \leqo v \big)$, which is nonincreasing in $v$. This completes the proof of $(ii)$.

We next prove $(iii)$. Let $Y,Y'$ be two independent r.v.s with probability density $\psi$. We assume that $Y, Y'$ and $\ff$ are independent. We first note that $(Y,Y', \ff) $ has the same law as $(Z,Z', \ff)$ under $\p (\, \cdot \, | \, Z \, \textrm{and} \, Z' \geko -\infty)$. Then by $(i)$, we get 
$\p \big(\! \log \ff \big(\ee^Y, \ee^{Y'}\big) \leqo v \big)\! =\!  \psib(v)\! -\! \e \big[ \Lambda_{\psi,\ff} (v)\big] $, because $\e \big[\epp (\ff) \big]\eqo 0$ and $\p(\epp (\ff) \eqo -1)\eqo \frac12$. Therefore, for all $v\ino \bbR$, 
$$ \p \big(\! \log \ff \big(\ee^Z, \ee^{Z'}\big) \leqo v \, ; \,  Z \, \textrm{and} \, Z' \geko -\infty\big)= (1\! -\! q)^2 \Big(  \psib(v)\! -\! \e \big[ \Lambda_{\psi,\ff} (v)  \big] \Big)\; .$$

We next recall from (\ref{contiexthomog}) how $\ff$ can be extended continuously on $\bbR_+^2$. We note that 
$ (Z, \ff (\ee^Z, 0) )$ under $\p (\, \cdot \, | \, Z\geko -\infty  \, \textrm{and} \, Z' \eqo -\infty)$ and $ (Z', \ff (0, \ee^{Z'}) )$ under $\p (\, \cdot \, | \, Z'\geko -\infty  \, \textrm{and} \, Z\eqo -\infty)$ have the same law as $(Y,\ff(\ee^Y,0))$. Note that 
$\ff (\ee^Y, 0)\eqo \ee^Y\un_{\{ \epp (\ff) = 1\}}$ and since $\log 0\eqo -\infty$,
we get $\p (\log \ff (\ee^Y, 0) \leqo v)\eqo \e \big[ \un_{\{\epp (\ff) =-1 \}} +  \un_{\{\epp (\ff) =1 \}} \psib(v) \big]\eqo \frac{1}{2} \big( 1+\psib(v)\big)$. 

Since $\log \ff(0,0)\eqo -\infty$, we finally get 
$$   \p \big(\! \log \ff\big(\ee^Z, \ee^{Z'}\big) \leqo v \big)= q^2 + 2q(1\! -\! q) \tfrac{1}{2} \big( 1+\psib(v)\big) + (1\! -\! q)^2 \big(  \psib(v)\! -\! \e \big[ \Lambda_{\psi,\ff } (v) \big] \big) , $$
which implies (\ref{centredLambda}). This completes the proof of the lemma. 
\cqfd 

\medskip

The next lemma discusses more specific properties of $\Lambda_{\psi, F}$ when $\psi$ is even or when it is derived from  even functions. 
\begin{lemma}
\label{propLambdaeven} For all $v\ino \bbR$ we set $\mathtt s (v)\eqo \un_{\{ v\geq 0\}}\! -\! \un_{\{ v< 0\}}$, which is taken as  the sign of $v$. 
Let $F\ino \mathscr H$. We set $\epp\eqo \epp (F)$ and $G\eqo G_F\eqo G_{F^*}$ as in (\ref{correspond}) and Proposition \ref{corrprop} respectively, and $A\eqo \{ (x,y)\ino \bbR_+^2\! : \! G(y\! -\! x) \geqo \min (x,y)\}$. 
Let $\phi, \phi_+, \phi_- \ino \mathscr D$  be even. We set 
$$\forall v\ino \bbR, \quad \psi (v)\eqo \phi_+(v) \un_{\{ v\geq 0\}}+\phi_-(v) \un_{\{ v< 0\}}. $$
Then the following holds true. 
\begin{compactenum}

\smallskip

\item[$(i)$] $\psi \ino \mathscr D$, $\psib(0)\eqo \phib_+(0)\eqo \phib_-(0)\eqo \frac{1}{2}$ and for all $v\ino \bbR$, $\psib(v)\eqo \phib_+(v) \un_{\{ v\geq 0\}}+\phib_-(v) \un_{\{ v< 0\}}$. 
 
\smallskip

\item[$(ii)$] For all $v\ino \bbR$, $\; \epp  \Lambda_{\phi, F} ( \epp  v )\eqo \iint_A \phi (x\! -\! v) \phi(y\! -\! v) \, \mathrm dx \mathrm dy \eqo \Lambda_{\phi, F^*} (v)$.

\smallskip

\item[$(iii)$] Let $v\ino \bbR$ be such that $\mathtt{s} (v) \eqo -\epp$, then 
$\Lambda_{\psi, F} (v) \eqo \Lambda_{\phi_{-\epp} , F} (v)$ (here $\phi_{-\epp}\eqo \phi_-$ if $\epp \eqo 1$ and $\phi_{-\epp}\eqo \phi_+$ if $\epp \eqo -1$).

\smallskip

\item[$(iv)$] Let $\eta \ino (0, \infty)$. We recall $M^{_{(\eta)}}_{\! F}$ from (\ref{Gamma}) and we assume that $M^{_{(\eta)}}_{\! F} \leko \infty $. We furthermore assume the following. 
\begin{compactenum}

\smallskip

\item[$(a)$] $\max_{z\in \bbR} \psi (z)\eqo \psi(0)\eqo \phi_+(0)\eqo \phi_-(0)$ (which implies $\max \phi_\pm = \psi(0)$). 
 
\smallskip

\item[$(b)$] There are $\sigma, K\ino (0, \infty)$ such that for all $z\ino \bbR$, $\; \phi_+(z)\! -\! \phi_-(z) \leq Kz^2\un_{\{ |z| <\sigma\}}$. 

\smallskip

\end{compactenum}
\noi
Then for all $v\ino \bbR$, 
\begin{equation}
\label{symmetrization}
 \Lambda_{\psi, F} (v) \geq \Lambda_{\phi_{+} , F} (v) \un_{\{ v\geq 0\}}+ 
\Lambda_{\phi_{-} , F} (v) \un_{\{ v<0\}} -K\psi(0) \sigma^{1-\eta} M^{_{(\eta)}}_{\! F}\; .
\end{equation} 
 \end{compactenum}
\end{lemma}
\noi
\textbf{Proof.} 
To prove $(i)$, we first observe that $\int_{\bbR_-}\!  \phi_{\pm}  (z) \, \mathrm dz \eqo \int_{\bbR_+} \! \phi_{\pm}  (z) \, \mathrm dz\eqo \frac12$ since $\phi_\pm$ are even probability densities. This implies that $\psi\ino \mathscr D$ and that $\psib(0)\eqo \phib_+(0)\eqo \phib_-(0)\eqo \frac{1}{2}$. With the usual convention on the sign of integrals, it also implies 
for all $v\ino \bbR$ that 
$\psib(v)\eqo \frac12+ \int_0^v \psi(w) \, \mathrm dw$ and that $\phib_\pm(v)\eqo \frac12+ \int_0^v \phi_\pm(w) \, \mathrm dw$, which implies the desired result. 

We next prove $(ii)$. By (\ref{Lambdadefeff}), $\epp \Lambda_{\phi, F} (\epp v)\eqo \iint_A \phi (\epp (v\! -\! x)) \phi (\epp (v\! -\! y))  \, \mathrm dx \mathrm dy$, which entails $(ii)$ because $\phi (\epp (v\! -\! z))\eqo \phi (z\! -\! v)$ for all $z\ino \bbR$ since $\phi$ is even.

Let us prove $(iii)$. Since we assume $\mathtt s(v)\epp \eqo -1$, we get $-\epp v\eqo |v|$.
Therefore (\ref{Lambdadefeff}) implies for all $v\ino \bbR$, 
$\epp \Lambda_{\psi, F} (v)\eqo \iint_A \psi ( -\epp (x+ |v|)) \psi (-\epp (y +|v|)) \, \mathrm dx \mathrm dy$. Similarly, by $(ii)$ we get  $\epp \Lambda_{\phi_{\pm}, F} (v)\eqo \iint_A \phi_\pm (x+ |v|) \phi_\pm(y +|v|) \, \mathrm dx \mathrm dy$. 
Since for all $z\ino \bbR_+$, 
$\psi (-\epp (z +|v|)))\eqo \phi_{-\epp} (z+ |v|) $, this implies $ \epp \Lambda_{\psi, F} (v)\eqo \epp \Lambda_{\phi_{-\epp}, F} (v)$,  which completes the proof of $(iii)$.

To prove $(iv)$ we fix $v\ino \bbR$ and set $\psi_{[\epp]} (z)\eqo \psi (-\epp z)$ for all $z\ino \bbR$. 
By (\ref{Lambdadefeff}), $\Lambda_{\psi, F} (v)\eqo  \epp \iint_A \psi_{[\epp]} (x\! -\! \epp v) \psi_{[\epp]} (y\! -\! \epp v)  \, \mathrm dx  \mathrm dy$. We next set $A'\eqo \{ (x,y) \ino \bbR^2 \! : \! (x+ \epp v, y+ \epp v) \ino A \}$ and by the change of variables $(x',y')\eqo (x\! -\! \epp v, y\! -\! \epp v)$ we get $\Lambda_{\psi, F} (v)\eqo  \epp \iint_{A'}  \psi_{[\epp]} (x')\psi_{[\epp]} (y') \mathrm dx'  \mathrm dy'$. Similarly, since the $\phi_\pm$ are even, we get 
$\Lambda_{\phi_\epp, F} (v)\eqo  \epp \iint_{A'}  \phi_{\epp} (x')\phi_{\epp} (y') \mathrm dx'  \mathrm dy'$. Therefore 
$$\Lambda_{\psi , F} (v)\! -\!  \Lambda_{\phi_\epp, F} (v)= \epp  \iint_{A'}  \big(  (\psi_{[\epp]}(x')\! -\! \phi_{\epp}(x')  )\psi_{[\epp]}(y') + \phi_{\epp} (x') (\psi_{[\epp]}(y')\! -\! \phi_{\epp}(y')  ) \big) \mathrm dx'  \mathrm dy'  .$$
For all $z\ino \bbR$, we see that $\psi_{[\epp]} (z) \eqo \phi_{-\epp \mathtt s (z) } (z)$, since $\phi_+$ and $\phi_-$ are even and since $\phi_{+} (0)\eqo \phi_+ (0) \eqo \psi (0)$. Thus 
$\psi_{[\epp]}(z)\! -\! \phi_{\epp}(z) \eqo -\epp (\phi_{+} (z) \! -\! \phi_- (z))\un_{\{ z\geq 0\}}$. Consequently 
\begin{eqnarray*} \Lambda_{\phi_\epp, F} (v)  \!\!\! \! & - &\! \!\!\! \Lambda_{\psi , F} (v) \\
 \!\! \!\! & =&  \!\!\!\!    \iint_{A'}  \!\!\!  \big( (\phi_{+}(x')\! -\! \phi_{-}(x')  )\psi_{[\epp]}(y')   \un_{\{ x'\geq 0 \}}+ \phi_{\epp} (x') (\phi_{+}(y')\! -\! \phi_{-}(y')  ) \un_{\{ y'\geq 0 \}} \big) \mathrm dx'  \mathrm dy'.
 \end{eqnarray*}
By Assumption $(a)$ and $(b)$ we then get 
\begin{eqnarray*}
(\phi_{+}(x')\! -\! \phi_{-}(x')  )\psi_{[\epp]}(y')\un_{\{ x'\geq 0 \}} \!\!\!\! &+& \!\!\!\! \phi_{\epp} (x') (\phi_{+}(y')\! -\! \phi_{-}(y')  ) \un_{\{ y'\geq 0 \}} \\
& \leq &  K(x')^2  \psi_{[\epp]}(y') \un_{\{ 0\leq x'< \sigma \}} +  K(y')^2  \phi_{\epp}(x') \un_{\{ 0\leq y'< \sigma \}}  \\
& \leq & K \psi (0) \big( (x')^2 \un_{\{ 0\leq x'< \sigma \}}+(y')^2 \un_{\{ 0\leq y'< \sigma \}}  \big).
\end{eqnarray*}
Thanks to the change of variables $(x,y)\eqo (x'+\epp v, y'+ \epp v)$, we obtain
$$  \Lambda_{\phi_\epp, F} (v)  \! -\! \Lambda_{\psi , F} (v) \leq  \!  K\psi (0)    \iint_{A}  \!\!\!  \big( (x\! -\! \epp v) ^2 
\un_{\{ 0\leq x-\epp v< \sigma \}}+(y\! -\! \epp v)^2 \un_{\{ 0\leq y-\epp v< \sigma \}}  \big) \mathrm dx  \mathrm dy.$$
If $\mathtt s(v) \eqo -\epp$, $(iii)$ shows that $\Lambda_{\phi_\epp, F} (v)  \! =\! \Lambda_{\psi , F} (v) $ 
and (\ref{symmetrization}) holds true trivially. Let us suppose that $\mathtt s(v)\eqo \epp$. Then $\epp v\eqo |v|$ and we get 
\begin{eqnarray*}
 (x\! -\! \epp v) ^2 
\un_{\{ 0\leq x-\epp v< \sigma \}}+(y\! -\! \epp v)^2 \un_{\{ 0\leq y-\epp v< \sigma \}} & = &  (x\! -\! |v| ) ^2 
\un_{\{ 0\leq x-|v|< \sigma \}}+(y\! -\! |v|)^2 \un_{\{ 0\leq y-| v|< \sigma \}}  \\
&\leq & \sigma^{1-\eta} \big( x^{1+\eta} + y^{1+\eta} \big). 
\end{eqnarray*}
By Lemma \ref{Gammabasicprop} $(ii)$ we thus get  $\Lambda_{\phi_\epp, F} (v)  \! -\! \Lambda_{\psi , F} (v) \! \leq  \!  K\psi (0) \sigma^{1-\eta} M^{_{(\eta)}}_{\! F}$, which completes the proof of $(iv)$. \cqfd

\subsection{Lower bounds for $\Lambda_{\phi, F}$ when $\phi$ is a quadratic Beta density}
\label{lowerLambdasec}

In this section we provide uniform lower bounds for $\Lambda_{\phi, F}$ when $\phi$ is an even quadratic Beta density, namely a probability density function of the form discussed in the following lemma. 
\begin{lemma}
\label{quadBeta1}
Let $\tau, \sigma, a\ino (0, \infty)$ such that $\tau \geko \sigma \geko \frac12\tau$. For all $v\ino \bbR$ we set 
\begin{equation}
\label{quadBetadef}
f(v) \eqo \tau^2 \! \! -\! v^2, \quad \phi (v) \eqo af(v)\un_{\{ |v| <\sigma\}} , \quad P(v) \eqo \tfrac34 \big(v\! -\! \tfrac{1}{3} v^3)  \quad \textrm{and} \quad \mathtt b(v)  \eqo  P'(v) \un_{\{|v| < 1\}} .  
\end{equation}
Then, the following holds true. 
\begin{compactenum}

\smallskip

\item[$(i)$] $\mathtt b\ino \mathscr D$  and $\overline{\mathtt b} (v) \eqo \frac12 + \frac34v -\frac14 v^3\eqo 
\frac12 +  P(v)$, for all  $v\ino [-1,1]$.

\smallskip

\item[$(ii)$] $\phi \ino \mathscr D$ if and only if $2a \big(  \tau^2\sigma \! -\! \frac13\sigma^3 \big)\eqo 1$, i.e.~$\frac83 a  \tau^3P(\frac{\sigma}{\tau})\eqo 1$. In this case it implies $\phib(v) \eqo \frac12+ \frac43 a\tau^3 P(\frac{v}{\tau})$ for all $v\ino [-\sigma, \sigma]$.

\smallskip

\item[$(iii)$] Let us assume $\phi \ino \mathscr D$. Then $\phib(v)\leqo a\tau (\tau+v)^2$ for all $v\ino [-\sigma, \sigma]$.

\smallskip

\item[$(iv)$] Let us assume $\phi \ino \mathscr D$. Then 
$\sup_{x\in \bbR} \big| \phib(x\tau)\! -\! \overline{\mathtt b} (x)\big| \leq 2 \big(1\! -\! \frac{\sigma}{\tau} \big)^2$. 
\end{compactenum}
\end{lemma}
\noi
\textbf{Proof.} Here $(i)$ and$(ii)$ follows from straightforward computations: we leave the proofs to the reader.  
Let us prove $(iii)$. We fix $v\ino [-\sigma, \sigma]$. By a change of variable $y\eqo 1+ \frac{w}{\tau}$, 
$$\frac{\phib(v)}{a\tau^3} \eqo \int_{-\sigma}^v \!\!\! \Big(1\! -\! \frac{w^2}{\tau^2} \Big) \frac{\mathrm dw}{\tau} 
\eqo \int_{1-\frac{\sigma}{\tau}}^{1+ \frac{v}{\tau}}\!\! \!\! \!\! \!\!  \!\!   (2\! -\! y)y \, \mathrm dy = \!  \int_{0}^{1+ \frac{v}{\tau}}\!\! \!\!  \!\!  \!\!  \!\!   (2\! -\! y)y \, \mathrm dy -\!  \int_0^{1-\frac{\sigma}{\tau}} \!\! \!\!  \!\!  \!\!  \!\!    (2\! -\! y)y \, \mathrm dy \; .$$
Since $|v| \leqo \sigma$, we have $0\leqo 1 \! -\! \frac{\sigma}{\tau} \leqo 1+ \frac{v}{\tau} \leqo 2$. Thus $\tfrac{1}{a\tau^3} \phib(v)\leqo \int_{0}^{1+ \frac{v}{\tau}}\!  2y \, \mathrm dy$, which implies $(iii)$. 

Let us prove $(iv)$. We note that both $\phib \! -\! \frac12$ and $\overline{\mathtt b}\! -\! \frac12$ are odd. So we only need to get an upper bound for $\sup_{x\in \bbR_+} \big| \phib(x\tau)\! -\! \overline{\mathtt b} (x)\big|$. 
If $x \geqo 1$, we note that $\phib(x\tau)\eqo \overline{\mathtt b} (x)\eqo 1$.  
If $\frac{\sigma}{\tau} \leqo x \leko 1$, then $\phib(x\tau)\eqo 1$, so by a change of variable $y\eqo 1\! -\! z$, 
$$0 \leqo \phib(x\tau)\! -\! \overline{\mathtt b} (x) \leqo 1 \! -\! \overline{\mathtt b} \big( \tfrac{\sigma}{\tau} \big) \eqo
 \int_{\frac{\sigma}{\tau}}^1 \!\! P'(y)\, \mathrm dy = \!  \tfrac{3}{4}\!    \int_0^{1-\frac{\sigma}{\tau}} \!\!\!\!\!\! \!\!\! (2\! -\! z)z \, \mathrm dz\leq \,   \tfrac{3}{4}\!    \int_0^{1-\frac{\sigma}{\tau}} \!\!\!\! \!\!\!\! 2 z \, \mathrm dz =\!   \tfrac{3}{4} \big(1\! -\! \tfrac{\sigma}{\tau} \big)^2.$$
If $x\ino [0, \frac{\sigma}{\tau} ]$, $(i)$ and $(ii)$ imply 
$\phib(x\tau) \! -\! \overline{\mathtt b} (x) \eqo (\frac43a\tau^3 \! -\! 
1) P(x)$. Note that $\max_{x\in [0, 1]} P(x)\eqo P(1) \eqo \frac12$. 
Thus $|\phib(x\tau) \! -\! \overline{\mathtt b} (x)| \leqo \frac12 | \frac43 a\tau^3 \! -\! 
1|$. We next observe that 
$$ 0\leqo \tfrac43 a\tau^3 \! -\! 1 
\eqo  \frac{1}{2P\big( \tfrac{\sigma}{\tau} \big)} -\frac{1}{2P(1)}\eqo \frac{P(1)\!- \! P\big( \tfrac{\sigma}{\tau} \big)}{2P(1)P(\tfrac{\sigma}{\tau})} 
\leq \tfrac{32}{11}\!  \int_{\frac{\sigma}{\tau}}^1 \!\! P'(y)\, \mathrm dy \leq  
  \tfrac{24}{11} \big(1\! -\! \tfrac{\sigma}{\tau} \big)^2.$$
since $P(\frac{\sigma}{\tau}) \geqo P(\frac12)\eqo \frac{11}{32}$ (because $\frac{\sigma}{\tau} \geqo \frac12$ and $P$ increases on $[0, 1]$) and $P'(y)\eqo \frac34 (1\! -\! y^2)\leqo \frac32 (1\! -\! y)$, for $y\ino [0, 1]$. This completes the proof of $(iv)$. 
\cqfd 

\medskip

The following lemma provides key bounds on $\Lambda_{\phi, F}$ when $\phi$ is a quadratic Beta density as in the previous lemma. It is the main technical point of the article. 
\begin{lemma}
\label{Lambdabounddeter} Let $\tau, \sigma, a\ino (0, \infty)$ be such that $\tau \geko \sigma \geko \frac12\tau\geko 2$. We set $D\eqo \tau \! -\! \sigma$. 
We define $f$ and $\phi$ as in (\ref{quadBetadef}) and we assume that $\phi \ino \mathscr D$, i.e.,~$2a(\sigma \tau^2 \! -\! \frac13 \sigma^3)\eqo 1$. Let $F\ino \mathscr H_+$ (i.e.~$\epp (F)\eqo 1$). 
We set $G\eqo G_F$ as in Proposition \ref{corrprop} $(i)$. 
Let $\eta\ino (0, 1)$. We recall from (\ref{Gamma}) the definition of $\mathtt c(F)$ and $M^{_{(\eta)}}_{\! F}$ and we assume that $ M^{_{(\eta)}}_{\! F} \leko \infty$, which implies that $\mathtt c(F) \leko \infty$ by Lemma \ref{Gammabasicprop} $(viii)$. We then set 
\begin{equation}
\label{Sfunc}
\forall v\ino \bbR, \quad S(v) \eqo f(v)^2 \Gamma^{0,1}_{\! F} + 2\mathtt c(F) vf(v) . 
\end{equation}
Then the following holds true. 
\begin{compactenum}

\smallskip

\item[$(i)$] $\frac12 D \! -\! \sigma \leko \sigma$ and for all $v \ino [\frac12 D \! -\! \sigma, \sigma]$, 
$ \big| a^{-2} \Lambda_{\phi, F} (v) \! -\! S(v)\big| \leq 69 M^{_{(\eta)}}_{\! F}\tau^2 (\sigma+ v)^{1-\eta} $.

\smallskip

\item[$(ii)$] For all $v\ino [-\sigma, \sigma]$, $  a^{-2} \Lambda_{\phi, F} (v) \leq S(v)+ 99M^{_{(\eta)}}_{\! F}\tau^2 (\tau+ v)^{1-\eta}$.  

\smallskip

\item[$(iii)$] Let $\gamma \ino (0, 1)$. We furthermore assume $D^2 \geko \tau$. Then 
$$ \inf_{v\in [\sigma, \sigma+ \gamma]}   a^{-2} \Lambda_{\phi, F} (v) \geq \tfrac{1}{4} \tau^2 D^2 \Big( \int_{\bbR}\big( G(z) \! -\! \gamma\big)_{\! +} \mathrm d z -73 M^{_{(\eta)}}_{\! F}\tau^{-\eta} \Big). $$

\end{compactenum}
\end{lemma}
\noi
\textbf{Proof.} Since $\tau \leko 2\sigma$, $\frac12 D\eqo \frac12(\tau \! -\! \sigma) \leko \frac12\sigma$ and $\frac12 D \! -\! \sigma \leko-\frac12 \sigma \leko 0$. We first prove $(i)$. We recall the notation 
$A \eqo \big\{ (x,y)\ino \bbR_+^2: G(y\! -\! x) \geq \min (x,y) \big\} $ and recall from Lemma \ref{propLambdaeven} $(ii)$ that $\Lambda_{\phi, F} (  v )\eqo \iint_A \phi (x\! -\! v) \phi(y\! -\! v) \, \mathrm dx \mathrm dy$ since $\epp (F)\eqo 1$. For all $v\ino \bbR$, we set 
\begin{equation}
\label{Bvdef}
B_v = \big\{(x,y)\ino \bbR_+^2 : (v\! -\! \sigma)_+ \leqo \min (x,y) \leqo \max (x,y) \leqo \sigma + v \big\} \; .
\end{equation} 
We first observe that $a^{-2}\phi (x\! -\! v) \phi (y\! -\! v) \eqo f (x\! -\! v) f (y\! -\! v)\un_{B_v} (x,y)$ for Lebesgue almost all $x,y\ino \bbR_+$ and all $v\ino \bbR$.  Thus we get 
$$a^{-2}\Lambda_{\phi, F} (v) \eqo \iint_{\! A\cap B_v} \!\!\!\! \!\! \!\!\!\!   f (x\! -\! v) f (y\! -\! v) \, \mathrm dx \mathrm dy . $$
We check that for all $x,y\ino \bbR_+$ and all $v\ino \bbR$, 
$$ f (x\! -\! v) f (y\! -\! v) \eqo f(v)^2\! + 2 (x+y) vf(v)+ 4xyv^2 \! -\! x^2(f(v)+ 2yv)\! -\! y^2(f(v)+2xv)+ x^2y^2. $$
We next suppose that $v\ino [-\sigma , \tau]$. Then $|v| \leqo \tau$, $|f(v)|\leqo \tau^2$ and 
\begin{eqnarray}
\label{squaboun}
LHS_{\eqref{squaboun}} \!\!\!\! & :=&\!\!\!\!  \big| f (x\! -\! v) f (y\! -\! v) \! -\!  f(v)^2\! -\!  2 (x+y) vf(v) \big|  \\
\!\!\!\!  &\leqo & \!\!\!\! 4\tau^2xy+ x^2 (\tau^2+2y\tau) + y^2 (\tau^2 + 2x\tau) + x^2y^2.  \nonumber
\end{eqnarray}
Let us furthermore assume that $(x,y)\ino A\cap B_v$. Then 
$0\leqo \min(x,y)\leqo \max(x,y) \leqo \sigma+v\leqo 2\tau$ and we get 
$$LHS_{\eqref{squaboun}}\leq 8\tau^2 xy + 5\tau^2 (x^2 \! + y^2) \leqo 9 \tau^2 (x^2\! + y^2) $$
because $ xy \leqo \frac12 (x^2+y^2)$. By Lemma \ref{Gammabasicprop} $(ii)$
we get 
$$ \iint_{A\cap B_v} \!\!\!\!\!\!\!\!\!\! (x^2\! + y^2) \, \mathrm dx \mathrm dy \leq (\sigma+v)^{1-\eta} \!\!  \iint_{A\cap B_v}  \!\!\!\!\!\!\!\!\!\!(x^{1+\eta}\! + y^{1+\eta}) \, \mathrm dx \mathrm dy 
\leq 
(\sigma+v)^{1-\eta} M_{\! F}^{_{(\eta)}}.$$
We set 
$$\forall v\ino \bbR, \quad S_0(v)\eqo f(v)^2 \! \iint_{\! A\cap B_v} \!\!\!\!\!\!\!\!\!  \mathrm dx \mathrm dy\, + 2vf(v) \! \iint_{\! A\cap B_v} \!\!\!\!\!\!\!\!\!  (x + y) \, \mathrm dx \mathrm dy\;, 
$$ and we have proved 
\begin{equation}
\label{S0bound}
\forall v\ino [-\sigma, \tau], \quad \big| a^{-2} \Lambda_{\phi, F} (v) \! -\! S_0(v) \big| \leq 9   M_{\! F}^{_{(\eta)}}\tau^2 (\sigma+v)^{1-\eta}. 
\end{equation}

We now suppose that $v\ino [ \frac12 D \! - \! \sigma, \sigma]$ as in $(i)$. We set 
\begin{equation}
\label{Bvprimedef}
B'_v= \big\{ (x,y)\ino \bbR_+^2: \max (x,y) \geko v+ \sigma \big\} \; .
\end{equation}
We recall $B_v$ from (\ref{Bvdef}) and since we assume $v\leqo \sigma$, we observe that 
$A\cap B^c_v\eqo A\cap B'_v$. We recall $S(v)$ from (\ref{Sfunc}). 
By Lemma \ref{Gammabasicprop} $(ii)$, $S(v)\eqo f(v)^2\! \iint_A \! \mathrm dx  \mathrm dy+ 2vf(v) \iint_A \! (x+y) \, \mathrm dx  \mathrm dy$.  
Thus $S(v) - S_0(v)\eqo f(v)^2 \iint_{A\cap B'_v}  \! \mathrm dx \mathrm dy \, + 2vf(v) \iint_{A\cap B'_v} \! (x + y) \, \mathrm dx \mathrm dy$. By Lemma \ref{Gammabasicprop} $(iv)$, 
\begin{eqnarray*}
 f(v)^2 \! \iint_{A\cap B'_v}  \!\!\!\!\!\!\!\!\!\!  \mathrm dx \mathrm dy \leq 
(\tau\! -\! v)^2 \! \frac{(\tau+v)^2}{(\sigma+v)^{1+\eta}} \!\!\!\! \!\!\!\! & & \!\!\!\! \!\!\!\!   \iint_A \!\!\!  \max (x,y)^{1+\eta}  \mathrm dx \mathrm dy \leq (\tau\! -\! v)^2 \! \frac{(\tau+v)^2}{(\sigma+v)^{1+\eta}} M_{\! F}^{_{(\eta)}} \\
 \textrm{and}  \quad 2|v| f(v)  \iint_{A\cap B'_v}  \!\!\!\!\!\!\!\!\!\!  (x+y) \, \mathrm dx \mathrm dy \!\!\! & \leq & \!\!\!   2(\tau\! -\! v) |v| \frac{\tau+v}{(\sigma+v)^\eta}  \iint_A \!\!\!  \max (x,y)^{\eta} (x+y) \mathrm dx \mathrm dy \\
 \!\!\! & \leq & \!\!\!  4(\tau\! -\! v) \tau \tfrac{\tau+v}{(\sigma+v)^\eta} M_{\! F}^{_{(\eta)}}.
\end{eqnarray*}
Since we assume that $v\ino [ \frac12 D \! - \! \sigma, \sigma]$, we get $|\tau \! -\! v| \leqo 2\tau$ and $\frac{\tau+v}{\sigma+v}= 1+ \frac{D}{\sigma+v} \leqo 3$. Therefore for all $v\ino [ \frac12 D \! - \! \sigma, \sigma]$,
$$ \big| S(v)\! -\! S_0(v)\big| \leq 4\tau^2. 9. (\sigma+v)^{1-\eta} M_{\! F}^{_{(\eta)}}+ 4.2\tau^2.3.(\sigma+v)^{1-\eta} M_{\! F}^{_{(\eta)}}=60M_{\! F}^{_{(\eta)}}\tau^2(\sigma+v)^{1-\eta},$$
which entails $(i)$ thanks to (\ref{S0bound}).

\smallskip

We next prove $(ii)$. Thanks to $(i)$, we only need to prove $(ii)$ for $v\ino [-\sigma, -\sigma + \frac12 D)$. 
There are two cases to consider: \emph{Case 1}: $G(0) \leko \frac12 D$ and  \emph{Case 2}: $G(0) \geqo \frac12 D$. 

\smallskip

\noi
$\bullet$  \emph{Case 1}: we suppose $G(0) \leko \frac12 D$ and we fix $v\ino [-\sigma, -\sigma + \frac12 D)$. 
By (\ref{S0bound}), 
\begin{equation}
\label{upp1}
a^{-2} \Lambda_{\phi, F} (v) \leq S_0(v) + 9   M_{\! F}^{_{(\eta)}}\tau^2 (\sigma+v)^{1-\eta}. 
\end{equation}
To simplify notation we set 
$$ w\eqo \sigma+v \in [0, \tfrac12 D], \quad \alpha_w\eqo -v\eqo \sigma-w \quad \textrm{and} \quad \gamma_w \eqo D+w \; .$$
We observe that $f(v)\eqo (\tau+ \sigma \! -\! w ) (\tau \! -\! \sigma +w)\eqo (2\alpha_w+ \gamma_w) \gamma_w$. Therefore, 
\begin{eqnarray*}
S_0 (v)\!\!\! \!\!  &= & \!\! \!\!\!  f(v)\!  \iint_{A\cap B_v} \!\!\!\!\!\!\! \!\!\! (f(v)\! + 2(x+y)v)  \mathrm dx \mathrm dy \\
\!\! \!\!\! \!\! =  \;\,  & &\!  \! \!\! \!\! \!\!\!\!\! \!\!\!\!\! \!\! (2\alpha_w\! +\gamma_w)\gamma_w \!\!   \iint_{A\cap B_v} \!\!\!\!\!\!\! \!\!\! \! \big( (2\alpha_w+\gamma_w)\gamma_w  \! -\!  2\alpha_w (x+y) \big)  \mathrm dx \mathrm dy\\
\!\! \!\!\! \!\! =  \;\,  & &\!  \! \!\! \!\! \!\!\!\!\! \!\!\!\!\! \!\! (2\alpha_w\! +\! \gamma_w)\gamma^3_w \!\!   \iint_{A\cap B_v} \!\!\!\!\!\!\! \!\!\! \!  \mathrm dx \mathrm dy \; + 2\alpha_w (2\alpha_w\! +\! \gamma_w)\gamma_w \!\! \iint_{A\cap B_v} \!\!\!\!\!\!\! \!\!\! (\gamma_w\! -\!  (x+y))  \mathrm dx \mathrm dy \\
\!\! \!\!\! \!\! =  \;\,  & &\!  \! \!\! \!\! \!\!\!\!\! \!\!\!\!\! \!\! (2\alpha_w\! +\! \gamma_w)\gamma^3_w \!\!   \iint_{A\cap B_v} \!\!\!\!\!\!\! \!\!\! \!  \mathrm dx \mathrm dy \; + 2\alpha_w (2\alpha_w\! +\! \gamma_w)\gamma_w\!\!  \iint_{A\cap B_v} \!\!\!\!\!\!\! \!\!\! \!\! 
\big( D \! -\! \min (x,y) + w\! -\! \max (x,y) \big)  \mathrm dx \mathrm dy.
\end{eqnarray*}
We set $C_w\eqo A\cap B_v \eqo \{ (x,y)\ino A: \max (x,y)\leqo w\}$. 
Let $w'\ino [w, \frac12 D]$. For $(x,y)\ino C_w \! \subset \! C_{w'}$, we have 
$\max (x,y)\leqo w'$ and since $\min (x,y) \leqo G(y\! -\! x) \leqo G(0) \leko \frac12 D$, 
we get $D \! -\! \min (x,y) + w\! -\! \max (x,y)\geqo D\! -\! G(0) + w\! -\! w'\geqo 0$. Since $\alpha_w\geqo 0$, $\gamma_w\geqo 0$, and thus $2\alpha_w+ \gamma_w \geqo 0$, we get 
\begin{eqnarray*}
S_0 (v)\!\!\! & \leq & \!\!\!  (2\alpha_w\! +\! \gamma_w)\gamma^3_w \!\! \iint_{C_{w'}} \!\!\!\!\!\!\! \mathrm dx \mathrm dy \\
\!\!\!\!& &\qquad \qquad \; + \; \,   2\alpha_w (2\alpha_w\! +\! \gamma_w)\gamma_w\!\!  \iint_{C_{w'}} \!\!\!\!\!\!\! 
\big( D \! -\! \min (x,y) + w\! -\! \max (x,y) \big)  \mathrm dx \mathrm dy \\
&= & f(v)^2\!  \iint_{C_{w'}} \!\!\!\!\! \mathrm dx \mathrm dy \, +\,  2vf(v)\!  \iint_{C_{w'}} \!\!\!\!\! (x+y) \, \mathrm dx \mathrm dy. 
\end{eqnarray*}
We recall that $C_{\! \frac12 D}\eqo A\cap B_{\! \frac12 D-\sigma}$ and that $S(v)\eqo f(v)^2\! \iint_A \! \mathrm dx  \mathrm dy+ 2vf(v) \iint_A \! (x+y) \, \mathrm dx  \mathrm dy$. 
So we apply the previous equality to $w'\eqo \frac12 D$ and get 
\begin{eqnarray}
S_0 (v) \!\!\!\! &\leq &\!\!\!\!  f(v)^2 \!\! \iint_{C_{\! \frac12 D}} \!\!\!\!\!\!\! \mathrm dx \mathrm dy\, + 2vf(v)\!  \iint_{C_{\! \frac12 D}} \!\!\!\!\!\!\! (x+y) \, \mathrm dx \mathrm dy  \nonumber \\
\!\!\!\! & = & \!\!\!\! S(v)- f(v)^2 \!\! \iint_{A\backslash C_{\! \frac12 D}} \!\!\!\!\!\!\! \!\!\!\!\!\!\! \mathrm dx \mathrm dy - 2vf(v) \!\!  \iint_{A\backslash C_{\! \frac12 D}} \!\!\!\!\!\!\!  \!\!\!\!\!\!\! (x+y) \, \mathrm dx \mathrm dy \nonumber \\
\!\!\!\! & \leq & \!\!\!\! S(v) + 2 |v| f(v) \!\!  \iint_{A\backslash C_{\! \frac12 D}} \!\!\!\!\!\!\! \!\!\!\!\!\!\!  (x+y) \, \mathrm dx \mathrm dy \label{S0uppbou}
\end{eqnarray}
since $-v\eqo |v|$ (because $v\ino [-\sigma, -\sigma + \frac12 D] \! \subset \! [-\sigma, 0] $) and since $f(v)\geko 0$. Moreover we observe that 
$2|v| f(v)\leq 4\tau^2(\tau+v)$ and $(\tau+v)/(\frac12 D) \leqo 3$. 
We recall from (\ref{Bvprimedef}) the definition of $B'_{\! \frac12 D-\sigma}$ and observe that $A\backslash C_{\! \frac12 D}\eqo A\cap B'_{\! \frac12 D-\sigma}$. As argued in the proof of $(i)$ and by Lemma \ref{Gammabasicprop} $(iv)$, we obtain the following    
\begin{eqnarray*}
 2|v|f(v)\!  \iint_{A\cap B'_{\! \frac12 D-\sigma}}  \!\!\!\!\!\!\! \!\!\!\!\!\!\! \!\!\!\!\! (x+y) \, \mathrm dx \mathrm dy &\leq &4\tau^2(\tau+v) \big( \tfrac{1}{2}D)^{-\eta} \!\!  \iint_{A} \!\!\! \max(x,y)^\eta (x+y) \, \mathrm dx \mathrm dy \\
&\leq&  8 M^{_{(\eta)}}_{\! F} \tau^2(\tau+v) \big( \tfrac{1}{2}D)^{-\eta} \leq 8 \! \cdot \! 3^\eta  M^{_{(\eta)}}_{\! F}\tau^2 (\tau+v)^{1-\eta}, 
\end{eqnarray*}
since $\frac{\tau+v}{\frac12 D}\leqo 3$. Therefore (\ref{S0uppbou}) entails $S_0 (v)\leqo S(v) + 24 M^{_{(\eta)}}_{\! F} \tau^2(\tau+v)^{1-\eta}$. Combined with (\ref{upp1}), it follows that if $G(0) \leko \frac12 D$, 
\begin{equation}
\label{G0moinsdemiD}
\forall v\ino [-\sigma, -\sigma +\tfrac{1}{2} D) , \quad  a^{-2}\Lambda_{\phi, F} (v) \leq S(v) +  33 M^{_{(\eta)}}_{\! F}\tau^2 (\tau+v)^{1-\eta}.
\end{equation}

\smallskip

\noi
$\bullet$  \emph{Case 2}: we suppose $G(0) \geqo \frac12 D$ and we fix $v\ino [-\sigma, -\sigma + \frac12 D)$. 
By Lemma \ref{elempropLambda} $(ii)$ $0\leqo \Lambda_{\phi, F} (v) \leqo \phib (v)^2$ and by Lemma \ref{quadBeta1} $(iii)$ we get $0\leqo a^{-2}\Lambda_{\phi, F} (v)\leqo \tau^2 (\tau+v)^4$. 
On the other hand, by 
Lemma \ref{Gammabasicprop} $(vii)$ 
we have $G(0)^{3+\eta} \leqo \Gamma_{\! F}^{0, 2+ \eta} \leqo M^{_{(\eta)}}_{\! F}$. Since $(\tau+v)/(\frac12 D)\leqo 3$, it follows that 
\begin{equation}
\label{G0plusdemiD}
a^{-2} \Lambda_{\phi, F} (v)\un_{\{G(0) \geq \frac12 D \}} \leq  \tau^2 (\tau+v)^4 \frac{G(0)^{3+\eta}}{(\tfrac{1}{2}D)^{3+\eta}} \leqo 81    M^{_{(\eta)}}_{\! F}\tau^2 (\tau+v)^{1-\eta} .
\end{equation}

\noi
\emph{End of the proof of $(ii)$.} We first note that (\ref{G0moinsdemiD}) and (\ref{G0plusdemiD}) imply 
\begin{equation}
\label{sparadrap}
\forall v\ino [-\sigma, -\sigma +\tfrac{1}{2} D) , \quad \; a^{-2}\Lambda_{\phi, F} (v) \leq S(v)\un_{\{ G(0) <\frac12 D\}} +  81 M^{_{(\eta)}}_{\! F}\tau^2 (\tau+v)^{1-\eta}.
\end{equation}
We observe that 
\begin{eqnarray}
S(v)\un_{\{ G(0) <\frac12 D\}} \!\!\! &=&  \!\!\! S(v)\! -\! f(v)^2 \Gamma^{0,1}_{\! F} \un_{\{ G(0) \geq \frac12 D\}}\! -\! 2vf(v) \mathtt c(F) \un_{\{ G(0) \geq \frac{1}{2} D\}} \nonumber \\
\!\!\!  &\leq &\!\!\!  S(v) + 2|v| f(v) \mathtt c(F)  \un_{\{ G(0) \geq \frac12 D\}} .
\label{ssppaa}
\end{eqnarray}
We recall that since $ v\ino [-\sigma, -\sigma +\tfrac{1}{2} D)$, we get $ 2|v| f(v)\leqo 4\tau^2(\tau+v)$ and 
$ (\tau+v)/(\frac12 D)\leqo 3$. By Lemma \ref{Gammabasicprop} $(ix)$ $\mathtt c(F) \eqo \Gamma^{1,1}_{\! F}+ \frac12 \Gamma^{0,2}_{\! F}\leq \tfrac32 M^{_{(\eta)}}_{\! F} G(0)^{-\eta}$. Therefore, 
\begin{eqnarray*}
 2|v| f(v) \mathtt c(F)  \un_{\{ G(0) \geq \frac12 D\}} \!\!\! & \leq &\!\!\! 4 \tau^2 (\tau+v) \! \cdot \!  \tfrac32   M^{_{(\eta)}}_{\! F} G(0)^{-\eta}\un_{\{ G(0) \geq \frac12 D\}} \\  
 \!\!\! & \leq & \!\!\! 6 M^{_{(\eta)}}_{\! F} \tau^2 (\tau+v)^{1-\eta} 
 \frac{(\tau+v)^{\eta}}{(\tfrac{1}{2}D)^\eta} \leq  18 M^{_{(\eta)}}_{\! F} \tau^2 (\tau+v)^{1-\eta} , 
\end{eqnarray*}
which completes the proof of $(ii)$ thanks to (\ref{sparadrap}) and (\ref{ssppaa}). 

\smallskip

It remains to prove $(iii)$. We assume that $D^2 \geko \tau$ and that $\frac12 \tau \geko 2$. Then $\frac12 D \geko 1\geko \gamma$. 
We fix $v\ino [\sigma, \sigma + \gamma]$. Then  $\sigma + \gamma \leqo \tau$, (\ref{S0bound}) applies and we get 
\begin{equation}
\label{minoexcess}
a^{-2} \Lambda_{\phi, F} (v)\geq S_0(v) -9 M^{_{(\eta)}}_{\! F} \tau^2 (\sigma+v)^{1-\eta} \geq f(v)^2 \!\!  \iint_{A\cap B_v} \!\!\! \!\!\! \!\!\! \! \mathrm dx \mathrm dy - 18M^{_{(\eta)}}_{\! F} \tau^{3-\eta} 
\end{equation}
since $2vf(v) \iint_{A\cap B_v}\! (x+y)  \mathrm dx \mathrm dy \geqo 0$ and since $\sigma+v\leqo 2\tau $. 

We next set $w\eqo v\! -\! \sigma \ino [0, \gamma]$. Then $f(v)^2\eqo (\tau\! -\! \sigma \! -\! w)^2 (\tau + \sigma + w)^2 \geqo (D\! -\! w)^2 \tau^2 \geqo \frac14 D^2\tau^2$ since $\gamma \leko 1\leko  \frac12 D$. 
We recall $B_v$ from (\ref{Bvdef}) and $B'_v$ from (\ref{Bvprimedef}), and observe that 
$$A\, \cap \, B_v \eqo \big\{ (x,y)\ino A: \min (x,y) \geqo w \big\}
\cap \big\{ (x,y)\ino A: \max (x,y) \leqo v+ \sigma \big\} \! = \!   (A\, \cap \,  [w, \infty)^2) \backslash (A\cap B'_{v}). $$
Thus 
$\iint_{A\cap B_v} \!  \mathrm dx \mathrm dy \geq   \iint_{A\cap [w, \infty)^2} \!\mathrm dx \mathrm dy-\iint_{A\cap B'_v} \! \mathrm dx \mathrm dy$. As proved in $(i)$, 
$$\iint_{A\cap B'_v} \!\!\!\!\!\!\!\! \!  \mathrm dx \mathrm dy \leq (\sigma+v)^{-1-\eta} \!\! \iint_A \! \max(x,y)^{1+\eta}  \mathrm dx \mathrm dy \leq M^{_{(\eta)}}_{\! F}  (\sigma+v)^{-1-\eta} \leq 4   M^{_{(\eta)}}_{\! F}  \tau^{-1-\eta}, $$
since we assume that $\sigma \geko \frac12 \tau$. On the other hand, 
by Lemma \ref{Gammabasicprop} $(v)$, $\iint_{A\cap [w, \infty)^2}  \mathrm dx \mathrm d y\eqo \int_\bbR \big( G(z) \! -\! w)_+ \mathrm dz \geqo \int_\bbR \big( G(z) \! -\! \gamma)_+ \mathrm dz$. Thus 
$$ f(v)^2 \!\!  \iint_{A\cap B_v} \!\!\! \!\!\! \!\!\! \! \mathrm dx \mathrm dy \geq \tfrac{1}{4} \tau^2 D^2 \Big(  \int_\bbR \big( G(z) \! -\! \gamma)_+ \mathrm dz - 4   M^{_{(\eta)}}_{\! F}  \tau^{-1-\eta}\Big) $$
and by (\ref{minoexcess}) 
\begin{eqnarray*}  
a^{-2} \Lambda_{\phi, F} (v) \!\!\! &\geq &\!\!\!  \tfrac{1}{4} \tau^2 D^2 \Big(  \int_\bbR \big( G(z) \! -\! \gamma)_+ \mathrm dz - 4   M^{_{(\eta)}}_{\! F}  \tau^{-1-\eta} \! -\!  72   M^{_{(\eta)}}_{\! F}  \tau^{1-\eta} D^{-2} \Big) \\
 \!\!\! &\geq &\!\!\!  \tfrac{1}{4} \tau^2 D^2 \Big(  \int_\bbR \big( G(z) \! -\! \gamma)_+ \mathrm dz - 73  M^{_{(\eta)}}_{\! F}  \tau^{-\eta}  \Big)
\end{eqnarray*}
because we assume $D^2\geko \tau \geko 4$. This completes the proof of $(iii)$. \cqfd

\medskip

We apply the previous lemma to an $\mathscr H$-valued random function $\ff$ satisfying the assumptions of Theorem \ref{t:main}. 
\begin{lemma}
\label{mainLambdaesti} Let $\tau, \sigma, a\ino (0, \infty)$ be such that $\tau \geko \sigma \geko \frac12\tau\geko 2$. We set $D\eqo \tau \! -\! \sigma$. 
We define 
$\phi$ as in (\ref{quadBetadef}) and assume that $\phi \ino \mathscr D$, i.e.~$2a(\sigma \tau^2 \! -\! \frac13 \sigma^3)\eqo 1$. Let $\ff$ be an $\mathscr H$-valued random function.  Let $\eta\ino (0, 1)$. We recall from (\ref{Gamma}) the definition of $\mathtt c(\ff)$ and $M^{_{(\eta)}}_{\! \ff}$. We assume (\ref{maintheohyp}), namely, 
\begin{equation}
\label{maintheohypbis}
\p \big(\mathbf F \! \notin \! \{ \varphi_+, \varphi_-\} \big)\geko 0, \quad K_\eta\! :=\! \e \big[ M^{(\eta)}_{\! \mathbf F}\big]\leko \infty \quad  \textrm{and} \quad \e[\sgn (\ff) ] \eqo \e \big[ \Gamma_{\! \ff}^{0,1} \sgn (\ff) \big]\eqo 0 ,
\end{equation}
which implies that $\e [\mathtt c(\ff)]  \leko \infty$ by Lemma \ref{Gammabasicprop}  $(viii)$. Then the following holds true. 
\begin{compactenum}

\smallskip

\item[$(i)$] $\frac12 D \! -\! \sigma \leko \sigma$ and for all $v \ino [\frac12 D \! -\! \sigma, \sigma]$, 
$$ a^{-2} \e[\Lambda_{\phi, \ff} (v) ]  \geq 2\e[\mathtt c(\ff)] v(\tau^2\! -\! v^2) -198  K_\eta \tau^{3-\eta}\; .$$

\smallskip

\item[$(ii)$] For all $v\leqo  \frac12 D \! -\!  \sigma$,  $  \; a^{-2} \e[ \Lambda_{\phi, \ff } (v)]  \geqo 
-198  K_\eta \tau^{3-\eta}- 15 \e \big[ \Gamma^{0,1}_{\! \ff} \! +\mathtt c(\ff) ]  \tau^2D^2$.  

\smallskip

\item[$(iii)$] If $D^2 \geko \tau$, then $\e\big[ \Gamma^{0,1}_{\! \ff} \un_{\{\epp (\ff)=1 \}}\big] \geko 0$ and 
there are $\gamma_0 \ino (0, 1)$ and $\tau_0\geko 4$, which only depend on $\eta$ and on the law of $\ff$, such that for all $\tau \geko \tau_0$, 
\begin{equation}
\label{lowerexcess} \inf_{v\in [\sigma, \sigma+ \gamma_0]}  \!\!\!\! a^{-2} \e \big[ \Lambda_{\phi, \ff}  (v)\big]  \geq \tfrac{1}{8} \tau^2 D^2 \e\big[ \Gamma^{0,1}_{\! \ff} \un_{\{\epp (\ff)=1 \}}\big]. 
\end{equation}
\end{compactenum}
\end{lemma}
\noi
\textbf{Proof.} Let us recall that $f(v)\eqo \tau^2 \! -\! v^2$, $v\ino \bbR$. We first prove $(i)$. By Lemma \ref{propLambdaeven} $(ii)$, for all $v\ino \bbR$, 
$\Lambda_{\phi, \ff} (v)\eqo \epp (\ff) \Lambda_{\phi, \ff^*} (\epp (\ff) v)$. We recall from the definition 
(\ref{TF}) of $T_\ff\eqo T_{\ff^*} $ and from (\ref{Gamma}) 
that $\Gamma^{a,b}_{\! \ff}\eqo \Gamma^{a,b}_{\! \ff^*}$ and thus 
$\mathtt c (\ff)\eqo  \mathtt c (\ff^*)$ and 
$M^{_{(\eta)}}_{\! \ff}\eqo M^{_{(\eta)}}_{\! \ff^*}$. By Lemma \ref{Lambdabounddeter} 
$(i)$ and $(ii)$, for all $v\ino [ \frac12 D \! -\! \sigma , \sigma]$,  
\begin{eqnarray*}
a^{-2}\Lambda_{\phi,  \ff} (v)  \!\! \!\!  & = &   \!\! \!\!  a^{-2}\Lambda_{\phi,  \ff^*} (v)   \un_{\{\epp (\ff)=1 \}}  -a^{-2}\Lambda_{\phi,  \ff^*} (-v)  \un_{\{\epp (\ff)=-1 \}}  \\
   \!\! \!\!  \!\! \!\!  & \geq & \!\! \!\!   \un_{\{\epp (\ff)=1 \}} 
\Big( f(v)^2\Gamma^{0,1}_{\! \ff} + 2 \mathtt c(\ff) vf(v)\Big) \! -\! 69 M^{_{(\eta)}}_{\! \ff} \tau^2(\tau+v)^{1-\eta} \un_{\{\epp (\ff)=1 \}}  \\
 \!\! \!\! \!\! & & \!\! \!\!  - \un_{\{\epp (\ff)=-1 \}}   \Big( f(-v)^2\Gamma^{0,1}_{\! \ff} - 2 \mathtt c(\ff) vf(-v)\Big) -  99 M^{_{(\eta)}}_{\! \ff} \tau^2(\tau-v)^{1-\eta} \un_{\{\epp (\ff)=-1 \}} .
\end{eqnarray*}
We observe that $|\tau \pm v|\leqo 2\tau$ and that $f(-v)\eqo f(v)$. Thus 
$$ a^{-2} \Lambda_{\phi,  \ff} (v) \geq f(v)^2
\epp (\ff) \Gamma^{0,1}_{\! \ff} + 2 \mathtt c(\ff) vf(v)-198 M^{_{(\eta)}}_{\! \ff} \tau^{3-\eta} $$
which implies $(i)$ by taking the expectation since we assume 
$\e \big[ \Gamma_{\! \ff}^{0,1} \sgn (\ff) \big]\eqo 0$.

\smallskip

Let us prove $(ii)$. We recall that $\Lambda_{\phi, \ff} (v)\eqo \epp (\ff) \Lambda_{\phi, \ff^*} (\epp (\ff) v)$. 
Since $\epp (\ff^*)\eqo 1$, Lemma \ref{elempropLambda} $(ii)$ applies and it asserts that 
$\Lambda_{\phi, \ff^*} (v) \geq 0$, that $v\mapsto \Lambda_{\phi, \ff^*} (v)$ vanishes on $(-\infty, -\sigma]$ and that $v\mapsto \Lambda_{\phi, \ff^*} (v)$ is nonincreasing on 
$[\sigma, \infty)$. Thus 
$$ \inf_{v\leq  \frac{_1}{^2} D-\sigma}  \!\! \!\! \! \Lambda_{\phi, \ff} (v) \geq -\un_{\{\epp (\ff)= -1 \}}   \!\! \!\! 
\sup_{v\in [-\sigma,  \frac{_1}{^2} D-\sigma]}  \!\! \!\!  \Lambda_{\phi, \ff^*} (-v)=  -\un_{\{\epp (\ff)= -1 \}}   \!\! \!\! \sup_{w\in [\sigma - \frac12 D, \sigma]} \!\! \!\! \!\! \!\!  \Lambda_{\phi, \ff^*} (w) . $$
By Lemma \ref{Lambdabounddeter} $(ii)$ applied to $\ff^*$, for all $w\in [\sigma \! -\!  \frac12 D, \sigma]$ we get 
$a^{-2}\Lambda_{\phi, \ff^*} (w) \leqo f(w)^2\Gamma^{0,1}_{\! \ff^*} + 2 \mathtt c(\ff^*) wf(w) + 99 M^{_{(\eta)}}_{\! \ff^*} \tau^2 (\tau+w)^{1-\eta}\eqo f(w)^2\Gamma^{0,1}_{\! \ff} + 2 \mathtt c(\ff) wf(w) + 99 M^{_{(\eta)}}_{\! \ff} \tau^2 (\tau+w)^{1-\eta}$, 
since, as already mentioned, $\Gamma^{a,b}_{\! \ff}\eqo \Gamma^{a,b}_{\! \ff^*}$ and thus 
$\mathtt c (\ff)\eqo  \mathtt c (\ff^*)$ and 
$M^{_{(\eta)}}_{\! \ff}\eqo M^{_{(\eta)}}_{\! \ff^*}$. 
Now observe that for all 
$w\in [\sigma \! -\!  \frac12 D, \sigma]$, we get $0\leqo f(w)\eqo (\tau\! -\! w)(\tau+w)\leqo (\tau \! -\! \sigma + \frac{_1}{^2} D)(\tau+ \sigma) \leqo \frac{_3}{^2} D\! \cdot \! 2\tau\eqo 3\tau D$ and $|\tau +w|\leqo 2\tau$. 
Thus 
 \begin{eqnarray*}
 \sup_{w\in [\sigma - \frac12 D, \sigma]}  \!\! \!\! \!  \!\!\! a^{-2}\Lambda_{\phi, \ff^*} (w) 
 &\leq & 9\tau^2 D^2\Gamma^{0,1}_{\! \ff} + 6 \mathtt c(\ff) \tau^2 D + 
 198 M^{_{(\eta)}}_{\! \ff} \tau^{3-\eta} \\
 &\leq &  198 M^{_{(\eta)}}_{\! \ff} \tau^{3-\eta} + 
 15( \Gamma^{0, 1}_{\! \ff} + \mathtt c(\ff) ) \tau^2 D^2  
 \end{eqnarray*}
 and we get 
 $$ a^{-2} \inf_{v\leq  \frac{_1}{^2} D-\sigma}  \!\! \!\! \! \e \big[ \Lambda_{\phi, \ff} (v)\big] \;  \geq 
  -198 \tau^{3-\eta} \e \big[ \un_{\{\epp (\ff)= -1 \}} M^{_{(\eta)}}_{\! \ff}   \big]- 15 \tau^2D^2  \e \big[ \big( \Gamma^{0, 1}_{\! \ff} +  \mathtt c(\ff)\big)\un_{\{\epp (\ff)= -1 \}}  \big]$$
which immediately implies $(ii)$.

\smallskip

To prove $(iii)$ we fix $\gamma_0\ino (0, 1)$, which is specified below, and fix $v\ino [\sigma, \sigma + \gamma_0]$. We recall, by Lemma \ref{propLambdaeven} $(ii)$, that
$\Lambda_{\phi, \ff} (v)\eqo \epp (\ff) \Lambda_{\phi, \ff^*} (\epp (\ff) v)$ and 
by Lemma \ref{elempropLambda} $(ii)$ we get $\Lambda_{\phi, \ff^*} (-v)\eqo 0$. 
Thus $\Lambda_{\phi, \ff} (v)\eqo\un_{\{ \epp (\ff)=1\}} \Lambda_{\phi, \ff^*} (v)$. 
Next we set 
$\mathbf G \eqo G_{\ff}\eqo G_{\ff^*}$ as in Proposition \ref{corrprop} and since $\epp (\ff^*)\eqo 1$, Lemma \ref{Lambdabounddeter} $(iii)$ applies 
and we get 
$$ \inf_{v\in [\sigma, \sigma+ \gamma_0]}  \!\!\! 
 \e \big[ a^{-2} \Lambda_{\phi, \ff} (v) \big] \geq \tfrac{1}{4} \tau^2 D^2 \Big( \e \Big[ 
 \un_{\{ \epp (\ff) =1\}}\!\!  \int_{\bbR}\! \big( \mathbf G (z) \! -\! \gamma_0 \big)_{\! +} \! \mathrm d z\Big] -73 
 \e \big[ \un_{\{ \epp (\ff)= 1\}}M^{_{(\eta)}}_{\! \ff} \big] \tau^{-\eta} \Big). $$

By (\ref{maintheohypbis}), 
$ \e \big[ \un_{\{ \epp (\ff) =1\}} \Gamma_{\! \ff}^{0,1} \big]\eqo 
 \e \big[ \un_{\{ \epp (\ff) =-1\}} \Gamma_{\! \ff}^{0,1} \big]$ and we denote by $c$ this quantity.  Therefore  
$\e \big[\Gamma_{\! \ff}^{0,1} \big]\eqo 2c$ and if $c\eqo 0$, then a.s.~$\Gamma_{\! \ff}^{0,1}\eqo 0$. 
By Lemma \ref{Gammabasicprop} $(ii)$ it entails $\ff \ino \{ \varphi_+, \varphi_-\}$ a.s., which contradicts Assumption (\ref{maintheohypbis}). Consequently $c\geko 0$.

Thanks to Lemma \ref{Gammabasicprop} $(v)$, we get $\lim_{\gamma \to 0} \e \big[  \un_{\{\epp (\ff)=1 \}} 
\int_\bbR (\mathbf G(z)\! -\! \gamma )_+  \mathrm dz\big] \eqo c $,  
by Fubini and monotone convergence. Thus there is $\gamma_0$, which only depends on the law of $\ff$, 
such that $\e \big[  \un_{\{\epp (\ff)=1 \}} 
\int_\bbR (\mathbf G(z)\! -\! \gamma_0 )_+  \mathrm dz\big] \geqo \frac34c$. 
On the other hand  there is $\tau_0$, which only depends on $\eta$ and on the law of $\ff$ such that  
$73 \e \big[ \un_{\{ \epp (\ff)= 1\}}M^{_{(\eta)}}_{\! \ff} \big] \tau^{-\eta} \leqo \frac14 c$ for all $\tau \geko \tau_0$. 
Consequently, we get (\ref{lowerexcess}), which completes the proof of $(iii)$. \cqfd 

\subsection{Estimates on a one parameter family of quadratic Beta density functions}
\label{estiBetasec}
Let $\tau, \sigma, a\ino (0, \infty)$ be such that $\tau \geko \max (1, \sigma) $. 
We define $f$ and $\phi$ as in (\ref{quadBetadef}) and we assume that $\phi \ino \mathscr D$, 
i.e.~$2a(\sigma \tau^2 \! -\! \frac13 \sigma^3)\eqo 1$. In the proof of Theorem \ref{t:main} we use probability densities $\phi$ that correspond to several distinct sets of parameters $(a, \sigma, \tau)$. However, it turns out that the mathematically 
relevant parameters are not $(a, \sigma, \tau)$ but rather $\theta \ino (0, 1)$ and $t\ino (0, \infty)$ given by 
\begin{equation}
\label{relevantparam}
\theta \! :=\!  -\frac{\log \big( 1\! -\! \frac{\sigma}{\tau}\big)}{\log \tau} \, , \, \textrm{i.e.,}\quad \sigma\eqo 
\tau (1\! -\! \tau^{-\theta}) \qquad \textrm{and} \qquad t \! :=\! \frac{3}{4a\tau^2}\eqo  
\frac{3}{4\phi(0)} \, . 
\end{equation}
We want to express $\tau$, and thus $\sigma$, $a$ and $\phi$, in terms of $\theta$ and $t$. To this end we first observe that $1\eqo 2a\tau^3 \Big( \frac{\sigma}{\tau}\! -\! 
\frac13 \big(\frac{\sigma}{\tau} \big)^3 \Big)\eqo 2a\tau^3 ( 1\! -\! \tau^{-\theta} \! -\! \frac{1}{3} (1\! -\! \tau^{-\theta})^3)$, i.e.,
\begin{equation}
\label{tentautheta}
t\eqo \tau \big( 1 \! -\!  \tfrac{3}{2} \tau^{-2\theta} \! + \tfrac{1}{2} \tau^{-3\theta} \big).
\end{equation}
We use the following lemma. 

\begin{lemma}
\label{parameter1} 
Let us fix $\theta\ino (0, 1)$. For all $\tau\ino [1, \infty)$, we set $h_\theta (\tau)\eqo 
\tau \! -\! \frac32\tau^{1-2\theta} \! + \frac12 \tau^{1-3\theta}$. Then the following holds true. 
\begin{compactenum}

\smallskip

\item[$(i)$] $h_\theta$ is an increasing bijection from $[1, \infty)$ onto $\bbR_+$. We denote by 
$\tau_\theta \! : \! \bbR_+ \! \to \! [1, \infty)$ its inverse, which is $C^1$ on $\r^*_+$. 

\smallskip

\item[$(ii)$] There is $t_\theta\ino (0, \infty)$, which only depends on $\theta$, such that for all 
$t\ino [t_\theta \, \infty) $ and all $s\ino \r^*_+$, 
$$\big| \tau_\theta (t+s)\! -\! \tau_\theta (t) \! -\! s \big|\leq \frac{2s}{t^{2\theta}} \quad \textrm{and} \quad 0\leqo \frac{\tau_\theta (t)}{t} \! -\! \frac{\tau_\theta (t+s)}{t+s} \leq \frac{4s}{t^{1+ 2\theta}} \; .$$

\smallskip

\item[$(iii)$] For all $t\ino \r^*_+$, we set $\sigma_\theta (t)\eqo \tau_\theta (t)\! -\! \tau_\theta (t)^{1-\theta}$, $a_\theta (t)\eqo \frac{3}{4 t\tau_\theta (t)^2}$ and for all $v\ino \bbR$, $\phi_{\theta,t} (v)\eqo a_\theta (t) \big( \tau_\theta (t)^2 \! -\! v^2 \big)\un_{\{ |v| <\sigma_\theta (t)\}}$. Then $\phi_{\theta, t}\ino \mathscr D$ and the following holds true as $t \to \! \infty$. 
\begin{compactenum}

\smallskip

\item[$\quad (iii1)$] $\tau_\theta (t)\eqo t \big( 1+ \frac32 t^{-2\theta} \! -\! \frac12 t^{-3\theta} + \mathcal O_\theta ( t^{-4\theta}) \big)$.

\smallskip

\item[$\quad (iii2)$] $\sigma_\theta (t)\eqo t \big( 1\! -\! t^{-\theta} + \frac32 t^{-2\theta} \! -\! \frac12 (4\! -\! 3\theta)t^{-3\theta} + \mathcal O_\theta ( t^{-4\theta}) \big)$. 

\smallskip

\item[$\quad (iii3)$] $a_\theta (t)\eqo \frac34 t^{-3} \big( 1\! -\!  3 t^{-2\theta}\! + t^{-3\theta} + \mathcal O_\theta ( t^{-4\theta}) \big)$. 

\end{compactenum}
Here $\mathcal O_\theta (f(t))$ means that there are $t'$ and $g$, which only depend on $\theta$, such that $ \mathcal O_\theta (f(t))\eqo g(t)f(t)$ for all $t\ino (t', \infty)$ and 
$\sup_{t\in (t', \infty)} |g(t)| \leko \infty$. 

\smallskip

\item[$(iv)$] Let $\theta'\ino (\theta, 1)$. Then there is $t_{\theta, \theta'}$, which only depends on $\theta$ and $\theta'$, such that for all $t\ino [t_{\theta, \theta'}\! ,\,  \infty)$, 
$$ \sigma_\theta (t)\leqo \sigma_{\theta'} (t) < \tau_{\theta'} (t) \leko \tau_\theta (t) \leqo \tau_{\theta'} (t)+ \tfrac{3}{2} t^{1-2\theta} \quad \textrm{and} \quad 0\leqo a_{\theta'} (t)\! -\! a_{\theta} (t) \leqo \tfrac{9}{4} t^{-3-2\theta} .$$

\smallskip

\item[$(v)$] Let $\theta'\ino (\theta, 1)$. Then there is $t_{\theta, \theta'}$, which only depends on $\theta$ and $\theta'$, such that for all $t\ino [t_{\theta, \theta'}, \infty)$ the following holds true.
\begin{compactenum}

\smallskip

\item[$\quad (v1)$] For all $v\ino \bbR$, $\; \phi_{\theta, t} (v) \! -\! \phi_{\theta'\! , t} (v) \leqo (a_{\theta'} (t) \! -\! a_\theta (t)) v^2 \un_{\{|v| <\sigma_\theta (t) \}}$. 
 
\smallskip

\item[$\quad (v2)$] For all $v\ino [-\sigma_\theta (t), \sigma_\theta (t)]$,  $\; \phi_{\theta'\! , t} (v)\leqo  \phi_{\theta, t} (v)$. 
\end{compactenum}
\end{compactenum}
\end{lemma}
\noi
\textbf{Proof.} We first prove $(i)$. We note that $h_\theta$ is clearly $C^2$ on $[1, \infty)$ and  for all $\tau\ino [1, \infty)$,
$$ h'_\theta (\tau)\eqo 1\! -\! \tfrac{3}{2} (1\! -\! 2\theta) \tau^{-2\theta} + \tfrac{1}{2} (1\! -\! 3\theta) \tau^{-3\theta} \quad \textrm{and} \quad h_\theta''(\tau)\eqo  \tfrac{3}{2}\theta  \tau^{-1-2\theta}
\big( 2 (1\! -\! 2\theta)\! -\! (1\! -\! 3\theta)\tau^{-\theta} \big). $$
We set $m\eqo \inf_{\tau \in [1, \infty)} h'_\theta(\tau)$ and we want to prove that $m\geko 0$. If $h'_\theta$ does not reach its minimum value on $[1, \infty)$ then $m\eqo \lim_{\tau \to \infty} h'_\theta (\tau)\eqo 1$, which is positive. If $m\eqo h'_\theta (1)$ then $m\eqo \frac32\theta \geko 0$. 
Otherwise there is $\tau_*\ino (1, \infty)$ such that $h'_\theta (\tau_*)\eqo m$ and 
$h_\theta''(\tau_*) \eqo 0$, namely $(1\! -\! 3\theta)\tau_*^{-\theta}\eqo 2 (1\! -\! 2\theta)$. Therefore 
$m\eqo h'_\theta (\tau_*) \eqo 1 \! -\! \frac12 (1\! -\! 2\theta)\tau_*^{-2\theta} \geko \frac12$. This proves that $h_\theta$ increases on $[1, \infty)$, which entails $(i)$ because $\lim_{\tau \to \infty} h_\theta (\tau)\eqo \infty$ and $h_\theta (1)\eqo 0$. 

To prove $(ii)$, we first observe that (\ref{tentautheta}) implies $\tau_\theta (t) \! \sim \! t$ 
as $t\! \to \! \infty$. Then we also derive from (\ref{tentautheta}) that the derivative in $t$ of $\tau_\theta $ 
is $\tau'_\theta (t)\eqo \big( 1\! -\! \frac32 (1\! -\! 2\theta)  \tau_{\theta} (t)^{-2\theta} \! 
+\tfrac{1}{2} (1\! -\! 3\theta) \tau_{\theta} (t)^{-3\theta} \big)^{-1}$. It implies $\lim_{t\to \infty} t^{2\theta} 
(\tau_\theta'(t) \! -\! 1)\eqo \frac32 (1\! -\! 2\theta)\ino [-\frac32, \frac32]$. Thus there is $t_\theta \ino \r_+^*$, 
which only depends on $\theta$ such that for all $t\ino [t_\theta, \infty)$ and all $s\ino \bbR^*_+$, 
$$ \big| \tau_\theta'(t) \! -\! 1 \big| \leq \frac{2}{t^{2\theta}} \quad \textrm{and} \quad 
\big| \tau_\theta (t+s)\! -\! \tau_\theta (t) \! -\! s \big|\leq 2 \! \int_t^{t+s}\!\!  \frac{\mathrm d r}{r^{2\theta}}
 \leq \frac{2s}{t^{2\theta}} .$$
We next use $\lim_{t\to \infty} \tau'_\theta (t)\eqo 1$ and (\ref{tentautheta}) to get 
$\lim_{t\to \infty} t^{1+2\theta} \frac{\mathrm d}{\mathrm dt} \frac{\tau_\theta (t)}{t}\eqo -3\theta$. 
Without loss of generality, we can assume that $t_\theta$ is sufficiently large to get for all $t\ino [t_\theta, \infty)$,  
$$ -\frac{4\theta}{t^{1+2\theta}} \leq \frac{\mathrm d}{\mathrm dt} \frac{\tau_\theta (t)}{t}\leq - \frac{2\theta}{t^{1+2\theta}} \quad \textrm{and thus} \quad 0\leq \frac{\tau_\theta (t)}{t} \! -\! \frac{\tau_\theta (t+s)}{t+s} \leq 4\theta \! \int_t^{t+s} \!\! \frac{\mathrm dr}{r^{1+2\theta}} \leq  \frac{4s}{t^{1+ 2\theta}} , $$
which completes the proof of $(ii)$.

The fact that $\phi_{\theta, t} \ino \mathscr D$ follows immediately from the definition. The asymptotic expansions in $(iii)$ are obtained by elementary computations: while $(iv)$ is an easy consequence of $(iii)$, we feel free to omit the details. 

It remains to prove $(v)$. By $(iv)$, $\sigma_{\theta} (t)\leqo   \sigma_{\theta'} (t) $ for large $t$, and since 
$a_{\theta} (t)\tau_{\theta} (t)^2 \eqo a_{\theta'} (t)\tau_{\theta'} (t)^2\eqo  \frac{3}{4t}$, we get for all $v\ino \bbR$,  $ \phi_{\theta, t} (v) \! -\! \phi_{\theta'\! , t} (v) \eqo (a_{\theta'} (t) \! -\! a_\theta (t)) v^2 \un_{\{|v| <\sigma_\theta (t) \}} - \phi_{\theta'\! , t} (v) \un_{ \{ \sigma_{\theta} (t) \leq |v| < \sigma_{\theta'} (t) \}}$, which entails $(v1)$ because $\phi_{\theta' , t} (v)$ is nonnegative. 

Let us now suppose that $v\ino [-\sigma_\theta (t), \sigma_\theta (t)]$. Then 
$\phi_{\theta'\!  , t} (v) \eqo \frac{3}{4t}  - a_{\theta'} (t) v^2\leqo  \frac{3}{4t} - a_{\theta} (t) v^2\eqo \phi_{\theta, t} (v) $, by $(iv)$. The proof of $(v)$ is completed. \cqfd

\begin{lemma}
\label{parameter2} 
Let $\theta, \theta'\ino (0, 1)$ be such that $\theta'\geko \theta$. Let $t\ino \r_+^*$. 
For all $v\ino \bbR$ we set 
\begin{equation}
\label{psitdef}
\psi_{\theta, \theta'\! ,\,  t} (v)\eqo \phi_{\theta, t} (v)\un_{\{v\geq 0 \}}+  \phi_{\theta'\! , t} (v)\un_{\{v< 0 \}}. 
\end{equation}
Then, the following holds true.
\begin{compactenum}

\smallskip

\item[$(i)$] $\psi_{\theta, \theta'\! , t}\ino \mathscr D$ and for all $v\ino \bbR$, $\; \psib_{\theta, \theta'\! , t} (v)\eqo \phib_{\theta, t} (v) \un_{\{v\geq 0 \}}+  \phib_{\theta'\! , t} (v)\un_{\{v< 0 \}}$. 

\smallskip

\item[$(ii)$] There is $t_{\theta, \theta'}\ino \bbR^*_+$, which only depends on $\theta$ and $\theta'$, such that for all $t\ino [t_{\theta, \theta'}, \infty)$ and all $v\ino (-\infty, \sigma_\theta (t)]$, $\; \psib_{\theta, \theta'\! , t} (v) \geqo \phib_{\theta'\! ,t}(v)$.

\smallskip

\item[$(iii)$] Let us recall $\mathtt b$ from (\ref{quadBetadef}). Let $t\ino \r^*_+ \! \mapsto \! u_t \ino \r^*_+$ be such that $\lim_{t\to \infty} u_t/t\eqo 1$. Then 
\end{compactenum}
\begin{equation}
\label{psitlim}
\forall x\ino \bbR , \qquad \lim_{t\to \infty} \psib_{\theta, \theta'\! , t} (u_t x)= \overline{\mathtt b} (x).
\end{equation} 

\end{lemma}
\noi
\textbf{Proof.} Since $\phi_{\theta, t}$ and $\phi_{\theta'\! , t}$ belong to $\mathscr D$ by Lemma \ref{parameter1} $(iii)$, Lemma \ref{propLambdaeven} $(i)$ immediately entails $(i)$. 
By $(i)$ we only need to prove $(ii)$ for $v\ino [0, \sigma_\theta (t)]$: by Lemma \ref{parameter1} $(v2)$, there is $t_{\theta, \theta'}\ino \bbR^*_+$,which only depends on $\theta$ and $\theta'$, such that for all $t\ino [t_{\theta, \theta'}, \infty)$, we have 
$\psib_{\theta, \theta'\! , t} (v)\eqo \frac12 + \int_0^v \! \phi_{\theta, t} (w) \, \mathrm d w\geqo  \frac12 + \int_0^v \! \phi_{\theta'\! , t} (w) \, \mathrm d w\eqo \phib_{\theta'\! , t} (v)$, which completes the proof of $(ii)$. 

  To prove $(iii)$, we first apply Lemma \ref{quadBeta1} $(iv)$ to $\phi_{\theta , t}$ 
and $\phi_{\theta'\! , t}$ and get 
$$ \sup_{y\in \bbR} \big|  \phib_{\theta, t} (y) \! -\! \overline{\mathtt b} 
\big( \tfrac{y}{\tau_{\theta}(t) }\big)\big| \leqo 2\Big( 1\! -\! 
\tfrac{\sigma_{\theta} (t)}{\tau_{\theta} (t) }\Big)^2\!\! \eqo 
2\tau_{\theta}(t)^{-2\theta} \; \,  \textrm{and} \; \,   \sup_{y\in \bbR} \big|  \phib_{\theta'\! , t} (y) 
\! -\! \overline{\mathtt b} \big( \tfrac{y}{\tau_{\theta'}(t) }\big)\big| \leqo 2\tau_{\theta'}(t)^{-2\theta'} . $$
Thus by $(i)$, $\sup_{x\in \bbR} \big|  \psib_{\theta,\theta'\! , t} (u_tx) \! -\! \overline{\mathtt b} 
\big( \tfrac{xu_t}{\tau_{\theta}(t) }\big)\un_{\{ x\geq 0\}} \! -\! \overline{\mathtt b}  \big( \tfrac{xu_t}{\tau_{\theta'}(t) }\big)
\un_{\{ x< 0\}}  \big| \leqo 2\tau_{\theta}(t)^{-2\theta}+ 2\tau_{\theta'}(t)^{-2\theta'}$, which implies 
(\ref{psitlim}) because $\tau_\theta (t)\sim t \sim \tau_{\theta'} (t) \sim u_t$ and since $\overline{\mathtt b}$ is continuous. 
 \cqfd 

\medskip

The following lemma provides estimates on $\phib_{\theta, t} (v)$ as $t$ \emph{and} $v$ vary, which is a key estimate in the proof of Theorem \ref{t:main}. 
\begin{lemma}
\label{distribcovar} Let $\theta, \delta \ino (0, 1)$ and let $c\ino \r^*_+$. For all $t\ino \r^*_+$, we recall from Lemma \ref{parameter1} $(i)$ and $(iii)$ the definitions of resp.~$\tau_\theta(t)$ and $\phi_{\theta, t} \ino \mathscr D$. Let $t\ino \r_+^* \! \mapsto \delta_t \ino (0, 1)$ and $t\ino \r_+^* \! \mapsto \Delta_t \ino (0, 1)$ be such that 
\begin{equation}
\label{hypodeltas}
\lim_{t\to \infty} \delta_t = \delta \quad \textrm{and} \quad c'\! := \sup_{t\in [1, \infty)}\!\!\!  t^3\big| t^2\Delta_t \! -\! \tfrac{1}{3} c\big| < \infty .
\end{equation}
Then there exists $t_0\ino \r_+^*$, which only depends on $\theta, \delta, c$ and $c'$, such that for all $t\ino [t_0, \infty)$, 
\begin{equation}
\label{coaccr}
\forall v\ino [-\sigma_\theta (t), \sigma_\theta (t)], \quad \phib_{\theta, t+ \Delta_t} 
(v+ \delta_t \Delta_t) -\phib_{\theta, t} (v) \geq  \frac{c}{4t^2 \tau_\theta (t)^4} \big(\tfrac12 \delta\tau_\theta (t)  \! -\! v \big)  \big( \tau_\theta (t)^2 \!\! -\! v^2\big). 
\end{equation}
\end{lemma}
\noi
\textbf{Proof.} To simplify notation we set $\; \Delta\! :=\! \Delta_t$, $\;  h \! :=\!  \delta_t \Delta_t$,  
$$\sigma\!: =\! \sigma_\theta (t), \quad   \tau \! :=\!  \tau_\theta (t), \quad x \! :=\!  \frac{v}{\tau},   \quad  t_*\! :=\!  t+ \Delta, \quad  \sigma_*\!: =\! 
 \sigma_{\theta} (t_*) , \quad  \tau_*\! :=\! \tau_\theta (t_*), \quad x_*\! :=\!  \frac{\tau x+h}{\tau_*}  . $$
We first derive from Lemma \ref{parameter1} $(ii)$ and $(iii)$  that 
$\tau_*\! -\! \tau \sim \sigma_* \! -\! \sigma \sim \Delta$ as $t\! \to \! \infty$. 
This implies that there is $t_{1} $, which only depends on $\theta, \delta , c$ and $c'$, such that for all 
$t\ino [t_{1} \, ,  \infty) $  and  all $v\ino [-\sigma, \sigma]$, 
\begin{equation}
\label{interlace}
-\sigma_*\leqo v+ h\leqo \sigma+ \sigma_* \! -\! \sigma \eqo \sigma_* \; .
\end{equation}
We also deduce from  Lemma \ref{parameter1} $(ii)$ that 
there is $t_{2}\ino (t_1, \infty)$, which only depends on $\theta$, such that for all $t\ino [t_2, \infty)$; 
\begin{equation}
\label{controltau0}
|\tau_* \! -\! \tau \!-\! \Delta| \leqo \frac{2\Delta}{t^{2\theta}} \qquad \textrm{and} 
\qquad 0\leq \frac{\tau}{t}  - \frac{\tau_*}{t_*} \leq \frac{4\Delta}{t^{1+ 2\theta}} .
\end{equation}
We want to bound $\tfrac{1}{t_*} (\tau_* \! -\! \tau) - \tfrac{c}{3t^2\tau}$. To this end, 
we note that (\ref{hypodeltas}) implies $\Delta \eqo \frac{c}{3t^2} (1+ \mathcal O_{c,c'} (t^{-3}))$ and then 
$t_* \eqo t + \Delta\eqo  t \big( 1+  \frac{c}{3t^3} + \mathcal O_{c, c'} (t^{-6}))$.
By (\ref{controltau0}), we get  
$\tau_* \! -\! \tau=  \frac{c}{3t^2} (1+ \mathcal O_{c, c'\! , \, \theta} (t^{-2\theta}))$. Therefore 
$\tfrac{1}{t_*} (\tau_* \! -\! \tau)\eqo  \frac{c}{3t^3} (1+ \mathcal O_{c, c'\! , \, \theta} (t^{-2\theta}))$. 
By Lemma \ref{parameter1} $(iii)$, we get $ \tfrac{c}{3t^2\tau}\eqo 
\frac{c}{3t^3} (1+ \mathcal O_{c, c'\! , \, \theta} (t^{-2\theta}))$ and thus finally 
$\tfrac{1}{t_*} (\tau_* \! -\! \tau) - \tfrac{c}{3t^2\tau} \eqo  \mathcal O_{c, c'\! , \, \theta} (t^{-3-2\theta})$. 
Consequently, by (\ref{controltau0}) and the previous arguments there are 
$ t_3\ino (t_2, \infty)$ and $c''\ino \r_+^*$, which only depend on $\theta$, $c$ and $c'$, 
such that for all $t\ino [t_3, \infty)$, 
\begin{equation}
\label{controltau}
\Delta +|\tau_* \! -\! \tau| \leq \frac{c''}{t^{2}}, \quad  \Big| \frac{\tau_* \! -\! \tau}{t_*} \!  - \! \frac{c}{3t^2\tau} \Big|
 \leq \frac{c''}{t^{3+2\theta}} 
\quad \textrm{and} \quad 0\leq \frac{\tau}{t}  \! -\!  \frac{\tau_*}{t_*} \leq\frac{c''}{t^{3+2\theta}} .
\end{equation}

Next, for all $y\ino \bbR$, we recall that $P(y)\eqo \frac34 (y\! -\! \frac13 y^3)$. For all $s\ino \r_+^*$ by  Lemma \ref{quadBeta1} $(ii)$ and since $\frac43 a_\theta (s) \tau_{\theta}(s)^3\eqo \tau_\theta (s)/s$, we get 
\begin{equation}
\label{rewritephib}
\forall w\ino [-\sigma_\theta (s), \sigma_\theta (s)], \quad \phib_{\theta, s} (w) = \tfrac12+ \tfrac{\tau_\theta (s)}{s}P \big( \tfrac{w}{\tau_\theta (s)}\big). 
\end{equation} 
Let us fix $v\ino [-\sigma, \sigma]$ and set 
$\Delta\phib\eqo \phib_{\theta, t + \Delta_t} 
(v+ \delta_t \Delta_t) \! -\! \phib_{\theta, t} (v) \eqo \phib_{\theta, t_*} (v+h) \! -\!  \phib_{\theta, t} (v) $.
Then by (\ref{interlace}), (\ref{rewritephib}) applies to $s\eqo t_*$ and $w\eqo v+h$ and to $s\eqo t$ and $w\eqo v$ and we get 
\begin{eqnarray}
\Delta \phib &=& \tfrac{\tau_*}{t_*} P\big( \tfrac{v+h}{\tau_*}\big) -\tfrac{\tau}{t} P\big( \tfrac{v}{\tau}\big)= \tfrac{\tau_*}{t_*} P\big( x_*\big) -\tfrac{\tau}{t} P\big( x\big) \nonumber \\
& =& 
\underbrace{\big(\tfrac{\tau_*}{t_*}  \! -\!  \tfrac{\tau}{t}  \big) P(x)}_{=: \, S_{1,t} (v)} \, + \, 
\underbrace{\, \tfrac{\tau_*}{t_*} \int_0^{x_*-x} \!\!\!\!\!\!\!\!\!\!\!\! \big( P'(y+x) \! -\! P'(x) \big) \mathrm dy\, }_{=: \, S_{2,t} (v)} \, + \, \underbrace{\, \tfrac{1}{t_*} \tau_*(x_* \! -\! x) P'(x)\, }_{=: \, S_{3,t} (v)} . \label{doubledeltas}
\end{eqnarray}

  We intend to prove that $S_{3,t} (v)$ is the dominant term. To this end we first observe that $P$ increases on $[-1, 1]$ from $-\frac12$ to $\frac12$. Thus  (\ref{controltau}) implies for all $t\ino [t_{3}, \infty)$ that 
\begin{equation}
\label{controlS1}
\sup_{v\in [-\sigma, \sigma]} |S_{1,t} (v)| \leq \frac{c''}{2t^{3+2\theta}} \; .
\end{equation}

 We next get an upper bound for  $|S_{2,t} (v)|$. We first check that $P'(y+x) \! -\! P'(x)\eqo -\frac34y(2x+y) $. Then we note that (\ref{interlace}) implies that $x_*\ino [-1, 1]$. 
 By assumption, $x\ino [-1, 1]$ too, and since $y$ in the integral of $S_{2, t} (v)$ is between $0$ and $x_*\! -\! x$, we get $|2x+y|\leqo 4$ and thus $|P'(y+x) \! -\! P'(x)|\leqo 3 |y|$. Therefore, for all $t\ino [t_3, \infty)$ and all $v\ino [-\sigma, \sigma]$, we get 
$$ \big| S_{2,t} (v)\big| \leq \tfrac{\tau_*}{t_*} \int_0^{|x_*-x|} \!\!\!\!\!\!\!\!\!\!\!\! 3 |y| \, \mathrm dy = \tfrac{3}{2t_*\tau_*} \big( \tau_*(x_* \! -\! x) \big)^2 \leq \tfrac{3}{2t^2} \big( \tau_*(x_* \! -\! x) \big)^2  .
$$

We observe that $\tau_* (x_*\! -\! x) \eqo \delta_t \Delta \! -\! (\tau_* \! -\! \tau) x$. Therefore by (\ref{controltau}), for all $t\ino [t_3, \infty)$, we get $| \tau_* (x_*\! -\! x)| \leqo c''t^{-2}$ and thus 
\begin{equation}
\label{controlS2}
\sup_{v\in [-\sigma, \sigma]} |S_{2,t} (v)| \leq \frac{3(c'')^2}{2t^6} \; .
\end{equation}

Finally we 
get a lower bound for $S_{3, t} (v)$ which 
holds true for all sufficiently large $t$ and for all $v\ino [-\sigma, \sigma]$. 
To this end we first observe that 
$$ \frac{\tau_*}{t_*} (x_* \! -\! x)= \frac{h}{t_*}  - \frac{\tau_* \! -\! \tau}{t_*} x = \; 
\underbrace{\Big( \delta_t \frac{\Delta}{t_*} \! -\! \frac{c\delta}{6t^2\tau} \Big)}_{=:V_{1,t}} \, +
 \frac{c}{3t^2\tau} \big(\tfrac{1}{2} \delta \! -\! x \big)  \, -\, \underbrace{\Big( \frac{\tau_* \! -\! \tau}{t_*} \!  
 - \! \frac{c}{3t^2\tau}\Big) x}_{=:V_{2,t} (v)} . $$
We note that $V_{1,t}$ does not depend on $v$ and that $V_{1,t}\sim \frac{c\delta}{6t^3}$ by (\ref{hypodeltas}), whereas by (\ref{controltau}) and since $|x| \leqo 1$, we get $\sup_{v\in [-\sigma, \sigma]} |V_{2,t} (v)| \leqo c''t^{-3-2\theta}$.  
Since $-\frac{\sigma}{\tau} \leqo x\eqo \frac{v}{\tau}  \leqo \frac{\sigma}{\tau} $, we get 
$\frac34 \geqo P'(x) \geqo P'(\frac{\sigma}{\tau}) \eqo \frac34(1+ \frac{\sigma}{\tau}) (1\! -\! \frac{\sigma}{\tau}) \geqo \frac34 \tau^{-\theta}$. Thus there is $t_4\ino (t_3, \infty)$, which only depends on $\delta, \theta$, $c$ and $c'$, such that for all 
$t\ino [t_4, \infty)$ and all $v\ino [-\sigma, \sigma]$, 
 \begin{eqnarray*}
S_{3,t} (v) \!\!\! &= &\!\!\!  V_{1, t}  P'(x)-    V_{2, t}  (v) P'(x)+ \frac{c}{3t^2\tau} \big(\tfrac12 \delta \! -\! x \big) P'(x) \\
\!\!\! &\geq &\!\!\! \frac{c\delta}{16t^{3+ \theta}} - \frac{3c''}{4t^{3+2\theta}}  +  \frac{c}{3t^2\tau} \big(\tfrac12 \delta \! -\! x \big) P'(x) 
\end{eqnarray*}
and by (\ref{doubledeltas}), (\ref{controlS1}) and (\ref{controlS2}), it implies that for all $v\ino [-\sigma, \sigma]$, 
\begin{equation}
\label{lowerbounddeltaphib}
\Delta \phib - \frac{c}{3t^2\tau} \big(\tfrac12 \delta \! -\! x \big) P'(x) \geq \frac{c\delta}{16t^{3+ \theta}} - \frac{3c''}{2t^{3+2\theta}}  -\frac{3(c'')^2}{2t^6} \, .
\end{equation}
Since the right-hand side of (\ref{lowerbounddeltaphib}), which does not depends on $v$, tends to $\infty$ as $t\! \to \! \infty$, there is $t_0 \geqo t_4$, which only depends on $\theta, \delta, c$ and $c'$, such that for all $t\ino [t_0, \infty)$ and all $v\ino [-\sigma, \sigma]$, $\Delta \phib \! -\!   \frac{c}{3t^2\tau} \big(\frac12 \delta \! -\! x \big) P'(x) \geqo 0$, which is exactly (\ref{coaccr}). This completes the proof of the lemma. \cqfd

\section{Proof of the main result}
\label{proofmainsec}

\subsection{A general recipe for lower bounds}
\label{s:proof1}

This section gives a general lower bound for 1-homogeneous systems upon a technical condition which involves the Lambda function $\Lambda_{\psi, \ff}$ as in Definition \ref{Lambdadef}.

Before stating and proving this result we recall basic properties of stochastic order for $[-\infty, \infty]$-valued r.v.s. More precisely let $X, Y$ be two  $[-\infty, \infty]$-valued r.v.s. We say that $X$ is stochastically smaller 
than $Y$ if for all $z\ino \bbR$, $\p (X \geko  z ) \leqo  \p (Y \geko z) $, which is denoted by $X\leqo_{\mathtt{st}} Y$. Let us suppose that $(\Omega, \mathscr F,  \p)$ is sufficiently rich to carry a r.v.~which is uniformly distributed on $[0, 1]$ and which is independent from $(X,Y)$. Then there are two r.v.s $X'$ and $Y'$ such that $X$ and $X'$ have the same law, $Y$ and $Y'$ have the same law and a.s.~$X'\leqo Y'$. The converse holds true trivially. 
Therefore if $f\! : \! [-\infty, \infty] \! \to \! [-\infty, \infty]$ is nondecreasing, then $f(X)\leq_{\mathtt{st}} f(Y)$. 

We prove the following theorem 
which is a general recipe to obtain a lower bound for critical $1$-homogeneous systems.   
\begin{theorem}
\label{p:Lambda_cond=>lb} Let $\ff$ be a random $\mathscr H$-valued function. Let $X_0$ be an $\r_+^*$-valued r.v.~and let $(X_n)_{n\in \bbN}$ be defined as in (\ref{syshom}). We suppose that 
\begin{equation}
\label{recipehyp1}
\p \big(\mathbf F \! \notin \! \{ \varphi_+, \varphi_-\} \big)\geko 0  \quad  \textrm{and} \quad \e[\sgn (\ff) ] \eqo  0 .
\end{equation}
We fix $\delta \ino (0, \infty)$, $\nu\ino \bbN$, $\frak{t}\ino (0, \infty)$ and for all $n\ino \bbN$, 
$$ \sigma_n \ino \bbR, \quad  \Delta_n \ino \r_+^*, \quad q_n \ino (0, 1) \quad \textrm{and} \quad \psi_n \ino \mathscr D \; \textrm{such that} \; \psib_n (\sigma_n)\eqo 1 . $$
We define $t_n \! :=\!  \frak t$ if $n\leqo \nu$ and 
$t_n \! :=\! \frak{t} + \sum_{\nu \leq  k<n } \Delta_k$ if $n\geko \nu$. For all $n\ino \bbN$ and for all $v\ino \bbR$ we also set 
$(\Delta\psib_n )(v) \eqo \psib_{n+1} (v+ \delta\Delta_n) \! -\! \psib_n(v) $. 
We assume that for all $n\geq \nu$,  
\begin{equation}
\label{Lambdahyp}
 \forall v\ino (-\infty ,  \sigma_{n+1} \! -\! \delta \Delta_n ), \; \,  \e \big[ \Lambda_{\psi_n, \ff} (v) \big] \! + \frac{1\! -\! q_{n+1}}{(1\! -\! q_n)^2} (\Delta \psib_n) (v) + \frac{q_{n+1} \! -\! q_n}{(1\! -\! q_n)^2} \big( 1\! -\! \psib_n(v)\big) \geq  0.
\end{equation}
Then there is a constant $c\ino \r_+^*$ which only depends on the law of $X_0$ and on 
$\sigma_{\nu}$ and $q_{\nu}$, such that 
\begin{equation}
\label{recipe}
\forall n\geqo \nu, \; \forall x\ino \bbR, \quad \p \big( \log X_{n-\nu} \geko xt_n\big) \geq 
( 1\! -\! q_n)\big( 1\! -\! \psib_n \big( (x+ \delta)t_n + c \big)\big).
\end{equation}
\end{theorem}
\noi
\textbf{Proof.} Let $Z_n, Z'_n$ be two $[-\infty, \infty)$-valued r.v.s.~which have the same law given by 
$$ \forall z \ino \bbR, \quad \p (Z_n \leqo z)\eqo \p(Z'_n \leqo z)\eqo q_n + (1\!-\! q_n)\psib_n(z). $$ 
In particular, $\p (Z_n \eqo -\infty)\eqo   \p (Z'_n \eqo -\infty)\eqo q_n$. We furthermore assume that $(Z_n)_{n\in \bbN}$, $(Z'_n)_{n\in \bbN}$ and $\ff$ are independent. We set for all $n\ino \bbN$, 
$$ V_n \eqo \log \ff \big( \ee^{Z_n}\! ,  \ee^{Z'_n}  \big) \; $$
with the conventions that $\ee^{-\infty}\eqo 0$ and $\log 0\eqo -\infty$, and where $\ff$ is extended continuously on $\bbR^2_+$ as explained in (\ref{contiexthomog}). By Assumption (\ref{recipehyp1}), we can apply (\ref{centredLambda}) in Lemma \ref{elempropLambda} $(iii)$: for all $v\ino \bbR$ we thus get 
\begin{eqnarray}
\label{leftmember}
\p (Z_{n+1} \leqo v+ \delta \Delta_n)  \!\!\!\!  &- & \!\!\!\!   \p (V_n \leqo v) \\
\eqo  (1\! -\! q_n)^2   \Big(  \e \big[  \!\!\!\!  \!\!\!\!   & &  \!\!\!\!  \!\!\!\!\  \Lambda_{\psi_n, \ff} (v) \big]+\frac{1\! -\! q_{n+1}}{(1\! -\! q_n)^2} (\Delta \psib_n) (v) + \frac{q_{n+1} \! -\! q_n}{(1\! -\! q_n)^2} \big( 1\! -\! \psib_n(v)\big) \Big)  . \nonumber
\end{eqnarray}
By Assumption (\ref{Lambdahyp}), this quantity is nonnegative if 
$v\leko \sigma_{n+1} \! -\! \delta \Delta_n$. If $v\geqo \sigma_{n+1} \! -\! \delta \Delta_n$, the 
left-hand side of (\ref{leftmember}) is clearly equal to $1\! -\! \p (V_n \leqo v) \geqo 0$ since 
$\psib_{n+1} (\sigma_{n+1}) \eqo 1$. Thus it implies that $\p (V_n \leqo v)  \leqo  \p (Z_{n+1}\! -\!  \delta \Delta_n \leqo v)$ 
for all $v\ino \bbR$, i.e., that for all $n\geqo \nu$, 
\begin{equation}
\label{interpretationLambdahyp}
Z_{n+1}\! -\!  \delta \Delta_n \leq_{\mathtt{st}} V_n\; .
\end{equation}

Since $\p (X_0 \geko 0)\eqo 1$, we see that 
$\lim_{\lambda \to 0^+} \p (X_0 \leko \lambda) \eqo 0 \leko q_{\nu}$ and there is $\lambda_*\ino (0, 1)$ such that $\p (X_0 \leko \lambda_*)  \leko q_{\nu}$. We set $c\! :=\! \sigma_{\nu} \! -\! \log \lambda_*$, which only depends on the law of $X_0$ and on 
$(\sigma_{\nu}, q_{\nu})$. 

We want to prove recursively for all $n\ino \bbN$ the property $\mathtt{Rec}_n$ given by 
\begin{equation}
\label{Recn}
\mathtt{Rec}_n: \qquad Z_{n+\nu}\! -\!  \delta (t_{n+\nu} \! -\! t_{\nu} ) \! -\! c \leq_{\mathtt{st}} \log X_n \; .
\end{equation}
\noi
\emph{Proof of $\; \mathtt{Rec}_0$.}  If $v\ino \r_+$, then $\p (\log X_0 \! -\! \log \lambda_* \leqo v) \leqo 1 \eqo \p (Z_{\nu} \! -\! \sigma_{\nu} \leqo v)$. If $v\ino (-\infty, 0)$, the very definition of $\lambda_*$ implies that $\p (\log X_0 \! -\! \log \lambda_* \leqo v) \leqo \p (\log X_0 \! -\! \log \lambda_* \leko 0) \leko q_{\nu} \leqo \p (Z_{\nu} \! -\! \sigma_{\nu} \leqo v)$. This completes the proof of $\mathtt{Rec}_0$. \cq 

\smallskip

\noi
\emph{Proof of $\; \mathtt{Rec}_n \Rightarrow \mathtt{Rec}_{n+1}$.} Let us assume $\mathtt{Rec}_n$. Since $Z_{n+\nu}$, $Z'_{n+\nu}$ have the same law, and since $Z_{n+\nu}, Z'_{n+\nu}$ and $\ff$ are independent, a standard argument implies that we can find independent r.v.s 
$X'_n$ and $X''_n$ which have the same law as $X_n$, which are independent of $\ff$ and such that 
a.s.~$\log X_n'\geqo Z_{n+\nu}\! -\!  \delta (t_{n+\nu} \! -\! t_{\nu} ) \! -\! c$ and  
$\log X''_n\geqo Z'_{n+\nu}\! -\!  \delta (t_{n+\nu} \! -\! t_{\nu} ) \! -\! c$. Since $\ff$ is nondecreasing in each coordinate and is $1$-homogenenous, it implies a.s., that 
\begin{eqnarray*}
\log  \ff (X'_n , X''_n) \!\!\! &\geq & \!\!\! \log \ff \big( \ee^{Z_{n+\nu}} e^{-\delta (t_{n+\nu}  - t_{\nu} )  - c }, \ee^{Z'_{n+\nu}} e^{-\delta (t_{n+\nu} - t_{\nu} )  - c }\big) \\
\!\!\! &= & \!\!\!  V_{n+ \nu} \! -\! \delta (t_{n+\nu} \! -\!  t_{\nu} ) \! -\! c .
\end{eqnarray*}
It entails $\mathtt{Rec}_{n+1}$ by (\ref{interpretationLambdahyp}) since $\log X_{n+1} \overset{\textrm{law}}{=} \log  \ff (X'_n , X''_n)$. \cq

\smallskip

It remains to prove that $\mathtt{Rec}_n$ implies (\ref{recipe}). Indeed, we see for all $n\ino \bbN$ and all $x\ino \bbR$ that  
\begin{eqnarray*}
\p \big( \log X_n \geko xt_{n+ \nu}\big)  \!\!\! &\geq & \!\!\! \p \big( Z_{n+\nu}\! -\!  \delta (t_{n+\nu} \! -\! t_{\nu} ) \! -\! c \geko x t_{n+ \nu}) \eqo \p \big( Z_{n+\nu} \geko c + (x+ \delta) t_{n+ \nu} \! -\! \delta t_{\nu} )\\
 \!\!\! &\geq & \!\!\! \p \big( Z_{n+\nu} \geko c + (x+ \delta) t_{n+ \nu}  )   \quad \textrm{(because $ \delta t_{\nu} \geko 0$)}\\
 \!\!\! &\eqo & \!\!\!  ( 1\! -\! q_{n+\nu})\big( 1\! -\! \psib_{n+\nu} \big( (x+ \delta)t_{n+\nu} + c \big)\big), 
\end{eqnarray*}
which is the desired result. \cqfd

\subsection{Proof of Theorem \ref{t:main}}
\label{proofmainthsec}
Let $\ff$ be an $\mathscr H$-valued random function and let $\eta\ino (0, 1)$. We recall from (\ref{Gamma}) that $M^{_{(\eta)}}_{\! \ff}\eqo 
\Gamma_{\! \ff}^{1+\eta,1}\! + \Gamma_{\! \ff}^{0,2+\eta} $ and 
we recall from (\ref{momentmono}) that 
for all $(a,b)\ino \r_+ \! \times \! [1, \infty)$ such that 
$a+b\leko 2+ \eta$, there is a constant $c_{a,b,\eta} \ino (0, \infty)$, which only depends on $a,b$ and $\eta$, such that $\Gamma_{{\! \ff}}^{{a,b}} \leq c_{a,b,\eta} \max \big(1, M^{(\eta)}_{\! \ff} \big)$. 
We assume that (\ref{maintheohyp}) holds true, namely, 
\begin{equation}
\label{maintheohypter}
\p \big(\mathbf F \! \notin \! \{ \varphi_+, \varphi_-\} \big)\geko 0, \quad \e \big[ M^{_{(\eta)}}_{\! \mathbf F}\big]\leko \infty \quad  \textrm{and} \quad \e[\sgn (\ff) ] \eqo \e \big[ \Gamma_{\! \ff}^{0,1} \sgn (\ff) \big]\eqo 0 .
\end{equation}
We also recall from (\ref{Gamma}) and  (\ref{cste}) that 
\begin{equation}
\label{cstarrecall}
\mathtt{c} (\ff)     \! :=\! 
\Gamma_{\! \ff}^{1,1}\! + \tfrac{1}{2}\Gamma_{\! \ff}^{0,2} \quad \textrm{and} \quad c_* \eqo 
\tfrac92 \e [\mathtt c (\ff) ]\, .
\end{equation}

\vspace{-4mm}

\begin{remark}
\label{hypinvF} 
 {\rm 
Since $\Gamma^{a,b}_{\! \ff}\eqo \Gamma^{a,b}_{\! \ff^*}\eqo \Gamma^{a,b}_{\! \ff^{\inv}}$
and $\epp (\ff^\inv)\eqo -\epp (\ff)$, we see that $\ff$ satisfies (\ref{maintheohypter}) if and only if $\ff^\inv$ also satisfies (\ref{maintheohypter}), and in this case $\mathtt c (\ff)\eqo 
\mathtt c (\ff^\inv)$. \cq 
 } 
\end{remark}
To prove Theorem \ref{t:main}, we apply Theorem \ref{p:Lambda_cond=>lb} to specific $\sigma_n, \Delta_n, q_n$ and $\psi_n$ that are defined as follows. We fix $\delta \ino (0, 1)$, which can be arbitrary small and we also fix the following parameters  $\theta, \theta'$ and $\kappa$ such that 
\begin{equation}
\label{paramrange}
0\leko \theta \leko \tfrac{\eta}{1+\eta} , \quad 0\leko 3\kappa \leko \eta
\quad \textrm{and} \quad \tfrac12 (1+ 3\kappa) \leko   \theta' \leko 1. 
\end{equation}
Note that $\theta'\geko \frac12\geko \theta$. We recall from Lemma \ref{parameter1} the notations  
$\tau_\theta (\cdot)$, $\tau_{\theta'} (\cdot)$, $\sigma_\theta (\cdot)$, $\sigma_{\theta'} (\cdot)$, etc. 
For all $n\ino \bbN^*$ we set 
\begin{eqnarray}
\label{setup}
t_n  \!\!\!& \! :=\!  &\!\!\! (c_*n)^{\frac{1}{3}}, \quad \Delta_n \! :=\! t_{n+1} \! -\! t_n  , \quad q_n \! :=\! \tfrac{1}{10}\delta \big( 1\! -\! n^{-\kappa} \big) \, ; \\
 \tau_n\!\!\!& \! :=\!  &\!\!\! \tau_{\theta} (t_n), \quad \sigma_n \! :=\!  \sigma_\theta (t_n),  \quad a_n \! :=\! a_\theta (t_n)= \frac{3}{4t_n\tau_n^2},  \quad\phi_n (\cdot) \! :=\! \phi_{\theta, t_n}  (\cdot)\, ;  \label{positivesetup} \\
\tau^\bullet_n  \!\!\!& \! :=\!  &\!\!\! \tau_{\theta'} (t_n), \quad \sigma^\bullet_n  \! :=\!  \sigma_{\theta'} (t_n) ,  \quad a^\bullet_n \! :=\! 
a_{\theta'} (t_n)\eqo \frac{3}{4t_n (\tau^\bullet_n)^2},   \quad \textrm{and}  \quad \phi^\bullet_n  (\cdot)  \! :=\!  \phi_{\theta'\! , \, t_n}  (\cdot)\, .
\label{negativesetup}
\end{eqnarray}
We define  $\psi_n$ by 
\begin{equation}
\label{psindef}
\forall v\ino \bbR, \quad \psi_n (v)= \phi_n (v)\un_{\{ v\geq 0 \}}+   \phi^\bullet_n (v)\un_{\{ v< 0 \}}. 
\end{equation}
Since $\phi_n$ and $\phi^\bullet_n$ belong to $\mathscr D$, so does $\psi_n$ by Lemma \ref{parameter2} $(i)$ because $\psi_n \eqo \psi_{\theta, \theta'\! , \, t_n}$. We first prove the following. 
\begin{lemma}
\label{param3}. We keep the previous notation. Then the following holds true. 
\begin{compactenum}

\smallskip

\item[$(i)$] $t_n \sim \tau_n \sim \tau^\bullet_n$. 

\smallskip

\item[$(ii)$] $0 \leqo t_n^3 \big( \tfrac13 c_* \! -\! t_n^2 \Delta_n\big) \leqo \tfrac19 c_*^2$.

\smallskip

\item[$(iii)$] $\Delta_n \sim \tau_{n+1} \! -\!  \tau_{n} \sim  \tau^\bullet_{n+1} \! -\! \tau^\bullet_n \sim \frac{1}{3} c_* t_n^{-2} \sim 
 \frac{1}{3} c_*\tau_n^{-2}  \sim \frac{1}{3} c_* (\tau^\bullet_n)^{-2} $. 
\end{compactenum}
\end{lemma}
\noi
\textbf{Proof.} $(i)$ is a direct consequence of Lemma \ref{parameter1} $(iii)$. To prove $(ii)$ we observe that 
$$ \Delta_n = (c_*n)^{\frac13} \Big(  (1+ \tfrac1n)^{\frac13} \! -\! 1\Big) = \frac{c_*}{3(c_*n)^{2/3}} - (c_*n)^{\frac13} \Big( 1+ \tfrac{1}{3n} \! -\!   (1+ \tfrac1n)^{\frac13} \Big). $$
To get $(ii)$, we use the following inequality $0\leqo  1\! +\!  \tfrac13 x \! -\! (1\! +\! x)^{1/3}
 \leqo \tfrac19 x^2$, $x\ino [0, 1]$.

Let us prove $(iii)$. We recall that $h_\theta$ and $h_{\theta'}$ are the inverse of resp.~$\tau_\theta$ and $\tau_{\theta'}$. Thus 
$h_\theta (\tau_n)\eqo h_{\theta'} (\tau^\bullet_n)\eqo t_n$ and 
$\Delta_n \eqo \int_{\tau_n}^{\tau_{n+1}} h'_\theta (\tau) \, \mathrm d\tau$. Since $\lim_{\tau \to \infty}h'_\theta (\tau) \eqo 1$, we get $\Delta_n \eqo t_{n+1} \! -\! t_n \sim\tau_{n+1} \! -\! \tau_n $. Similarly, we get $\Delta_n\sim  \tau^\bullet_{n+1} \! -\! \tau^\bullet_n$. On the other hand, we derive from $(ii)$ that 
$\Delta_n \sim \tfrac13 c_*t_n^{-2}$, which completes the proof of $(iii)$ thanks to $(i)$.  \cqfd 

\begin{proposition}
\label{p:main} 
Let $\eta, \delta \ino (0, 1)$. Let $\ff$ be an $\mathscr H$-valued random function. Let $X_0$ be an $\r_+^*$-valued r.v.~and let $(X_n)_{n\in \bbN}$ be defined as in (\ref{syshom}). We assume that $\ff$ satisfies (\ref{maintheohyp}) (which is recalled in (\ref{maintheohypter})). 
Let $\theta, \theta', \kappa $ be as in (\ref{paramrange}). Let $\sigma_n, \Delta_n, \psi_n, q_n$, 
$n\ino \bbN^*$, be as in (\ref{setup}), (\ref{positivesetup}), (\ref{negativesetup}) and (\ref{psindef}). Then we can find $\nu\ino \bbN^*$, which only depends on the parameters $\eta, \delta, \theta, \theta'\! $, $\kappa$ and on the law of $\ff$, such that Assumption (\ref{Lambdahyp}) in Theorem \ref{p:Lambda_cond=>lb} holds true with $\frak t\eqo t_{\nu}$. 
\end{proposition}
\noi
\textbf{Proof.} Proposition \ref{p:main} is the main technical part of the proof of Theorem \ref{t:main}. Its proof is postponed to Section \ref{Lambdachecksec}. \cqfd 

\medskip

\noi
\textbf{Proof of Theorem \ref{t:main}}. We admit Proposition \ref{p:main} and use it to complete the proof of  Theorem \ref{t:main}. We fix $\delta\ino (0, 1)$ that can be arbitrarily close to $0$. Since $\ff$ satisfies (\ref{maintheohyp}) and Proposition \ref{p:main}, we can apply Theorem \ref{p:Lambda_cond=>lb} with $\frak t \eqo t_{\nu}$ and we get for all $n\geqo \nu$ and all $x\ino \bbR$ 
$$ \p \big( \log X_{n-\nu} \geko xt_{n-\nu}\big) \geq \p \big( \log X_{n-\nu} \geko xt_n\big) \geq 
( 1\! -\! q_n)\big( 1\! -\! \psib_n \big( (x+ \delta)t_n + c \big)\big) \, .$$
We note that $\lim_{n\to \infty} (1\! -\! q_n) \eqo 1\! -\! \frac{1}{10}\delta$. We apply (\ref{psitlim}) in Lemma \ref{parameter2} $(iii)$ to $\psi_n \eqo \psi_{\theta, \theta'\! , \, t_n}$ and $u_{t_n}\eqo t_n + \frac{c}{x+ \delta} \sim t_n$. This yields that 
if $\ff$ satisfies (\ref{maintheohyp}), then for all $\delta\ino (0, 1)$ and all $x\ino \bbR$, 
$$ \liminf_{n\to \infty} \p \big(  \log X_{n} \geko xt_{n}\big) \geq \big(1\! -\! \tfrac{1}{10} \delta \big) \big( 1\! -\! \overline{\mathtt b} (x+\delta) \big)$$
where $\mathtt b \ino \mathscr D$ is as in (\ref{quadBetadef}) in Lemma \ref{quadBeta1}. Since $\delta$ can be aribtrarily close to $0$ and since $\overline{\mathtt b}$ is continuous, we have proved that for all $\mathscr H$-valued random functions $\ff$ satisfying (\ref{maintheohyp}), 
\begin{equation}
\label{onesidelimit}
\forall x\ino \bbR, \quad   \liminf_{n\to \infty} \p \big(  \log X_{n} \geko xt_{n}\big) \geq 1\! -\! \overline{\mathtt b}(x) .
\end{equation}
We observe that (\ref{syshom}) implies $\frac{1}{X_{n+1}} \overset{\textrm{law}}{=} \ff^\inv (\frac{1}{X_n}, \frac{1}{X'_n})$.  By Remark \ref{hypinvF}, $\ff^\inv$ also satisfies (\ref{maintheohyp}). Therefore we get 
$$ \forall x\ino \bbR, \quad   \liminf_{n\to \infty} \p \big(  -\log X_{n} \geko -xt_{n}\big) \geq 1\! -\! \overline{\mathtt b}(-x) . $$
Since $\liminf_{n\to \infty} \p \big(  -\log X_{n} \geko -xt_{n}\big)\eqo 1\! -\! \limsup_{n\to \infty} \p \big(  \log X_{n} \geqo xt_{n}\big)$ and $1\! -\! \overline{\mathtt b} (-x)\eqo \overline{\mathtt b} (x)$, we get $\limsup_{n\to \infty} \p \big(  \log X_{n} \geqo xt_{n}\big)\leqo 
1\! -\! \overline{\mathtt b} (x)$ for all $x\ino \bbR$. Combined with (\ref{onesidelimit}) it implies (\ref{wcv}) and the proof of Theorem \ref{t:main} is completed. \cqfd 
\subsection{Proof of Proposition \ref{p:main}}
\label{Lambdachecksec} 
The proof of Proposition \ref{p:main} is the main technical part of the proof of Theorem 
\ref{t:main}. It heavily relies on the estimates proved in Section \ref{prelimsec}. Here we fix 
$\eta, \delta \ino (0, 1)$. Let $\ff$ be an $\mathscr H$-valued random function. Let $X_0$ be an $\r_+^*$-valued r.v.~and 
let $(X_n)_{n\in \bbN}$ be defined as in (\ref{syshom}). We assume that $\ff$ satisfies (\ref{maintheohyp}) 
(which is recalled in (\ref{maintheohypter})). 
Let $\theta, \theta', \kappa $ be as in (\ref{paramrange}). 
Let $\sigma_n, \Delta_n, \psi_n, q_n$, 
$n\ino \bbN^*$, be as in (\ref{setup}), (\ref{positivesetup}), (\ref{negativesetup}) and (\ref{psindef}). 
In what follows we say that \emph{an integer only depends on our set of parameters} 
to mean that it only depends on $\eta, \delta, \theta, \theta'\! , \kappa$ and on the law of $\ff$. 
To simplify notation we set $K_\eta \eqo \e \big[ M^{_{(\eta)}}_{\! \ff }\big]$. 

We first prove that there is $n_1\ino \bbN^*$, which only depends on our set of parameters, 
such that for all $n\geqo n_1$, 
\begin{equation}
\label{symmetryfication}
\forall v\ino \bbR, \quad \e \big[ \Lambda_{\psi_n, \ff} (v) \big]\, \geq \, 
\e \big[ \Lambda_{\phi_n, \ff} (v) \big]\un_{\{ v\geq 0 \}} +\e \big[ \Lambda_{\phi^\bullet_n, \ff} (v) \big] \un_{\{v<0 \}} - 2K_\eta t_n^{-3-2\theta-\eta}.
\end{equation}
\noi
\textbf{Proof of (\ref{symmetryfication})}. 
We first apply Lemma \ref{parameter1} $(v1)$ to $ \phi_{\theta, t_n}\eqo \phi_n$ and $ \phi_{\theta'\! , \, t_n}\eqo \phi^\bullet_n$: namely, there is $n_0\ino \bbN$, which only depends on our set of parameters, such that for all $n\geq n_0$ and all $v\ino \bbR$, we get $\phi_n (v)\! -\! \phi^\bullet_n (v) \leqo (a^\bullet_n \! -\! a_n) v^2 \un_{\{ |v| < \sigma_n\}}$. We note that $\psi_n, \phi_n$ and 
$\phi^\bullet_n$ are maximal at $0$, so Lemma \ref{propLambdaeven} $(iv)$ applies to 
$F\eqo \ff$, $\psi\eqo \psi_n$, 
$\phi_+\eqo \phi_n$, $\phi_-\eqo \phi_n^\bullet$, $\sigma\eqo \sigma_n$ and $K\eqo  a^\bullet_n \! -\! a_n$ and (\ref{symmetrization}) asserts that for all $n\geqo n_0$ and all $v\ino \bbR$, 
 \begin{equation}
\label{symmetryfication1}
\e \big[ \Lambda_{\psi_n, \ff} (v) \big]\, \geq \, 
\e \big[ \Lambda_{\phi_n, \ff} (v) \big]\un_{\{ v\geq 0 \}} +\e \big[ \Lambda_{\phi^\bullet_n, \ff} (v) \big] \un_{\{v<0 \}} - (a^\bullet_n\! -\! a_n) \psi_n (0) \sigma_n^{1-\eta} K_\eta  \, . 
\end{equation}
By Lemma \ref{parameter1} $(iv)$, for all sufficiently large $n$ we get $0\leqo 
(a^\bullet_n\! -\! a_n) \psi_n (0) \sigma_n^{{1-\eta}}
\leqo \frac94 t_n^{-3-2\theta}\psi_n (0) \sigma_n^{1-\eta}$. 
Since $\psi_n (0)\eqo \phi_n(0)\eqo 3/(4t_n)$, and $\tau_n \sim t_n\sim \sigma_n$, we get 
$  \frac94  t_n^{-3-2\theta}\psi_n (0) \sigma_n^{1-\eta} \sim \frac{27}{16} t_n^{-3-2\theta-\eta} $. So there is an integer $n_1\geqo n_0$, which only depends on our parameters, such that  $(a^\bullet_n\! -\! a_n) \psi_n (0) \sigma_n^{1-\eta}\leqo 2 t_n^{-3-2\theta-\eta}$ for all $n\geqo n_1$ and it implies (\ref{symmetryfication}) thanks to (\ref{symmetryfication1}). \cq  

\medskip

For all $n\ino \bbN$ and all $v\ino \bbR$, we recall from Theorem \ref{p:Lambda_cond=>lb} the notation $ (\Delta\psib_n )(v) \eqo \psib_{n+1} (v+ \delta \Delta_n) \! -\! \psib_n(v)$ and set, moreover,  
$$(\Delta\phib_n) (v) \! :=\!  \phib_{n+1} (v+ \delta\Delta_n) \! -\! \phib_n(v) \quad \textrm{and} \quad (\Delta\phib^{_\bullet}_n) (v) \! :=\!  \phib^{_\bullet}_{n+1} (v+ \delta\Delta_n) \! -\! \phib^{_\bullet}_n(v) .$$
We prove that there is an integer $n_2\geqo n_1$, which only depends on our parameters, such that for all $n\geqo n_2$, 
\begin{equation}
\label{reductio1}
\delta \Delta_n \leko \sigma_{n+1} \; \quad \textrm{and} \quad   \forall v\ino \bbR, \quad 
(\Delta\psib_n) (v) \geq (\Delta\phib_n) (v) \un_{\{ v\geq 0\}} + (\Delta\phib^{_\bullet}_n) (v)\un_{\{ v< 0\}}. 
\end{equation}
\noi
\textbf{Proof of (\ref{reductio1})}. Since $\lim_{n\to \infty}\delta\Delta_n  \eqo 0$ and $\lim_{n\to \infty} \sigma_{n+1}\eqo \infty$, there is $n_2\geqo n_1$, which only depends on our parameters, such that $\delta \Delta_n \leko \sigma_{n+1}$ for all $n\geqo n_2$. 

We first assume that $v\ino \bbR_+$ and observe that $v+ \delta \Delta_n \geqo 0$. We apply Lemma \ref{parameter2} $(i)$ 
to get  $\psib_{n+1} (v+ \delta\Delta_n) \eqo \phib_{n+1} (v+ \delta\Delta_n)$ and $\psib_n(v)\eqo \phib_n(v)$, which 
implies that $(\Delta \psib_n)(v)\eqo (\Delta \phib_n)(v)$. 

Let us suppose next that $v+ \delta  \Delta_n \leko 0$. Then $v\leko 0$, and 
Lemma \ref{parameter2} $(i)$ implies that $\psib_{n+1} (v+ \delta\Delta_n) \eqo \phib^{_\bullet}_{n+1} (v+ \delta\Delta_n)$ and $\psib_n(v)\eqo \phib^{_\bullet}_n(v)$, which entails $(\Delta \psib_n)(v)\eqo (\Delta \phib^{_\bullet}_n)(v)$. 

It remains to consider the case where $v+ \delta\Delta_n \geqo 0 \geko v$. 
Since $\delta\Delta_n\leko \sigma_{n+1}$, Lemma \ref{parameter2} $(ii)$ applies and we get 
$\psib_{n+1} (v+ \delta \Delta_n) \geqo \phib^{_\bullet}_{n+1} (v+ \delta\Delta_n)$.  Since $v\leko 0$, 
Lemma \ref{parameter2} $(i)$ entails $\psib_n(v)\eqo \phib^{_\bullet}_n(v)$. Thus $(\Delta \psib_n)(v)\geqo (\Delta \phib^{_\bullet}_n)(v)$, which completes the proof of (\ref{reductio1}). \cq

\medskip

As soon as $\tau_n\geqo 1$, we have 
$$\sigma_{n+1} -\delta\Delta_n \eqo \sigma_n + \tau_{n+1} \! -\! \tau_n + \tau_{n}^{1-\theta} \! -\! \tau_{n+1}^{1-\theta} -\delta\Delta_n \leko \sigma_n + \tau_{n+1} \! -\! \tau_n \! -\! \delta \Delta_n . $$
By Lemma \ref{param3} $(iii)$, there is $n_3\geqo n_2$, which only depends on our parameters, such that 
\begin{equation}
\label{debord}
\forall n\geqo n_3, \quad 
\sigma_{n+1} -\delta\Delta_n \leko \sigma_n +  2(1\! -\! \delta) \Delta_n . 
\end{equation}

Therefore by (\ref{symmetryfication}), (\ref{reductio1}) and (\ref{debord}), and since $1 - \psib_n (v) \geqo \tfrac12$ if 
$v\leqo 0$, Proposition \ref{p:main}, i.e.~(\ref{Lambdahyp}), will be a straightforward consequence of the following limits: 
\begin{eqnarray} 
\label{positivenewlim}
 \lim_{n\to \infty} \; t_n^{3+2\theta+\eta} \!\!\!\!\!\!\!\!\! \!\! & &\!\!\!\!\!\!\!\! \inf_{v\in [0, \sigma_n+ 2(1-\delta)\Delta_n]}
  \Big( \e \big[\Lambda_{\phi_n, \ff} (v) \big]  + \frac{1\! -\! q_{n+1}}{(1\!-\! q_n)^2} (\Delta \phib_n)(v)\Big) = \infty
   \quad \textrm{and} \\
\label{negativenewlim}
\lim_{n\to \infty} \; t_n^{3+2\theta+\eta} \!\!\!\!\!\!\!\!  & &\!\!\!\!\!\!  \inf_{v\in (-\infty, 0] } \Big( \e \big[\Lambda_{\phi^\bullet_n, \ff} (v) 
\big]  + \frac{1\! -\! q_{n+1}}{(1\!-\! q_n)^2} (\Delta \phib^{_\bullet}_n) (v) + \frac{q_{n+1}\! -\! q_n }{2(1\!-\! q_n)^2}\Big) = \infty .
\end{eqnarray}
\noi
\textbf{Proof of (\ref{positivenewlim})}. We consider two cases: $v\ino [\sigma_n, \sigma_n + 2(1\! -\! \delta) \Delta_n ]$ 
and $v\ino [0, \sigma_n)$. 

\medskip

\noi
\textbf{Case 1:} we suppose that $v\ino [\sigma_n, \sigma_n + 2(1\! -\! \delta) \Delta_n ]$. To simplify notation we set $D_n \! :=\! \tau_n \! -\! \sigma_n \eqo \tau_n^{1-\theta}$. We want to apply Lemma \ref{mainLambdaesti} $(iii)$ to $\tau\eqo \tau_n$, $\sigma\eqo \sigma_n$, $a\eqo a_n$ (and $D\eqo D_n$ and $\phi\eqo \phi_n$, necessarily). 
To this end, 
we first observe that the first inequality in (\ref{paramrange}) implies that 
$\frac12 \leko \frac{1}{1+ \eta} \leko 1\! -\! \theta \leko 1$. Therefore there is $n_4 \geqo n_3$, which only depends on our parameters, such that for all $n\geqo n_4 $
\begin{equation}
\label{largeD}
\tau_n^{2(1-\theta) }\eqo D_n^2 \geko \tau_n \geko \sigma_n \geko \tfrac12 \tau_n \geko 2 \quad \textrm{and} \quad 2(1\! -\! \delta) \Delta_n \leko \gamma_0, 
\end{equation}
where $\gamma_0$ is as in Lemma \ref{mainLambdaesti} $(iii)$, which applies and asserts the following: 
$c_1\! :=\! \e \big[ \Gamma^{0, 1}_{\! \ff} \un_{\{ \epp (\ff)= 1\}} \big]$ is a positive quantity and for all 
$n\geqo n_4$, (\ref{lowerexcess}) holds true. Thus, 
\begin{equation}
\label{lowerexcessbis} \inf_{v\in [\sigma_n, \sigma_n + 2(1- \delta) \Delta_n ]}  \!\!\!\!  \e \big[ \Lambda_{\phi_n, \ff}  (v)\big]  \geq \tfrac{1}{8} c_1 a_n^2 \tau_n^2 D_n^2 = \tfrac{1}{8} c_1 a_n^2\tau_n^{2+ 2(1-\theta)} .
\end{equation}
  We note that for all $v\ino [\sigma_n, \sigma_n + 2(1\! -\! \delta) \Delta_n ]$, $\phib_n (v)\eqo \phib_n (\sigma_n)\eqo 1$. Therefore 
 \begin{equation}
 \label{stepounet}
 v\ino [\sigma_n, \sigma_n + 2(1\! -\! \delta) \Delta_n ], \quad 
( \Delta\phib_n) (v)\geqo \phib_{n+1} (\sigma_n+ \delta \Delta_n) \! -\! \phib_n (\sigma_n)\eqo (\Delta \phib_n) (\sigma_n).
 \end{equation}

  To get a lower bound for $(\Delta \phib_n) (\sigma_n)$, we want to apply Lemma \ref{distribcovar} to $t\eqo t_n$, $\tau_{\theta} (t_n)\eqo \tau_n$, $\phi_{\theta, t_n}\eqo \phi_n$, $\Delta_{t_n}\eqo t_{n+1} \! -\! t_n\eqo \Delta_n$, $\delta_{t_n} \eqo \delta $ and $c\eqo c_*$. 
To this end, we use Lemma \ref{param3} to see that (\ref{hypodeltas}) is satisfied: 
\begin{equation}
\label{hypodeltasbis}
\lim_{n\to \infty} \delta_{t_n} = \delta \quad \textrm{and} \quad  \sup_{n\geq 1}\;  t_n^{3} \big| t_n^2\Delta_{n}\!  \! -\! \tfrac13 c_* \big| \leqo \tfrac19 c_*^2 .
\end{equation}
Since $\phi_{\theta, t_n + \Delta_{n}}\eqo \phi_{\theta, t_{n+1}} \eqo \phi_{n+1}$, we have 
$(\Delta \phib_n) (v)\eqo \phib_{\theta, t_n + \Delta_{n}} (v+ \delta_{t_n} \Delta_{n} )   \! -\! \phib_{\theta, t_n} (v)$. Lemma  \ref{distribcovar} applies, namely, there is $n_5 \geqo n_4$, which only depends on our parameters, such that (\ref{coaccr}) holds true, i.e., for $n\geqo n_5$, 
\begin{equation}
\label{coaccrbis}
\forall v\ino [-\sigma_n, \sigma_n], \quad (\Delta \phib_n) (v)\geq  \frac{c_*}{4t_n^2 \tau_n^4} \big(\tfrac12 \delta\tau_n  \! -\! v \big)  \big( \tau_n^2 \!\! -\! v^2\big)= \tfrac49 c_* a_n^2 \big(\tfrac12 \delta\tau_n  \! -\! v \big)  \big( \tau_n^2 \!\! -\! v^2\big), 
\end{equation}
since $1/(t_n\tau_n^2)\eqo \tfrac43a_n$ by (\ref{positivesetup}). 

We now apply (\ref{coaccrbis}) to 
$v\eqo \sigma_n$. First note that $\tfrac12 \delta\tau_n  \! -\! \sigma_n\geqo -(1\! -\! \frac12 \delta) \tau_n$ and that $\tau_n^2 \! -\! \sigma_n^2 \eqo D_n (\tau_n + \sigma_n)\leqo 2D_n\tau_n\eqo 2\tau_{n}^{2-\theta}$. By (\ref{coaccrbis}), we get 
\begin{equation}
\label{coaccrsigman}
\forall n\geqo n_5 , \quad  (\Delta \phib_n) (\sigma_n)\geq  -\tfrac89c_* a_n^2\big( 1 \! -\! \tfrac12 \delta \big) \tau_n^{3-\theta} .
\end{equation}

On the other hand, we see that $\lim_{n\to \infty} (1\!-\! q_{n+1})/(1\! -\! q_n)^2\eqo 1/(1\! -\! \frac{1}{10}\delta) \leko 2$. Thus by (\ref{lowerexcessbis}),  (\ref{stepounet})  and  (\ref{coaccrsigman}), there is $n_6\geqo n_5$, which only depends on our parameters, such that for all $n\geqo n_6$, 
$$ \inf_{v\in [\sigma_n, \sigma_n + 2(1- \delta) \Delta_n ]}  \!\!  \Big(  \e \big[ \Lambda_{\phi_n, \ff}  (v)\big]  + \frac{1\! -\! q_{n+1}}{(1\! -\! q_n)^2} (\Delta \phib_n) (v)  \Big) \geq  \tfrac{1}{8} c_1 a_n^2\tau_n^{2+ 2(1-\theta)} -\tfrac{16}{9}c_* a_n^2\big( 1 \! -\! \tfrac12 \delta \big) \tau_n^{3-\theta} .$$
Since $a_n \eqo \frac{3}{4 t_n \tau^2_n}$ and since $t_n \! \sim \! \tau_n$, we get 
$a_n^2\tau_n^{2+ 2(1-\theta)} \! \sim \! \frac{9}{16} t_n^{-2-2\theta} $ and  
$a_n^2 \tau_n^{3-\theta} \! \sim \frac{9}{16} t_n^{-3-\theta}\eqo  o \big(t_n^{-2-2\theta} \big)$ since $\theta \leko 1$. Since $t_n^{3+2\theta+\eta}  t_n^{-2-2\theta}\eqo t_n^{1+\eta} \! \to \! \infty$, we get the desired estimate:
\begin{equation}
\label{case1}
\lim_{n\to \infty} \; t_n^{3+2\theta+\eta} \!\!\!\!\!  \inf_{v\in [\sigma_n, \sigma_n+ 2(1-\delta)\Delta_n]} \Big( \e \big[\Lambda_{\phi_n, \ff} (v) \big]  + \frac{1\! -\! q_{n+1}}{(1\!-\! q_n)^2} (\Delta \phib_n)(v)\Big) = \infty.
\end{equation}

\noi
\textbf{Case 2:} we suppose that $v\ino [0, \sigma_n )$. We recall from (\ref{cstarrecall}) that $\frac29 c_* \eqo  \e [\mathtt c (\ff) ]$ and 
we apply Lemma \ref{mainLambdaesti} $(i)$ to $\tau\eqo \tau_n$, $\sigma\eqo \sigma_n$, $a\eqo a_n$ (and $D\eqo D_n$ and $\phi\eqo \phi_n$, necessarily) and get the following: for all 
$n\geqo n_6$, 
$$ \forall v \ino [0, \sigma_n), \quad  \e \big[ \Lambda_{\phi_n, \ff} (v) \big] \geq \tfrac49c_* a_n^2v(\tau_n^2 \! -\! v^2)-198 K_\eta a_n^2 \tau_n^{3-\eta} . $$
There exists $n_7\geqo n_6$, which only depends on our parameters, such that for all $n\geqo n_7$, we have $1\geko \frac{(1-q_{n})^2}{1-q_{n+1}} \geko 1\! -\! \frac19 \delta$. Thus,  
\begin{equation}
\label{posbulk1}
 \forall v \ino [0, \sigma_n), \quad  \frac{(1\! -\! q_{n})^2}{1\! -\! q_{n+1}}\e \big[ \Lambda_{\phi_n, \ff} (v) \big] \geq \tfrac49c_* ( 1\! -\! \tfrac19 \delta) a_n^2v(\tau_n^2 \! -\! v^2)-198 K_\eta a_n^2 \tau_n^{3-\eta} . 
\end{equation}
Since $v\ino [0, \sigma_n)$,  $\tfrac49c_* ( 1\! -\! \tfrac19\delta ) a_n^2v(\tau_n^2 \! -\! v^2) \geqo  \tfrac49c_*  a_n^2v(\tau_n^2 \! -\! v^2) -   \tfrac{4}{81} c_* \delta a_n^2\tau_n(\tau_n^2 \! -\! v^2)$. 
By (\ref{coaccrbis}) and (\ref{posbulk1}), for all $n\geqo n_7$ and all $v\ino [0, \sigma_n)$, we get 
\begin{eqnarray*}
\frac{(1\! -\! q_{n})^2}{1\! -\! q_{n+1}}\e \big[ \Lambda_{\phi_n, \ff} (v) \big] + (\Delta \phib_n)(v) \!\!\!\! &\geq & \!\!\!\! 
\tfrac{4}{9} c_* a_n^2 \tfrac12 \delta\tau_n   \big( \tau_n^2 \!\! -\! v^2\big) - 
  \tfrac{4}{81} c_* \delta a_n^2\tau_n(\tau_n^2 \! -\! v^2) -198K_\eta a_n^2 \tau_n^{3-\eta} \\
\!\!\!\! & \geq & \!\!\!\! \tfrac{14}{81} c_* \delta a_n^2 \tau_n (\tau_n^2 \! -\! v^2) -198 K_\eta a_n^2 \tau_n^{3-\eta}\\
\!\!\!\! & \geq & \!\!\!\! \tfrac{14}{81} c_* \delta a_n^2 \tau_n^{3-\theta} -198K_\eta a_n^2 \tau_n^{3-\eta}
\end{eqnarray*}
since $ \tau_n^2 \! -\! v^2\eqo (\tau_n + v) (\tau_n \! -\! v) \geqo \tau_n (\tau_n\! -\! \sigma_n)\eqo \tau_n^{2-\theta}$. 
We recall from (\ref{paramrange}) that $\theta\leko \eta$. Therefore there is $n_8\geqo n_7$, 
which only depends on our parameters, such that for all $n\geqo n_8$, 
$$ \tfrac{14}{81} c_* \delta a_n^2 \tau_n^{3-\theta} -198 K_\eta a_n^2 \tau_n^{3-\eta} \geko \tfrac{7}{81} c_* \delta \tfrac{9}{16} t_n^{-6} t_n^{3-\theta} \geko \tfrac{1}{32} c_* \delta t_{n}^{-3-\theta} \quad \textrm{and} \quad \frac{1\! -\! q_{n+1}}{(1\! -\! q_n)^2} \geko 1, $$ 
which implies 
\begin{equation}
\label{case2}
t_n^{3+2\theta+\eta} \!\!\!\!  \inf_{v\in [0,\sigma_n)} \Big( \e \big[\Lambda_{\phi_n, \ff} (v) \big]  + \frac{1\! -\! q_{n+1}}{(1\!-\! q_n)^2} (\Delta \phib_n)(v)\Big) \geq \tfrac{1}{32} c_* \delta  t_n^{\theta +\eta} \underset{n\to\infty}{\longrightarrow} \infty.
\end{equation}
This, combined with (\ref{case1}), completes the proof of (\ref{positivenewlim}). \cq

\bigskip

\noi
\textbf{Proof of (\ref{negativenewlim})}. For brevity we set $D^\bullet_n\eqo \tau_n^\bullet \! -\! \sigma^\bullet_n\eqo (\tau_n^\bullet)^{1-\theta'}$ and we consider two cases: $v\ino [-\sigma^\bullet_n + \frac12 D^\bullet_n, 0]$ and $v\in (-\infty, -\sigma^\bullet_n + \frac12 D^\bullet_n]$. 

\medskip

\noi
\textbf{Case 3:} we suppose $v\ino [-\sigma^\bullet_n + \frac12 D^\bullet_n\, , \, 0]$. We recall from (\ref{cstarrecall}) that $\frac29 c_* \eqo  \e [\mathtt c (\ff) ]$. 
There is $n_9\geqo n_8$, which only depends on our parameters, such that for all $n\geqo n_9$, 
$$ \sigma^\bullet_n \geko \tfrac12 \tau^\bullet_n  \geko 2 .$$
In particular, $ \frac12 D^\bullet_n \! -\! \sigma_n^\bullet  \leko 0$. Thus, we can apply Lemma \ref{mainLambdaesti} $(i)$ to $\tau\eqo \tau^\bullet_n$, $\sigma\eqo \sigma^\bullet_n$, $a\eqo a^\bullet_n$ (and $D\eqo D^\bullet_n$ and $\phi\eqo \phi^\bullet_n$, necessarily) to get the following: there exists $n_{10}\geqo n_9$, which only depends on our parameters, such that for all  
$n\geqo n_{10}$ 
$$ \forall v \ino [-\sigma^\bullet_n + \tfrac12 D^\bullet_n \, ,\,  0], \quad  \e \big[ \Lambda_{\phi^\bullet_n, \ff} (v) \big] \geq \tfrac49c_* (a^{\bullet}_n)^2v \big((\tau^\bullet_n)^2 \! -\! v^2 \big)-198 K_\eta (a^\bullet_n)^2 (\tau^\bullet_n)^{3-\eta} . $$
Moreover, since $n_{10}\geqo n_{8}$, 
we get $\frac{(1-q_{n})^2}{1-q_{n+1}} \leko 1$. Thus for all  
$n\geqo n_{10}$ and all $v\ino [ \frac12 D^\bullet_n \! -\! \sigma^\bullet_n \, , 0]$, 
\begin{equation}
 \label{negbulk}
 \frac{(1\! -\! q_{n})^2}{1\! -\! q_{n+1}}\e \big[ \Lambda_{\phi^\bullet_n, \ff} (v) \big] 
 \geq \tfrac49c_*  (a^\bullet_n)^2v\big( (\tau^\bullet_n)^2 \! -\! v^2\big)-198 K_\eta (a^\bullet_n)^2 (\tau^\bullet_n)^{3-\eta}
\end{equation}
because $v\big( (\tau^\bullet_n)^2 \! -\! v^2 \big) \leqo 0$. 

   We want to apply Lemma \ref{distribcovar} to $t\eqo t_n$, $\tau_{\theta'} (t_n)\eqo 
\tau^\bullet_n$, $\phi_{\theta'\! , \, t_n}\eqo \phi_n^\bullet$, $\Delta_{t_n}\eqo t_{n+1} \! -\! t_n\eqo \Delta_n$, 
$\delta_{t_n} \eqo \delta$ and $c\eqo c_*$. We recall from (\ref{hypodeltasbis})
that (\ref{hypodeltas}) is satisfied.  
We also observe that $\phi_{\theta'\! ,\,  t_n + \Delta_{n}}\eqo \phi_{\theta'\! , \, t_{n+1}} \eqo \phi^\bullet_{n+1}$. Therefore $(\Delta \phib^{_\bullet}_n) (v)\eqo \phib_{\theta'\! ,\,  t_n + \Delta_{n}} (v+ \delta_{t_n} \Delta_{n} )   \! -\! \phib_{\theta'\! , \, t_n} (v)$. 
Thus Lemma \ref{distribcovar} applies and asserts that 
there is $n_{11} \geqo n_{10}$, which only depends on our parameters, such that (\ref{coaccr}) holds true for all $n\geqo n_{11}$, namely, for all $v\ino [-\sigma^\bullet_n, \sigma^\bullet_n]$, 
\begin{equation}
\label{coaccrter}
(\Delta \phib^{_\bullet}_n) (v)\geq  \frac{c_*}{4t_n^2 (\tau^\bullet_n)^4} \big(\tfrac12 \delta\tau^\bullet_n  \! -\! v \big)  \big( (\tau^\bullet_n)^2 \!\! -\! v^2\big)= \tfrac49 c_* (a^\bullet_n)^2 \big(\tfrac12 \delta\tau^\bullet_n  \! -\! v \big)  \big( (\tau^\bullet_n)^2 \!\! -\! v^2\big), 
\end{equation}
since $1/(t_n(\tau^\bullet_n)^2)\eqo \tfrac43a^\bullet_n$ by (\ref{negativesetup}). 

   As $n\! \to \! \infty$ we observe  that 
$$ \frac{q_{n+1} \! -\! q_n}{2(1\! -\! q_{n+1})} \sim \frac{\delta\kappa n^{-1-\kappa}}{2(10 \! -\! \delta)}\sim 
 \frac{\delta \kappa (c_*)^{1+\kappa}}{2(10  \! -\!  \delta)} \big( (c_*n)^{\frac{1}{3}} \big)^{-3(1+\kappa)} =2c_\bullet t_n^{-3(1+\kappa)}, $$
where we have set $c_\bullet\eqo \frac{\delta \kappa (c_*)^{1+\kappa}}{4(10  -  \delta)} \ino (0, \infty)$ to simplify the notation. Thus, 
there exists $n_{12}\geqo n_{11}$, which only depends on our parameters, such that for $n\geqo n_{12}$, 
\begin{equation}
\label{thirdmembre}
\forall n\geqo n_{12}, \qquad \frac{q_{n+1} \! -\! q_n}{2(1\! -\! q_{n+1})} \geko c_\bullet t_n^{-3(1+\kappa)} \, .
\end{equation}
By (\ref{negbulk}), (\ref{coaccrter}) and (\ref{thirdmembre}), we get for all $n\geqo n_{12}$ and all $v\ino [\frac12 D^\bullet_n\!  -\! \sigma^\bullet_n \,  , 0]$, 
\begin{eqnarray}
 \frac{(1\! -\! q_n)^2}{1\! -\! q_{n+1}} \e \big[ \Lambda_{\phi^\bullet_n, \ff} (v) \big] \!\! \!\! &+ & \!\! \!\! (\Delta \phib^{_\bullet}_n) (v)+ \frac{q_{n+1} \! -\! q_n}{2(1\! -\! q_{n+1})} \nonumber \\
 \!\! \!\! &\geq &\!\! \!\!   \tfrac{2}{9}  c_* \delta (a^\bullet_n)^2 \tau^\bullet_n  \big( (\tau^\bullet_n)^2 \!\! -\! v^2\big) + c_\bullet t_n^{-3(1+\kappa)} -198K_\eta (a^\bullet_n)^2 (\tau^\bullet_n)^{3-\eta}\nonumber \\
 \!\! \!\! &\geq &\!\! \!\!   c_\bullet t_n^{-3(1+\kappa)} -198 K_\eta (a^\bullet_n)^2 (\tau^\bullet_n)^{3-\eta} \label{lower3}
\end{eqnarray} 
because $c_* \delta (a^\bullet_n)^2 \tau^\bullet_n  \big( (\tau^\bullet_n)^2 \!\! -\! v^2\big) \geko 0$ if 
$v\ino [\frac12 D^\bullet_n\!  -\! \sigma^\bullet_n \,  , 0]$. By Lemma \ref{param3}, 
$a^\bullet_n\sim \frac34t_n^{-3}$, $(a^\bullet_n)^2 (\tau^\bullet_n)^{3-\eta}  \sim \frac{9}{16} t_{n}^{-3-\eta}\eqo o (t_n^{-3(1+\kappa)})$, since $3\kappa \leko \eta$ by (\ref{paramrange}). 
Thus, there is $n_{13}\geqo n_{12}$, which only depends on our parameters, such that for all $n\geqo n_{13}$, 
\begin{equation}  
\label{lowerremainingpart}
c_\bullet t_n^{-3(1+\kappa)} -198 K_\eta (a^\bullet_n)^2 (\tau^\bullet_n)^{3-\eta} \geko \tfrac12 c_\bullet t_n^{-3(1+\kappa)} \quad \textrm{and} \quad \frac{1\! -\! q_{n+1}}{(1\!-\! q_n)^2} \geko 1\; .
\end{equation}
Thus by (\ref{lower3}) and (\ref{lowerremainingpart}), we get 
\begin{equation}\label{case3}
t_n^{3+2\theta+ \eta} \!\!\!\!\!\!\! \!\! \! \!\! \inf_{v\in [\frac12 D^\bullet_n -\sigma_n^\bullet\, ,\,  0]}\! \!  \Big( \e \big[\Lambda_{\phi^\bullet_n, \ff} (v) \big]  \! +\!  \frac{1\! -\! q_{n+1}}{(1\!-\! q_n)^2} (\Delta \phib^{_\bullet}_n)(v)\!  +\!  \frac{q_{n+1}\! -\! q_n }{2(1\!-\! q_n)^2} \Big) \geqo \tfrac12 c_\bullet t_n^{2\theta+ \eta-3\kappa} \!\!\! \underset{n\to \infty}{\longrightarrow}\!  \infty
 \end{equation} 
because $2\theta+ \eta-3\kappa >0$ by (\ref{paramrange}).

\medskip

\noi

\noi
\textbf{Case 4:} we now suppose $v\ino (-\infty,   \frac12 D^\bullet_n\! -\! \sigma^\bullet_n ]$. We first observe that for all $n\geqo n_{13}$, 
\begin{equation}
\label{Deltaposi}
\forall v \ino (-\infty, \tfrac12 D^\bullet_n \! -\! \sigma_n^\bullet ] , \quad (\Delta\phib_n^{_\bullet}) (v) \geq 0 .
\end{equation}
\emph{Indeed,} this is a consequence of (\ref{coaccrter}) if 
$v\ino [-\sigma_n^\bullet,  \tfrac12 D^\bullet_n -\sigma_n^\bullet ]$, whereas if $v\leko -\sigma_n^\bullet$, we observe that $\phib_n^{_\bullet} (v)\eqo 0$ and therefore $(\Delta\phib_n^{_\bullet}) (v)\eqo \phib_{n+1}^{_\bullet} (v+ \delta \Delta_n)\geq 0$, which completes the proof of (\ref{Deltaposi}). 

\smallskip

We observe that 
$$ \frac{q_{n+1} \! -\! q_n}{2(1\! -\! q_n)^2} \sim \frac{\delta \kappa n^{-1-\kappa}}{20 (1\! -\! 
\tfrac{1}{10}\delta)^2} =c' t_n^{-3(1+ \kappa)} \quad \textrm{where} \quad c'\! :=\! 
\frac{\delta \kappa c_*^{1+\kappa}}{20 (1\! -\! \tfrac{1}{10}\delta)^2} \geko 0 \; .$$
Thus there exists $n_{14}\geqo n_{13}$, which only depends on our parameters, such that 
\begin{equation}
\label{maincase4}
\forall n\geqo n_{14}, \qquad  \frac{q_{n+1} \! -\! q_n}{2(1\! -\! q_n)^2} \ > \tfrac12 c' t_n^{-3(1+\kappa)} \; .
\end{equation}

For brevity, we set $c''\! :=\! \e \big[ \Gamma^{0, 1}_{\! \ff} + \mathtt c (\ff) \big]$. 
We apply Lemma \ref{mainLambdaesti} $(ii)$ to $\tau\eqo \tau^\bullet_n$, $\sigma\eqo \sigma^\bullet_n$, $a\eqo a^\bullet_n$ (and $D\eqo D^\bullet_n$ and $\phi\eqo \phi^\bullet_n$, necessarily) to get the following: there exists $n_{15}\geqo n_{14}$, which only depends on our parameters, such that for all  
$n\geqo n_{15}$, 
\begin{equation}
\label{neglow}
\forall v \ino (-\infty, \tfrac{1}{2}D^\bullet_n \! -\! \sigma^\bullet_n], 
\quad  \e \big[ \Lambda_{\phi^\bullet_n, \ff} (v) \big] \geq -198 K_\eta (a^\bullet_n)^2 (\tau^\bullet_n)^{3-\eta} - 15c'' (a^{\bullet}_n)^2 (\tau^\bullet_n)^{4-2\theta'} . 
\end{equation}
because $D_n^\bullet\eqo (\tau_n^\bullet)^{1-\theta'}$. 

By Lemma \ref{param3} and (\ref{negativesetup}), we get $\tau^\bullet_n \sim t_n$ and 
$a^\bullet_n\sim \frac34 t_n^{-3}$. Therefore, 
 \begin{equation}
\label{whoisbigger}  
(a^{\bullet}_n)^2 (\tau^\bullet_n)^{4-2\theta'} \! \sim \! \tfrac{9}{16} t_{n}^{-2-2\theta'}\eqo o \big( t_n^{-3(1+\kappa)}\big) \; \,  \textrm{and} \; \,  (a^\bullet_n)^2 (\tau^\bullet_n)^{3-\eta}\!  \sim \! \tfrac{9}{16} t_{n}^{-3-\eta}\eqo o \big( t_n^{-3(1+\kappa)}\big)
\end{equation}
since $3+ 3\kappa \leko 2+ 2\theta'$ and $3\kappa \leko \eta$ by (\ref{paramrange}). By (\ref{Deltaposi}), (\ref{maincase4}), (\ref{neglow}) and (\ref{whoisbigger}), there is $n_{16}\geqo n_{15}$, which only depends on our parameters, such that for all $n\geqo n_{16}$, 
\begin{equation}
\label{case4}
t_n^{3+2\theta+ \eta} \!\!\!\!\!\!\! \! \inf_{v \leq \frac12 D^\bullet_n -\sigma_n^\bullet}\!  \Big( \e \big[\Lambda_{\phi^\bullet_n, \ff} (v) \big]  + \frac{1\! -\! q_{n+1}}{(1\!-\! q_n)^2} \Delta \phib^{_\bullet}_n(v) + \frac{q_{n+1}\! -\! q_n }{2(1\!-\! q_n)^2} \Big) \geqo \tfrac14 c' t_n^{2\theta+ \eta-3\kappa} \!\!\! \underset{n\to \infty}{\longrightarrow}\!  \infty
\end{equation}
because $2\theta+ \eta-3\kappa >0$ by (\ref{paramrange}). This combined with (\ref{case3}) completes the proof of (\ref{negativenewlim}). \cq 

\smallskip

We have proved (\ref{positivenewlim}) and (\ref{negativenewlim}), which imply Proposition \ref{p:main} and the proof of Theorem \ref{t:main} is now completed. \cqfd

\end{document}